\input amstex
\documentstyle{amsppt}
\let \aro=@
\NoBlackBoxes
\mag1100
\def\C{{\Bbb C}}
\def\R{{\Bbb R}}
\def\Z{{\Bbb Z}}

\def\N{{\Bbb N}}
\def\diff{\mathop{\sim}\limits}
\def\ed{{\ \diff_{eq}\ }}
\def\sc{\mathop{\#}\limits}
\font\nom=cmcsc10 scaled 1200
\newcount\numcount
\def\numerote{\global\advance\numcount by 1 \the\numcount}
\def\lastnum[#1]{{\advance \numcount by #1(\the\numcount)}}
\input psbox
\def \scaledpicture #1 by #2 (#3 scaled #4)
{\dimen0=#1 \dimen1=#2
\divide\dimen0 by 1000\multiply\dimen0 by #4
\divide\dimen1 by 1000\multiply\dimen1 by #4
$$\psboxto(\dimen0;\dimen1){#3.ps}$$}
\topmatter
\title Real quadrics in $\Bbb C^n$, complex manifolds and convex polytopes
\endtitle 
\rightheadtext{Real quadrics, complex manifolds and convex polytopes}
\author
Fr\'ed\'eric Bosio, Laurent Meersseman
\endauthor
\address 
{Fr\'ed\'eric Bosio}\hfill\hfill\linebreak 
\indent{Universit\'e de Poitiers}\hfill\hfill\linebreak 
\indent{UFR Sciences SP2MI}\hfill\hfill\linebreak 
\indent{D\'epartement de Math\'ematiques}\hfill\hfill\linebreak 
\indent{Boulevard Marie-et-Pierre-Curie}\hfill\hfill\linebreak
\indent{T\'el\'eport 2--BP 30179}\hfill\hfill\linebreak
\indent{86962 Futuroscope Chasseneuil Cedex}\hfill\hfill\linebreak
\indent{FRANCE}\hfill\hfill   
\endaddress 
\email Frederic.Bosio\@math.univ-poitiers.fr
\endemail 
\address 
{Laurent Meersseman}\hfill\hfill\linebreak
\indent{I.R.M.A.R.}\hfill\hfill\linebreak
\indent{Universit\'e de Rennes I}\hfill\hfill\linebreak
\indent{Campus de Beaulieu}\hfill\hfill\linebreak
\indent{35042 Rennes Cedex}\hfill\hfill\linebreak
\indent{FRANCE}\hfill\hfill   
\endaddress
\email laurent.meersseman\@univ-rennes1.fr \endemail
\date 19 March 2004 \enddate

  
\thanks 
The second author would like to thank Santiago L\'opez de Medrano for many fruitful discussions. Thanks also to Alberto Verjovsky
for giving the reference \cite{Al}.
\endthanks

\dedicatory
dedicated to Alberto Verjovsky on his $60^{\text{th}}$ birthday
\enddedicatory
\abstract
In this paper, we investigate the topology of a class of non-K\"ahler compact complex manifolds generalizing that of Hopf
and Calabi-Eckmann manifolds. These manifolds are diffeomorphic to
special systems of real quadrics in $\Bbb C^n$ which are invariant with respect to the natural action of the real torus $(\Bbb S^1)^n$
onto $\Bbb C^n$. The quotient space is a simple convex polytope. The problem reduces thus to the study of the topology of certain real algebraic sets and can be handled
using combinatorial results on convex 
polytopes. We prove that the homology groups of these compact complex manifolds can have arbitrary amount of torsion so that their
topology is extremely rich. We also resolve an associated wall-crossing problem by introducing holomorphic equivariant 
elementary surgeries related to some transformations of the simple convex polytope. Finally, as a nice consequence, we obtain that affine
non K\"ahler compact complex manifolds can have arbitrary amount of torsion in their homology groups, contrasting with the K\"ahler 
situation. 
\endabstract

\endtopmatter
\loadbold
\loadmsbm
\loadmsam
\nonfrenchspacing
\document
\numcount=0

\head
{\bf Introduction}
\endhead

This work explores the relationships existing between three classes of objects, coming from different domains of mathematics, namely:
\medskip

\noindent (i) {\it Real algebraic geometry}: the objects here are what we call {\it links}, that is transverse intersections in $\Bbb C^n$ of real quadrics of the
form
$$
\sum_{i=1}^n a_i\vert z_i\vert ^2=0
\eqno a_i\in \Bbb R
$$
with the unit euclidean sphere of $\Bbb C^n$.

\noindent (ii) {\it Convex geometry}: the class of simple convex polytopes.

\noindent (iii) {\it Complex geometry}: the class of non-K\"ahler compact complex manifolds of \cite{Me1}.
\medskip

The natural connection between these classes goes as follows. First, a link is invariant by the standard action of the real torus $(\Bbb S^1)^n$
onto $\Bbb C^n$ and the quotient space is easily seen to identify with a simple convex polytope (Lemma 0.11). 
Secondly, as a direct consequence of the construction of \cite{Me1},
each link (after taking the product by a circle in the odd-dimensional case) can be endowed with a complex structure of a manifold of \cite{Me1} (Theorem 12.2).
\smallskip

The aim of the paper is to describe the topology of the links and to apply the results to address the following question

\proclaim{Question}
How complicated can be the topology of the compact complex manifolds of \cite{Me1}?
\endproclaim

This program is achieved by making a reduction to combinatorics of simple convex polytopes: a simple convex polytope encodes completely the topology of the associated
link.
\smallskip

As shown by the question, the main motivation comes from complex geometry. Let us explain a little more why we find important to know the topology of the manifolds
of \cite{Me1}.

Complex geometry is concerned with the study of (compact) complex manifolds. Nevertheless, no general theory exists and only special
classes of complex manifolds as projective or K\"ahler manifolds or complex manifolds which are at least 
bimeromorphic to projective or K\"ahler ones are well understood. 
Moreover, except for the case of surfaces, there are few explicit examples having none of these properties; 
explicit meaning that it is possible to work with and to compute things on it. 
Indeed, the two classical families are the Hopf manifolds 
(diffeomorphic to $\Bbb S^1\times \Bbb S^{2n-1}$, see \cite{Ho}) and
the Calabi-Eckmann manifolds (diffeomorphic to $\Bbb S^{2p-1}\times \Bbb S^{2q-1}$, see \cite{C-E}). 

In \cite{LdM-Ve}, \cite{Me1} and \cite{Bo}, a new class of examples was provided. In particular, the class of \cite{Me1}
is explicit in the previous sense; the main complex geometrical properties (algebraic dimension,
generic holomorphic submanifolds, local deformation space, ...) of these objects are established in \cite{Me1}. 

Besides, it is proved
in \cite{Me2} that they
are small deformations of holomorphic principal bundles over projective toric varieties with fiber a compact complex torus. 
In this sense, they constitute a natural 
generalization of Hopf and Calabi-Eckmann manifolds, which can be deformed into compact complex manifolds fibering in elliptic
curves over the complex projective space $\Bbb P^{n-1}$ (Hopf case) or over the product of projective spaces
$\Bbb P^{p-1}\times\Bbb P^{q-1}$ (Calabi-Eckmann case).
One of the main interest of these manifolds however is that
they have a richer topology, since it is also proved in \cite{Me1} that complex structures on certain connected sums of
products of spheres can be obtained by this process. 

Nevertheless, these examples of connected sums constitute very particular cases of the construction and the problem of
describing the topology in the other cases was left wide open in \cite{Me1}. 
Of course, due to the lack of examples of non K\"ahler and non 
Mo\"{\i}shezon compact complex manifolds, the more intricate this topology is,
the more interesting is the class of \cite{Me1}. 
This is the starting point and motivation for this work and leads to the question stated above.

In \cite{Me1}, it was conjectured that they are all diffeomorphic to products of connected sums of spheres products
and odd-dimensional spheres. 

On the other hand, it follows from the construction that a manifold $N$ of \cite{Me1} is entirely characterized by a set $\Lambda$
of $m$ vectors of $\Bbb C^n$ (with $n>2m$). Moreover, a homotopy of $\Lambda$ in $\Bbb C^n$ gives rise to a deformation of $N$
as soon as an open condition is fullfilled at each step of the homotopy. 
If this condition is broken during the homotopy, the diffeomorphism type of the
new complex manifold $N'$ is different from that of $N$. In other words, there is a natural wall-crossing problem, and this leads to:

\proclaim{Problem}
Describe the topological and holomorphic changes occuring after a wall-crossing.
\endproclaim

This wall-crossing problem is linked with the previous question, since knowing how the topology changes after a wall-crossing,
one can expect describe the most complicated examples. But it has also a holomorphic part, since the initial and final manifolds are
complex.
 \smallskip

In this article, we address these questions and give a description as complete as we can of
the topology of these compact complex manifolds: 

\item{$\bullet$} Concerning the question above, the very surprising answer is that the topology of the complex manifolds of \cite{Me1}
is much more complicated than expected. Indeed, their homology groups can have arbitrary amounts of torsion (Theorem 14.1).
Counterexamples are given in Section 11, as well as a constructive way of obtaining these arbitrary amounts of torsion.
 
\item{$\bullet$} Concerning the wall-crossing problem, we show that crossing a wall means performing 
a complex surgery and describe precisely 
these surgeries from the topological and the holomorphic point of view (Theorems 5.4 and 13.3).
      
As an easy but nice consequence, we obtain that affine compact complex manifolds (that is compact complex manifolds with an affine
atlas) can have arbitrary amount of torsion. It is thus not possible to classify, up to diffeomorphism, affine compact complex 
manifolds or manifolds having a holomorphic affine connection in high dimensions ($\geq 3$).

It is interesting to compare this result with the K\"ahler case: it is known that affine K\"ahler manifolds are covered by complex 
tori (see \cite{K-W}), so the difference here is striking. Notice also that a statement similar to Theorem 14.1 is unknown for
K\"ahler manifolds.   

The paper is organized as follows. In Section 0, we collect the basic facts about the links. In particular, we introduce the simple convex polytope
associated to a link as well as a subspace arrangement whose complement has the same homotopy type as the associated link. We also recall
the previously known cases studied in \cite{LdM1} and \cite{LdM2}.

In part I, we prove that the classes of links up to equivariant diffeomorphism (equivariant with respect to the action of the real torus)
and up to product by circles are in $1:1$ correspondence with the combinatorial classes of simple convex polytopes (Rigidity Theorem
4.1). This is the first main result of this part. 
It allows us to translate topological problems about the links entirely in the world of combinatorics of simple convex polytopes. 
In particular, we recall the notion of 
flips of simple polytopes of \cite{McM} and \cite{Ti} in Section 2 and prove some auxiliary results. We define in
Section 3 a set of equivariant elementary surgeries on the links and prove in Section 4 (Theorem 4.7) that performing a flip
on a simple convex polytope means performing an equivariant surgery on the associated link. Finally, we introduce in Section 5 the 
notion of wall-crossing of links and prove the second main Theorem of this part (Wall-crossing Theorem 5.4): crossing a wall for a link
is equivalent to performing a flip for the associated simple convex polytope and therefore the wall-crossing can be described in terms
of elementary surgeries. As a consequence, we generalize a result of Mac Gavran (see \cite{McG}) and describe explicitely the
diffeomorphism type of certain families of links in Section 6.

In part II, we give a formula for computing the cohomology ring of a link in terms of subsets of the associated simple convex
polytope. To do this, we apply the Goresky-MacPherson formula \cite{G-McP} 
and the cohomology product formula of De Longueville \cite{DL} 
on the subspace arrangement mentioned earlier.
We rewrite them in terms of the simple polytope. The existence of such a formula is rather mysterious. Indeed it is somewhat miraculous
that Goresky-MacPherson and De Longueville formulas can be rewritten on the convex polytope and that they become so easy
in this new form. For example, it is rather difficult to check with the Goresky-MacPherson formula that the homology groups of a
link satisfy Poincar\'e duality; with this new formula, Poincar\'e duality is given by Alexander duality on the boundary of the 
simple convex polytope (seen as a sphere). The proof of this formula is long and technically difficult. It is a matter of taking 
explicit Alexander duals of cycles in simplicial spheres. The formula is stated up to sign (for the cohomology product) 
in Section 7 as Cohomology Theorem 7.6. and is
proved in Sections 7 and Sections 9 after some preliminary material about orientation and explicit Alexander duals in Section 8 and 9.
The sign is made precise in Section 10. Finally, applications and examples are given in Section 11, and it is proved that the
homology groups of a link can have arbitrary torsion (Torsion Theorem 11.11).

In part III, we apply the previous results to the family of compact complex manifolds of \cite{Me1}. In Section 12, we recall
very briefly their construction and prove that an even-dimensional link admits such a complex structure as well as the product of an
odd-dimensional link by a circle. We resolve the holomorphic wall-crossing problem in Section 13 (Theorem 13.3). 
Finally, in Section 14, we obtain as an easy consequence
of Theorem 11.11 that the homology groups of a compact complex manifold of \cite{Me1} can have arbitrary amount of torsion, and
as easy consequence of the construction that such a statement is true for affine compact complex manifolds.

Although the main motivation comes from complex geometry, part I (especially Section 6) should also be of interest for readers working on smooth actions of the torus on
manifolds. It can be seen as a continuation of \cite{LdM1}, \cite{LdM2} and \cite{McG}. On the other hand, the cohomology formula of Part II has its own interest as
a nice simplification of the Goresky-Mac Pherson and De Longueville formulas for a special class of subspace arrangements.

Notice that the smooth manifolds that we call links appear (but with a different definition, in particular not as intersection of
quadrics) in the study of toric or quasitoric manifolds (see \cite{D-J} and \cite{B-P}). In a sense, some results of this paper are
complementary to that of \cite{D-J} and \cite{B-P}.

\head
{\bf 0. Preliminaries}
\endhead

In this Section, we give the basic definitions, notations and lemmas. 
Some of the results are stated and sometimes proved in \cite {Me1} or \cite{Me2},
but in different versions; in this case we give the original reference, but at the same time, we give at least some indication about
the proof to be self-contained.

In this paper, we denote by $\Bbb S^{2n-1}$ the unit euclidean sphere of $\Bbb C^n$, 
and by $\Bbb D^{2n}$ (respectively $\overline{\Bbb D^{2n}}$) the unit euclidean open (respectively closed) ball of $\Bbb C^n$.

\definition{Definition 0.1}
A {\it special real quadric in $\C^n$} is a set of points $z\in\C^n$ satisfying :
$$
\sum_{i=1}^n a_i \vert z_i\vert ^2=0
$$
for some fixed $n$-uple $(a_1,\hdots, a_n)$ in $\R^n$.
\enddefinition 

We are interested in the topology of the intersection of a finite (but arbitrary) number of special real quadrics in $\C^n$ with the euclidean
unit sphere. We call such an intersection the {\it link} of the system of special real quadrics.

Let $A\in M_{np}(\R)$, that is $A$ is a real matrix with $n$ columns and $p$ rows. We write $A$ as $(A_1,\hdots, A_n)$. 
To $A$, we may associate $p$ special real quadrics in
$\Bbb C^n$ and a link, which we denote by $X_A$. The corresponding system of equations, that is :
$$
\left\{ \eqalign{\sum_{i=1}^n A_i\cdot \vert z_i\vert ^2&=0 \cr
\sum_{i=1}^n \vert z_i\vert ^2&=1}
\right .
$$
will be denoted by $(S_A)$. 

Notice that we include the special case $p=0$. In this situation, $A=0$ is a matrix of $M_{n0}(\R)$ and $X_A$ is $\Bbb S^{2n-1}$.

\definition{Definition 0.2} 
Let $A\in M_{np}(\R)$. We say that $A$ is {\it admissible} if it gives rise to a link $X_A$ whose system $(S_A)$ is non degenerate at 
every point of $X_A$. We denote by $\Cal A$ the set of admissible matrices.
\enddefinition

In this paper, we restrict ourselves to the case where $A$ is admissible. A link is thus a smooth compact
manifold of dimension $2n-p-1$ without boundary. Moreover it has trivial normal bundle in $\C^n$, so is orientable.

We denote by $\Cal H(A)$ the convex hull of the vectors $A_1$,
..., $A_n$ in $\R^p$.

\proclaim{Lemma 0.3 (cf [Me2], Lemma 1.1)}
Let $A\in M_{np}(\R)$. Then $A$ is admissible if and only if it satisfies :

\noindent (i) The Siegel condition : $0\in\Cal H(A)$.

\noindent (ii) The weak hyperbolicity condition : $0\in\Cal H(A_i\quad\vert\quad i\in I)\Rightarrow \text{\rom {cardinal}} (I)>p$.

\endproclaim

\demo{Proof}

Clearly $X_A$ is non vacuous if and only if the Siegel condition is satisfied.

Let $z\in X_A$ and let 
$$
I_z=\{1\leq i\leq n \quad\vert\quad z_i\not =0\}=\{i_1,\hdots,i_q\}\ .
\tag\numerote
$$

The system $(S_A)$ is non degenerate at $z$ if and only if the matrix :
$$
\tilde A_z=
\pmatrix
A_{i_1} &\hdots & A_{i_q} \cr
1 &\hdots & 1
\endpmatrix
$$
has maximal rank, i.e. rank $p+1$.

Assume the weak hyperbolicity condition. As $z\in X_A$, we have $0\in\Cal H((A_i)_{i\in I_z})$.
By Carath\'eodory's Theorem (\cite{Gr}, p.15), there exists a subset $J=\{j_1,\hdots,j_{p+1}\}\subset I_z$ 
such that $0$ belongs to 
$\Cal H((A_i)_{i\in J})$. Moreover,  
$(A_{j_1},\hdots,A_{j_{p+1}})$ has rank $p$, otherwise, still by Carath\'eodory's Theorem, $0$ would
be in the convex hull of $p$ of these vectors, contradicting the weak hyperbolicity condition.

As a consequence of these two facts, the vector space of linear relations between $(A_{j_1},\hdots,A_{j_{p+1}})$
has dimension one and is generated by a solution with all coefficients nonnegative. Assume that
$\tilde A_z$ has rank strictly less than $p+1$. Then, there is a
non-trivial linear relation between $(A_{j_1},\hdots,A_{j_{p+1}})$ with the additional property that the sum of the 
coefficients of this relation is zero. Contradiction.

Conversely, assume that the weak hyperbolicity condition is not satisfied. For example, assume that $0$ belongs to $\Cal H(A_1,\hdots,
A_p)$ and let $r\in(\R^+)^p$ such that :
$$
\sum_{i=1}^p r_i\cdot A_i=0,\quad \sum_{i=1}^p r_i=1\ .
$$

Then $z=(\sqrt r_1,\hdots,\sqrt r_p, 0,\hdots, 0)$ belongs to $X_A$ and rank $\tilde A_z$ is at most $p$ so $A$ is not admissible.
$\square$
\enddemo

Note that the intersection $\Cal A\cap M_{np}(\R)$ is open in $M_{np}(\R)$.
\medskip

Let us describe some examples.

\example{Example 0.4}
{\bf Let} $\boldkey p \boldkey = \boldkey 1$. 
Then the $A_i$ are real numbers. The weak hyperbolicity condition implies that none of the $A_i$ is zero.
Let us say that $a$ of the $A_i$ are strictly positive whereas $b=n-a$ 
of the $A_i$ are strictly negative. The Siegel condition implies that $a$ and 
$b$ are strictly positive. There is just one special real quadric, which is the equation of a cone over a product of spheres $\Bbb S^{2a-1}
\times\Bbb S^{2b-1}$. As we take the intersection of this quadric with the unit sphere, we finally obtain that $X_A$ is diffeomorphic 
to $\Bbb S^{2a-1}\times\Bbb S^{2b-1}$.
\endexample
\medskip

\example{Example 0.5}
{\bf Let} $\boldkey p \boldkey =\boldkey 2$. Then the $A_i$ are points in the plane containing $0$ in their convex hull (Siegel condition). The weak hyperbolicity 
condition implies that $0$ is not on a segment joining two of the $A_i$. Here are two examples of admissible configurations.

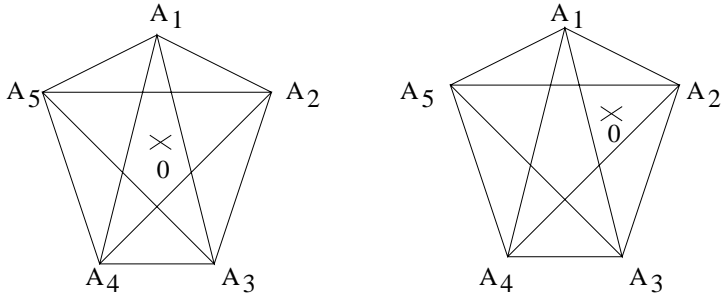
\hfil\scaledpicture 6.2in by 2.6in (fig1 scaled 550) \hfil
\vskip -.5cm

Assume that we perform a smooth homotopy  $(A^t)_{0\leq t\leq 1}$ between $A^0=A$ and $A^1$ in $\R^2$ such that $A^t$ {\it still
satisfies the Siegel and the weak hyperbolicity conditions for any } $t$. Then the union of the $X_{A^t}$ (seen as a smooth submanifold 
of $\C^n\times \R$) admits a submersion onto $[0,1]$ with compact fibers. Therefore, by Ehresmann's Lemma, this submersion is a locally
trivial fiber bundle and $X_{A^1}$ is diffeomorphic to $X_{A^0}=X_A$. Using this trick, it can be proven that $X_A$ is diffeomorphic
to $X_{A'}$, where $A'$ is a configuration of an odd number $k=2l+1$ of distinct points with weights $n_1,\hdots, n_k$ (see \cite{LdM2}).
The result of such an homotopy on the two configurations of the previous picture is represented below. The arrows indicate the
homotopy and the numbers appearing on the circles are the weights of the final configuration.

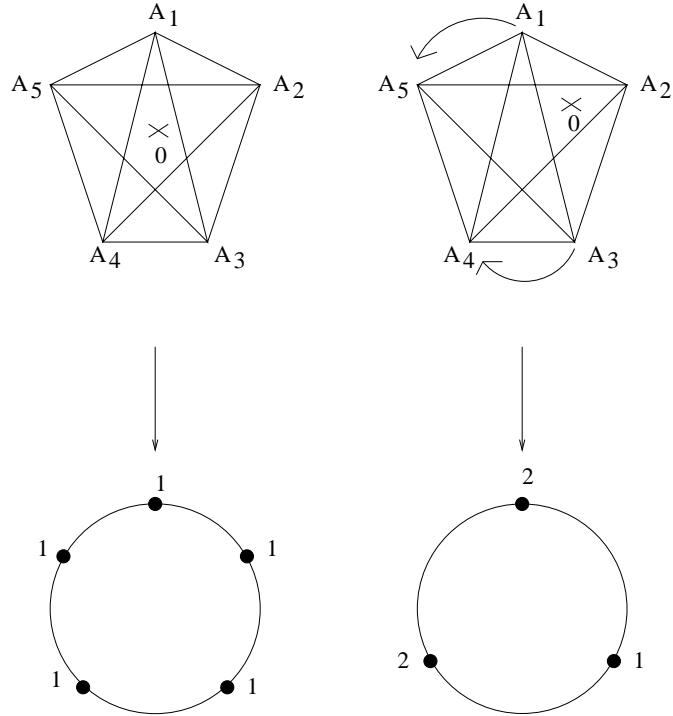
\hfil\scaledpicture 6.3in by 6.9in (fig2 scaled 500) \hfil
\null

\vskip -.5cm
These weights encode the topology of the links.

\proclaim{Theorem [LdM2]}
Let $p=2$ and let $A\in\Cal A$. Assume that $A$ is homotopic (in the sense given just above) to a reduced configuration of 
 $k=2l+1$ distinct points with weights $n_1,\hdots, n_k$. Then

\noindent (i) If $l=1$, then $X_A$ is diffeomorphic to $\Bbb S^{2n_1-1}\times \Bbb S^{2n_2-1}\times \Bbb S^{2n_3-1}$.

\noindent (ii) If $l>1$, then $X_A$ is diffeomorphic to
$$
\sc_{i=1}^k \Bbb S^{2d_i-1}\times \Bbb S^{2n-2d_i-2}
$$
where $\sc$ denotes the connected sum and where $d_i=n_i+\hdots+n_{i+l-1}$ (the indices are taken modulo $k$).
\endproclaim

In particular, $X_A$ is diffeomorphic to $\Bbb S^3\times\Bbb S^3\times\Bbb S^1$ 
for the configuration on the right of the previous figures, and diffeomorphic to $\sc (5) \Bbb S^3\times \Bbb S^4$ (that is the
connected sum of five copies of $\Bbb S^3\times \Bbb S^4$) for the configuration on the left.
\endexample
\medskip

\example{Example 0.6}
{\bf Products}. Let $A$ and $B$ be two admissible configurations of respective dimensions $(n,p)$ and $(n',p')$. Set
$$
C=\pmatrix
A &0 \\
-1 \hdots -1 & 1 \hdots 1 \\
0 & B
\endpmatrix
$$
Then it is straightforward to check that $C$ is admissible and that $X_C$ is diffeomorphic to the product $X_A\times X_B$. In
other words, {\it the class of links is stable by direct product}. In particular, the product of a link with an odd-dimensional sphere
is a link. For example, letting
$$
C=\pmatrix
A &0 \\
-1 \hdots -1 & 1
\endpmatrix
$$
then $X_C$ is diffeomorphic to $X_A\times\Bbb S^1$.
\endexample
\medskip

Let $\Cal L_A$ denote the complex coordinate subspace arrangement of $\C^n$ defined as follows :
$$
L_I=\{z\in\C^n \quad\vert\quad z_i=0 \text{ for } i\in I\}\in \Cal L_A \iff L_I\cap X_A =\emptyset
\tag\numerote
$$
and let $\Cal S_A$ be its complement in $\C^n$. In other words,
$$
\Cal S_A=\{z\in\C^n \quad\vert\quad 0\in\Cal H((A_i)_{i\in I_z})\}
\tag\numerote
$$
where $I_z$ is defined as in \lastnum[-2]. We have:

\proclaim{Lemma 0.7}
The sets $X_A$ and $\Cal S_A$ have the same homotopy type.
\endproclaim

\demo{Proof}
This is an argument of foliations and convexity already used in \cite{C-K-P}, \cite{LdM-V}, \cite{Me1} and \cite{Me2}. 
We sketch the proof and refer to these articles for more details.

Let $\Cal F$ be the smooth foliation of $\Cal S_A$ given by the action:
$$
(z,T)\in \Cal S_A\times\R^p\longmapsto \bigm (z_i\cdot \exp \langle A_i, T\rangle \bigm )_{i=1}^n\in \Cal S_A\ .
$$

Let $z\in \Cal S_A$ and let $F_z$ be the leaf passing through $z$. Consider now the map:
$$
f_z\ :\ w\in F_z\longmapsto \Vert w\Vert^2=\sum_{i=1}^n \vert w_i\vert ^2
$$

Using the strict convexity of the exponential map, it is easy to check that each critical point of $f_z$ is indeed a local minimum and
that $f_z$ cannot have two local minima and thus cannot have two critical points (see \cite{C-K-P} for more details). Now as $z\in
\Cal S_A$, then, by definition, $0$ is in the convex hull of $(A_i)_{i\in I_z}$. This implies that $F_z$ is a closed leaf and does not
accumulate onto $0\in\C^n$ (see \cite{Me1} and \cite{Me2}, Lemma 2.12 for more details). 
Therefore, the function $f_z$ has a global minimum, which is unique by the previous argument. Finally, a straightforward
computation shows that the minimum of $f_z$ is the point $w$ of $F_z$ such that:
$$
\sum_{i=1}^n A_i\vert w_i\vert ^2=0
$$

In particular $w/\Vert w\Vert$ belongs to $X_A$. 

As a consequence of all that, the foliation $\Cal F$ is trivial and its leaf space can be identified with 
$X_A\times \R^+_*\times\Bbb R^p$. More precisely, the map:
$$
\Phi_A\ :\ (z,T,r)\in X_A\times\R^p\times \R^+_*\longmapsto r\cdot \bigm (z_i\cdot \exp \langle A_i, T\rangle \bigm )_{i=1}^n\in \Cal S_A
\tag\numerote
$$
is a global diffeomorphism.
$\square$
\enddemo

Let $A\in\Cal A$. The real torus $(\Bbb S^1)^n$ acts on $\C^n$ by :
$$
(u,z)\in (\Bbb S^1)^n\times \C^n \longmapsto (u_1\cdot z_1,\hdots, u_n\cdot z_n)\in \C^n\ .
\tag\numerote
$$

Let $X$ be a subset of $\C^n$, which is invariant by the action \lastnum[0]. We define the {\it natural torus action on } $X$ as the restriction of \lastnum[0] to $X$.
In particular, every link $X_A$ for $A\in\Cal A$ is endowed with a natural torus action, as well as $\Bbb S^{2n-1}$, $\Bbb D^{2n}$ and $\overline{\Bbb D^{2n}}$.

\definition{Definition 0.8}
Let $A\in\Cal A$ and $B\in\Cal A$. We say that $X_A$ and $X_B$ are {\it equivariantly diffeomorphic} and we write $X_A\ed 
X_B$ if there
exists a diffeomorphism between $X_A$ and $X_B$ respecting the natural torus actions on $X_A$ and $X_B$.

More generally, we say that $X_A$ and $X_B\times (\Bbb S^1)^k$
are {\it equivariantly diffeomorphic} and we write $X_A\ed X_B\times (\Bbb S^1)^k$ if 
 there
exists a diffeomorphism between $X_A$ and $X_B\times (\Bbb S^1)^k$ respecting the natural torus actions on $X_A$ and 
on $X_B\times (\Bbb S^1)^k$ (seen as a subset of $\C^n\times \C^k$). 
\enddefinition

\proclaim{Lemma 0.9}
There exists $k\in \N$ and $B\in\Cal A$ such that $X_A$ is equivariantly diffeomorphic to $X_B\times (\Bbb S^1)^k$ and $X_B$
is 2-connected.
\endproclaim

\demo{Proof}
Assume that $X_A\cap\{z_1=0\}$ is vacuous. Let $A_i=\pmatrix a_i \cr \tilde A_i\endpmatrix$. As $A_1$ is not zero by weak hyperbolicity
condition, we may assume without loss of generality that $a_1\not =0$. Then, there exists an equivariant diffeomorphism :
$$
z\in X_A\longmapsto \left ( \dfrac{z_1}{\vert z_1\vert}, \dfrac{z_2}{\sqrt{1-\vert z_1\vert ^2}},\hdots, 
\dfrac{z_n}{\sqrt{1-\vert z_1\vert ^2}}\right )\in \Bbb S^1\times X_B
$$
where $B$ is defined as :
$$
B=\left ( \tilde A_2-\tilde A_1\dfrac{a_2}{a_1},\hdots,\tilde A_n-\tilde A_1\dfrac{a_n}{a_1}\right )\ .
$$

Now, $B$ is admissible since, at each point, the system $(S_B)$ has rank $p$. We may continue this process until we have
$X_A\ed X_B\times (\Bbb S^1)^k$ where the manifold $X_B\subset \C^{n-k}$ intersects each coordinate hyperplane of $\C^{n-k}$
(note that $X_B$ may be reduced to a point). This means that the subspace arrangement
$\Cal L_B$ has complex codimension at least $2$ in $\Bbb C^n$ and thus, by transversality, $\Cal S_B$ is 2-connected. By Lemma 0.7, 
this implies that $X_B$ is 2-connected. 
$\square$
\enddemo

We will denote by $\Cal A_0$ the set of admissible
matrices giving rise to a 2-connected link. More generally, let $k\in\N$. We will denote by $\Cal A_k$ the set of admissible matrices
giving rise to a link with fundamental group isomorphic to $\Z^k$. Of course, by Lemma 0.9, the set $\Cal A$ is the disjoint union
of all of the $\Cal A_k$ for $k\in\N$. Still from Lemma 0.9, observe that $k$ is exactly the number of coordinate hyperplanes of $\C^n$
lying in $\Cal L_A$.
\medskip

The action \lastnum[0] induces the following action of $\Bbb S^1$ onto a link $X_A$:
$$
(u,z)\in \Bbb S^1\times X_A \longmapsto u\cdot z \in X_A
\tag\numerote
$$

We call this action the {\it diagonal action} of $\Bbb S^1$ onto $X_A$. We have

\proclaim{Lemma 0.10}
Let $A\in\Cal A$. Then the Euler characteristic of $X_A$ is zero.
\endproclaim

\demo{Proof}
The diagonal action is the restriction to $X_A$ of a free action of $\Bbb S^1$ onto $\Bbb S^{2n-1}$, so is free. Therefore, we may construct a smooth non vanishing vector
field on $X_A$ from a constant unit vector field on $\Bbb S^1$.
$\square$
\enddemo
  
The quotient space of $X_A$ by the natural torus action is given by the positive solutions of the system
$$
A\cdot r=0 \qquad \sum_{i=1}^n r_i=1
\tag\numerote
$$
By the weak hyperbolicity condition, it has maximal rank. We may thus parametrize its set of solutions by
$$
r_i=\langle v_i, p\rangle+\epsilon_i \qquad p\in\Bbb R^{n-p-1}
\tag\numerote
$$
for some $v_i\in\Bbb R^{n-2p-1}$ and some $\epsilon_i\in\Bbb R$. Projecting onto $\Bbb R^{n-p-1}$, this gives an identification
of the quotient of $X_A$ by \lastnum[-3] as
$$
K_A=\{u\in\Bbb R^{n-p-1} \quad\vert\quad \langle v_i,u\rangle \geq -\epsilon_i\}
\tag\numerote
$$

\proclaim{Lemma 0.11}
Let $A\in\Cal A_k$.
The set $K_A$ is a (full) simple convex polytope of dimension $n-p-1$ with $n-k$ facets.
\endproclaim
\demo{Proof}
As $K_A$ is the quotient space of the compact manifold $X_A$ by the action of a compact torus, it is a compact subset
of $\Bbb R^{n-p-1}$.

Using \lastnum[0], $K_A$ is a bounded intersection of half-spaces, i.e. a (full) convex polytope of dimension $n-p-1$.

For every subset $I$ of $\{1,\hdots, n\}$, let :
$$
Z_I=\{z\in\C^n \quad\vert\quad z_i=0\text{ if } i\in I,\ z_i\not =0 \text{ otherwise }\}
\tag\numerote
$$

Let $z\in X_A$ and define $I_z$ as in \lastnum[-9]. Then, for every $z'$ belonging to the orbit of $z$, we have $I_z=I_{z'}$ and thus
the action respects each set $Z_{I_z}$. Moreover, the action induces a trivial foliation of $X_A\cap Z_{I_z}$.

It follows from all this that each $k$-face of $K_A$ corresponds to a set of orbits of points $z$ with fixed $I_z$, i.e. to a set 
$X_A\cap Z_{I_z}$. In particular, there is a numbering of the faces of $K_A$ such that each $j$-face is numbered by the $(n-p-1-j)$-uple
$I$ of the corresponding $Z_I$. As a first consequence, the number of facets of $K_A$ is exactly equal to the number of coordinate
hyperplanes of $\C^n$ whose intersection with $X_A$ is non vacuous, that is is equal to $n-k$ (see the remark just after the proof of
Lemma 0.9). As a second consequence of this numbering, each vertex $v$ corresponds to a $(n-p-1)$-uple $I$ and each facet having 
$v$ as vertex corresponds to a singleton of $I$: each vertex is thus 
attached to exactly $n-p-1$ facets, i.e. the convex
polytope is simple.
$\square$
\enddemo

We will call the set $K_A$ the {\it associate polytope} of $X_A$. We will denote by $P_A$ the combinatorial type of $K_A$ and
by $P^*_A$ the dual of $P_A$, which is thus the combinatorial type of a simplicial polytope.

Following the numbering introduced in the proof of the previous Lemma, we will see $P_A$ as a poset whose elements are subsets of
$\{1,\hdots,n\}$ satisfying :
$$
I\in P_A \iff L_I\cap X_A \not = \emptyset \iff Z_I\subset \Cal S_A\iff 0\in\Cal H((A_i)_{i\in I^c})
\tag\numerote
$$
where $I^c=\{1,\hdots, n\}\setminus I$.
We equip $P_A$ with the order coming from the inclusion
of faces. Of course $P^*_A$ will be seen as the same set but with the reversed order.

Let $(v_1,\hdots,v_n)$ be a set of vectors of some $\Bbb R^q$. Following \cite{B-L}, we call Gale diagram of $(v_1,\hdots, v_n)$ a set of
points $(w_1,\hdots, w_n)$ in $\R ^{n-q-1}$ verifying for all $I\subset \{1,\hdots, n\}$ :
$$
0\in \text{Relint }(\Cal H(w_i)_{i\in I}) \iff
\Cal H(v_i)_{i\in I^c}\text{ is a face of } \Cal H(v_1,\hdots, v_n)
\tag\numerote
$$
where Relint denotes the relative interior of a set.

Now, consider $K_A$. 
Notice that we may assume that the $\epsilon_i$ are positive, taking as $(\epsilon_1,\hdots,\epsilon_n)$
a particular solution of \lastnum[-5]. 
Under this assumption, let $B_i=v_i/\epsilon_i$ for $i$ between $1$ and $n$. The convex hull of $(B_1, \hdots, B_n)$
is a realization of $P^*_A$. Using \lastnum[0] and the weak hyperbolicity condition, it is easy to prove the following result.

\proclaim{Lemma 0.12 (cf [Me1], Lemma VII.2)}
The set $(B_1,\hdots, B_n)$ is a Gale diagram of $(A_1,\hdots, A_n)$.
\endproclaim
Notice that two Gale diagrams of the same set are combinatorially equivalent.
We finish this part with a realization theorem.

\proclaim{Realization Theorem 0.13 (see [Me1], Theorem 14)}
Let $P$ be the combinatorial type of a simple convex polytope. 
Then, for every $k\in \Bbb N$ there exists $A(k)\in\Cal A_k$ such that $P_{A(k)}=P$. In particular,
every combinatorial type of simple convex polytope can be realized as the associate polytope of some 2-connected link.
\endproclaim

\demo{Proof}
Let $P$ be the combinatorial type of a simple polytope and let $P^*$ be its dual. Realize $P^*$ in $\R^q$ (with $q=\dim P^*$)
as the convex hull of its vertices $(v_1,\hdots, v_n)$. 

Let us start with $k=0$. By Lemma 0.12, it is sufficient to find $A(0)\in \Cal A_0$ such
that $P^*$ is a Gale diagram of $A(0)$.

This can be done by taking a Gale transform (\cite {Gr}, p.84) of $(v_1,\hdots,v_n)$, that is by taking the transpose of a basis
of the solutions of :
$$
\left\{\eqalign{\sum_{i=1}^n x_iv_i&=0 \cr
\sum_{i=1}^n x_i&=0}\right .
$$
We thus obtain $n$ vectors $(A_1,\hdots, A_n)$ in $\R^{n-q-1}$. Set $A(0)=(A_1,\hdots, A_n)$. We have now to check that $A(0)\in \Cal A_0$.
By an immediate computation, the Gale transform $(A_1,\hdots, A_n)$ satisfies the Siegel condition. Assume that $0$ belongs to
$\Cal H(A_i)_{i\in I}$ for some $I=\{i_1,\hdots, i_p\}$. Then $\Cal H(v_i)_{i\in I^c}$ is a face of $P^*$ of dimension less than
$n-p-2$ with $n-p$ vertices. This face cannot be simplicial. Contradiction. The weak hyperbolicity condition is fulfilled.

Finally, as $P^*=P^*_{A(0)}$ has $n$ vertices, the link $X_{A(0)}$ intersects each coordinate hyperplane of $\C^n$ so is 2-connected 
(see Lemma 0.7).

Now, using the construction detailed in Example 0.6, we can find $A(k)\in\Cal A_k$
for every $k$ such that $P_{A(k)}=P$.
$\square$
\enddemo

Note that, when $P^*$ is the $n$-simplex, the previous construction (for a 2-connected link) yields $p=0$ and the corresponding $X_A$ is the
standard sphere of $\C^{n-1}$.
\vfill\eject
\specialhead
{\nom Part I: Elementaries surgeries, flips and wall-crossing}
\endspecialhead
\bigskip

\head
{\bf 1. Submanifolds of $X_A$ given by a face of $P_A$}
\endhead

Let $A\in\Cal A$ and let $F$ be a proper face of $P_A$ numbered by $I$. Then, we may associate to $F$ and $A$ a link which we will denote by
$X_F$ (by a slight abuse of notation), smoothly embedded in $X_A$. To do this, just recall by \lastnum[-1] that 
$$
B=(A_j)_{j\in I^c}
$$
is admissible and thus gives rise to a link $X_B$ in $\C^{n-b}$ where $b$ is the cardinal of $I$. Now, $X_B$ is naturally embedded into $X_A$ as $X_F$ by defining:
$$
X_F=L_I\cap X_A
\tag\numerote
$$
where $L_I$ was defined in \lastnum[-11]. Moreover, the natural torus action of $(\Bbb S^1)^n$ onto $X_A$ gives by restriction to
$L_I$ the natural torus action of $(\Bbb S^1)^{n-b}$ onto $X_F\ed X_B$.

We have

\proclaim{Proposition 1.1}
Let $A\in\Cal A$ and let $F$ be a face of $P_A$ of codimension $b$. Then,

\noindent (i) $X_F$ is a smooth submanifold of codimension $2b$ of $X_A$ which is invariant under the natural torus action.

\noindent (ii) The quotient space of $X_F$ by the natural torus action is $F\subset K_A$.

\noindent (iii) $X_F$ has trivial invariant tubular neighborhood in $X_A$.
\endproclaim

\demo{Proof}
The points (i) and (ii) are direct consequences of the definition \lastnum[0] of $X_F$. Let us prove (iii).
 For $\epsilon>0$, define:
$$
L_I^{\epsilon}=\{z\in\C^n \quad\vert\quad \sum_{i\in I} \vert z_i \vert ^2< \epsilon\}\ .
$$
and
$$
W_F^{\epsilon}=X_A\cap L_I^{\epsilon}\ .
$$

For simplicity, assume that $I=\{1,\hdots,b\}$. Set $y_j=z_j$ for $1\leq j\leq b$ and $w_j=z_{b+j}$ for $1\leq j \leq n-b$. For $\epsilon>0$ sufficiently small, the map
$$
\pi \ : \ (y,w)\in W^{\epsilon}_F \longmapsto \dfrac{1}{\sqrt \epsilon}\cdot y \in \Bbb D^{2b}
$$
is a smooth submersion. Indeed, a straightforward computation shows that
the previous map is a submersion as soon as $W_F^{\epsilon}$ does not intersect any of the sets
$$
\{w_j=0 \quad\vert\quad b+j\in J\}
$$
for $J$ satisfying $F\cap F_J=\emptyset$ (cf the proof of Lemma 0.3). As this
submersion has compact fibers, it is a locally trivial fiber bundle by Ehresmann's Lemma. It is even a trivial
bundle, since $\Bbb D^{2b}$ is contractible. Notice now that the action of $(\Bbb S^1)^n$ onto $W^{\epsilon}_F$ 
can be decomposed into an action
of $(\Bbb S^1)^b$ leaving fixed the $y$-coordinates and an action of $(\Bbb S^1)^{n-b}$ leaving fixed the $w$-coordinates.
The fibers of the previous submersion are invariant with respect to the action of $(\Bbb S^1)^{n-b}$ whereas the disk $\Bbb D^{2b}$ 
is invariant
with respect to the action of $(\Bbb S^1)^b$. All this implies that $W^{\epsilon}_F$ is equivariantly 
diffeomorphic to $X_F\times \Bbb D^{2b}$ endowed with its natural torus action. 
$\square$
\enddemo

In the case where $F$ is a simplicial face, then we can identify precisely $X_F$.

\proclaim{Proposition 1.2}
Let $A\in\Cal A_0$. The following statements are equivalent:

\noindent (i) $X_A$ is equivariantly diffeomorphic to the unit euclidean sphere $\Bbb S^{2n-1}$ of $\C^n$ equipped with the action
induced by the standard action of $(\Bbb S^1)^n$ on $\C^n$.

\noindent (ii) $X_A$ is diffeomorphic to $\Bbb S^{2n-1}$.

\noindent (iii) $X_A$ has the homotopy type of $\Bbb S^{2n-1}$.

\noindent (iv) $P_A$ is the $(n-1)$-simplex.
\endproclaim

\demo{Proof}
When $p=0$, the link $X_A$ is the unit euclidean sphere $\Bbb S^{2n-1}$ of $\C^n$ and the natural torus action comes from
the standard action of $(\Bbb S^1)^n$ on $\C^n$. On the other hand, when $P_A$ is the $(n-1)$-simplex, we have $p=0$,
since the dimension of $P_A$ is $n-p-1$; in this way, we get an equivalence between (i) and (iv).

Of course, (i) implies (ii) and (ii) implies (iii). So assume now that $X_A$ is a homotopy sphere of dimension $2n-1$.
Recall that a polytope with $n$ vertices is $k$-{\it neighbourly} if its $k$-skeleton coincides with the $k$-skeleton of a $(n-1)$-simplex
(cf \cite{Gr}, Chapter 7). In particular, a $(n-1)$-simplex is $(n-2)$-neighbourly. We will use the following Lemma:

\proclaim{Lemma 1.3} 
Let $A\in \Cal A_0$. The link $X_A$ is $(2k)$-connected if and only if $P^*_A$ is the combinatorial type of a $(k-1)$-neighbourly 
polytope.
\endproclaim

\demo{Proof of Lemma 1.3}
Assume that $P^*_A$ is $(k-1)$-neighbourly. This means that every subset of $\{1,\hdots,n\}$ of cardinal less than $k$ numbers a face
of $P^*_A$. Using \lastnum[-11] and \lastnum[-2], 
this means that every coordinate subspace of $\Cal L_A$ has at least complex codimension $k+1$.
By transversality, this implies that $\Cal S_A$ is $(2k)$-connected and thus, by Lemma 0.7, the link $X_A$ is $(2k)$-connected.

Now, assume moreover that $P^*_A$ is {\it not} $k$-neighbourly. Then, there exists a coordinate subspace $L_I$ in $\Cal L_A$ of
codimension $k+1$. The unit sphere $\Bbb S^{2k+1}$ of the complementary coordinate subspace $L_{I^c}$ lies in $\Cal S_A$ and is not
null-homotopic in $\Cal S_A$. Therefore, $\Cal S_A$ and thus $X_A$ are not $(2k+1)$-connected.
$\square$
\enddemo

Applying this Lemma gives that $P^*_A$ is $(n-2)$-neighbourly. But its dimension being $n-p-1$, this implies that $p$ equals $0$ and
that it is the 
$(n-1)$-simplex. Therefore (iii) implies (iv).
$\square$
\enddemo

\proclaim{Corollary 1.4}
Let $A\in\Cal A$. Then $P_A$ is the $(n-p-1)$-simplex if and only if $X_A$ is equivariantly diffeomorphic to $\Bbb S^{2n-2p-1}\times
(\Bbb S^1)^p$.
\endproclaim

\demo{Proof}
Assume that $P_A$ is the $(n-p-1)$-simplex.
The polytope $P_A$ having $n-p$ facets, we know that $A\in\Cal A_p$. By Lemma 0.9, there exists $B\in\Cal A_0$ such that 
$X_A\ed X_B\times (\Bbb S^1)^p$. Now, this
implies that $P_B=P_A$, so that $P_B$ is the $(n-p-1)$-simplex. We conclude by Proposition 1.2. 

The converse is obvious by Proposition 1.2.
$\square$
\enddemo

\proclaim{Corollary 1.5}
Let $F$ be a simplicial face of $P_A$ of codimension $b$. Then $X_F$ is equivariantly diffeomorphic to $\Bbb S^{2n-2p-2b-1}
\times (\Bbb S^1)^p$.
\endproclaim

\head
{\bf 2. Flips of simple polytopes}
\endhead

We will make use of the notion of flips of simple polytopes. This Section is deeply inspired from \cite{Ti}, \S3 (see also \cite{McM}). 
The main difference is
that we only deal with combinatorial types of simple polytopes. Recall that two convex polytopes are combinatorially equivalent
if there exists a bijection between their posets of faces which respects the inclusion. Two combinatorially equivalent convex
polytopes are PL-homeomorphic and the classes of convex poytopes up to combinatorial equivalence coincide with the classes of 
convex polytopes up to PL-homeomorphism. In the sequel, 
{\it we make no distinction between a convex polytope and its combinatorial class.}
No confusion should arise from this abuse.

\definition{Definition 2.1}
Let $P$ and $Q$ be two simple polytopes of same dimension $q$. Let $W$ be a simple
polytope of dimension $q+1$. We say that $W$ is a {\it cobordism} between $P$ and $Q$ if $P$ and $Q$ are disjoint facets of $W$.

In addition, if $W\setminus (P\sqcup Q)$ contains no vertex, we say that $W$ is a trivial cobordism; if $W\setminus (P\sqcup Q)$ contains a unique vertex, 
we say that $W$ is an {\it elementary cobordism} between $P$ and
$Q$.
\enddefinition

In the next Section, we will relate this notion of cobordism of polytopes to the classical notion of cobordism of manifolds (here
of links) via the Realization Theorem 0.13. This will justify the terminology.

Notice that the existence of a trivial cobordism between $P$ and $Q$ implies $P=Q$; notice also that a cobordism of simple polytopes may be decomposed into
a finite number of elementary cobordisms. 

Now, let $W$ be an elementary cobordism between $P$ and $Q$ and let $v$ denote the unique vertex of 
$W\setminus (P\sqcup Q)$. An edge attached to $v$ has another vertex which may belong to $P$ or $Q$. Let us say that, among the $(q+1)$ edges attached to $v$,
then $a$ of them join $P$ and $b$ of them join $Q$.

\definition{Definition 2.2 (compare with [Ti], \S 3.1)}
We call {\it index of} $v$ or {\it index of the cobordism} 
the couple of integers $(a,b)$ such that $a$ (respectively $b$) denotes the number of edges of $W$ attached
to $v$ and joining $P$ (respectively $Q$).

Let $P$ and $Q$ be two simple polytopes of same dimension $q$. Assume that there exists an elementary cobordism $W$ between them and let
$(a,b)$ denote its index. Then we say that $Q$ is obtained from $P$ by performing on $P$
a {\it flip} of type $(a,b)$, or that $P$ undergoes a {\it flip} of type $(a,b)$.
\enddefinition

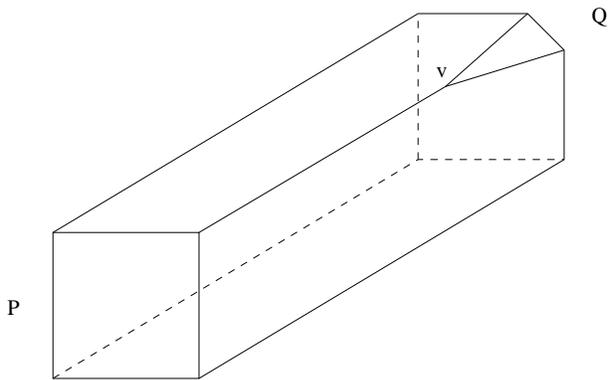
\hfil\scaledpicture 4.1in by 2.6in (flip1 scaled 700) \hfil

The previous picture is an example of a flip of type $(1,2)$.

Notice that if $Q$ is obtained from $P$ by a flip of type $(a,b)$, then obviously $P$ is obtained from $Q$ by a flip of type $(b,a)$.
Note also that we have the obvious relations $a+b=q+1$ and $1\leq a\leq q$ and $1\leq b\leq q$.

\proclaim{Lemma 2.3}
Every simple convex $q$-polytope can be obtained from the $q$-simplex by a finite number of flips.
\endproclaim

\demo{Proof}
Let $P$ be a simple convex $q$-polytope. Consider the product $P\times [0,1]$ and cut off one vertex of $P\times \{1\}$ by a
generic hyperplane.
The resulting polytope, let us call it $W$, is simple and realizes a cobordism between $P$ (seen as $P\times\{0\}$) and the $q$-simplex
(seen as the simplicial facet created by the cut). As observed above, this cobordism may be decomposed into a finite number of 
elementary cobordisms, that is of flips.
$\square$
\enddemo

Following \cite{Ti}, \S3.2, it is possible to give a more precise description of a flip of type $(a,b)$. 
We use the same notations as before. Let 
$F_1,\hdots, F_{q+1}$ be the facets of $W$ attached to the vertex $v$. As $W$ is simple, a sufficiently small neighborhood of $v$
in $W$ is PL-isomorphic to the neighborhood of a point in a $(q+1)$-simplex. As a consequence, each facet $F_i$ contains all the
edges attached to $v$ but one. Assume that $(F_1,\hdots F_b)$ contain all the edges joining $P$, whereas $(F_{b+1},\hdots, F_{q+1})$
contain all the edges joining $Q$.

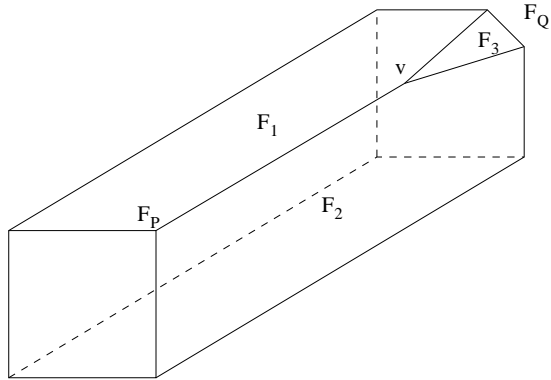
\hfil\scaledpicture 3.7in by 2.6in (flip1bis scaled 700) \hfil

Let $F_P=P\cap F_1\cap\hdots\cap F_b$ and $F_Q=Q\cap F_{b+1}\cap\hdots\cap F_{q+1}$. The face $F_1\cap\hdots\cap F_b$ (respectively
$F_{b+1}\cap\hdots\cap F_{q+1}$) is a pyramid with base $F_P$ (respectively $F_Q$) and apex $v$. As these faces are simple as convex 
polytopes, this implies that $F_P$ and $F_Q$ are simplicial. More precisely, if $a=1$ (respectively $b=1$), 
then $F_P$ (respectively $F_Q$) is
a point and $F_P\cap F_{q+1}=\emptyset$ (respectively $F_Q\cap F_{1}=\emptyset$). Otherwise $F_P$ is a simplicial face of
strictly positive dimension $q-b=a-1$ with facets $F_P\cap F_{b+1},\hdots,F_P\cap F_{q+1}$ (respectively
$F_Q$ is a simplicial face of strictly positive
dimension $b-1$ with facets
$F_Q\cap F_1,\hdots, F_Q\cap F_b$).

In the previous picture, $F_P$ is a point and $F_Q$ is a segment. There are three facets $F_1$, $F_2$, $F_3$ containing $v$.

The flip destroys the face $F_P$ and creates the face $F_Q$ in its place.
Continuously, the face $F_P$ is
homothetically reduced to a point and then this point is inflated to the face $F_Q$. In a more static way of thinking, a trivial
neighborhood of $F_P$ in $P$ is cut off and a closed trivial neighborhood of $F_Q$ in $Q$ is glued. In particular, the simple polytope obtained from $P$ by cutting off
a neighborhood of $F_P$ by a hyperplane and the polytope obtained from $Q$ by cutting off a neighborhood of $F_Q$ by a hyperplane 
are the same (up to combinatorial equivalence). Let us denote by $T$ this polytope. 

\definition{Definition 2.4}
The simple convex polytope $T$ will be called the {\it transition polytope} of the flip between $P$ and $Q$.
\enddefinition

\remark{Remark 2.5}
This definition is not the same as the definition of transition polytope of \cite{Ti}.
\endremark
\medskip

Notice that $T$ has just one extra facet (with respect to $P$ and $Q$), except for the special case of index $(1,1)$. Let us call it $F$.

The following picture describes a flip of type $(2,2)$. We simply drew the initial state $P$ and the final state $Q$ and indicated
the two edges $F_P$ of vertices $A$ and $B$ and $F_Q$ of vertices $A$ and $B'$.

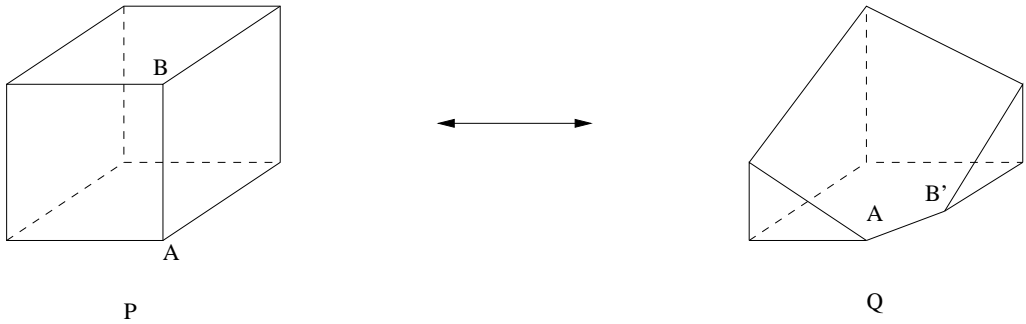
\hfil\scaledpicture 6.5in by 2.0in (flip2 scaled 750) \hfil
\vskip -.5cm

 To visualize the $4$-dimensional cobordism between $P$ and $Q$, just perform the following homotopy:
move the hyperplane supporting the upper facet of the cube to the bottom in order to 
contract the edge $AB$ to its lower vertex $A$; then move the hyperplane supporting the right facet of the cube
to the right in order
to inflate the transverse edge $AB'$, keeping $A$ fixed. The transition polytope $T$ is:

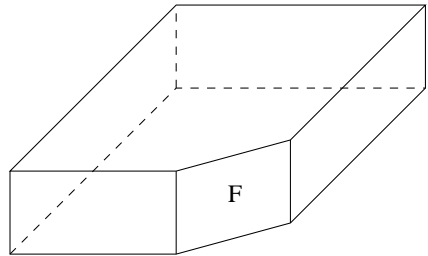
\hfil\scaledpicture 2.5in by 1.9in (flip3 scaled 800) \hfil
 \vskip -.5cm

\proclaim{Proposition 2.6}

\noindent (i) The extra facet $F$ of $T$ is combinatorially equivalent to $F_P\times F_Q$, that is to a product of a $(a-1)$-simplex by a
$(b-1)$-simplex.

\noindent (ii) A neighborhood of $F_P$ in $P$ (respectively $F_Q$ in $Q$) is combinatorially equivalent to $F_P\times \Cal C(F_Q)$
(respectively $\Cal (F_P)\times F_Q$), where $\Cal C(F_P)$ (respectively $\Cal C(F_Q)$) denotes the pyramid with base $F_P$
(respectively $F_Q$).
\endproclaim

\demo{Proof}
Assume that $P$ is a simplex. Cut off a neighborhood of $F_P$ by a hyperplane. The created facet is combinatorially equivalent
to a product of the simplex $F_P$ by a simplex $S$ of complementary dimension, whereas the cut part is combinatorially equivalent
to $F_P\times\Cal C(S)$, with the notation introduced in the statement of the Proposition.
Both statements follows then since 
the neighborhood of a simplicial face in a simple convex polytope is PL-homeomorphic to the neighborhood of
a face of same dimension in a simplex.
$\square$
\enddemo

In particular, the combinatorial types of $P$ and $Q$ can be recovered from that of $T$ (up to exchange of $P$ and $Q$): 
the face poset of $P$ is obtained
from that of $T$ by identifying two faces $A\times B$ and $A\times B'$ of $F_P\times F_Q$ and the face poset of $Q$ is obtained
from that of $T$ by identifying two faces $A\times B$ and $A'\times B$ of $F_P\times F_Q$.

Combining this observation with Proposition 2.6 yields

\proclaim{Corollary 2.7 (Rigidity of a flip)}
Let $Q$ and $Q'$ be obtained from $P$ by a flip of type $(a,b)$ along the same simplicial face $F_P$. Then $Q$ and $Q'$ are
combinatorially equivalent.
\endproclaim

Given a simple convex polytope $T$ with a facet $F$ combinatorially equivalent to a product of simplices $S_{a-1}\times S_{b-1}$,
we may define two posets from the poset of face of $T$ making the identifications explained just before Corollary 2.7. These
two posets {\it may or may not be} the face posets of some simple convex polytopes $P$ and $Q$ (see the examples below).
In the case they are, we write $P=F/S_{a-1}$ and $Q=F/S_{b-1}$. 
Of course, in the case of a flip, with the same notations as before, we have $P=T/F_P$ and $Q=T/F_Q$. The next 
Corollary is a reformulation of Corollary 2.7 which will be useful in the sequel.

\proclaim{Corollary 2.8}
Let $Q$ be obtained from $P$ by a flip along $F_P$ and let $T$ be the transition polytope. 
Let $P'$ and $Q'$ be two simple convex polytopes satisfying $P'=P/F_P$ and $Q'=Q/F_Q$. Then
$P$ and $P'$ are combinatorially equivalent as well as $Q$ and $Q'$.
\endproclaim

Let us describe another way of visualizing a flip. Let $P$ be a simple polytope and $F_P$ a simplicial face of dimension
$a-1$ of $P$. Let $Q$ be a simple polytope and assume that $Q$ is obtained from $P$ by performing a flip on $F_P$. 
Cut off $F_P$ by a hyperplane, you obtain the transition polytope $T$. Consider now a simplex $\Delta$ of same dimension
as $P$ and a $(a-1)$-face $F'$ of $\Delta$. Cut off $F'$ by a hyperplane, you obtain, with the notations of Proposition 2.6, the polytope
$\Cal F'\times S$, where $S$ is the maximal simplicial face of $\Delta$ without intersection with $F'$. It follows from Proposition 2.6
and Corollary 2.8 that the polytope $Q$ is combinatorially equivalent to the gluing of $T=P\setminus F_P\times \Cal C(S)$ and
of $\Delta\setminus F_P\times\Cal C(S)=\Cal (F')\times S$.

Finally, from all that preceeds, a complete combinatorial characterization of a flip may easily be derived. In the following statement, 
we consider also flips of type $(q+1,0)$, that is destruction of a $q$-simplex.

\proclaim{Proposition 2.9 ([Ti], Theorem 3.4.1)}
Let $Q$ be a simple polytope obtained from $P$ by a flip of type $(a,b)$. Using the same notations as before, we have

\noindent (i) If $a\not =1$, the facets $P\cap F_{b+1},\hdots, P\cap F_{q+1}$ undergo flips of index $(a-1,b)$.

\noindent (ii) The facets $P\cap F_1,\hdots, P\cap F_b$ undergo flips of index $(a,b-1)$.

\noindent (iii) The other facets keep the same combinatorial type.
\endproclaim

It is however important to remark that the notion of ``combinatorial flip'' is not well defined in the class of simple polytopes: the result of cutting off
a neighborhood of a simplicial face of a simple polytope and gluing in its place the neighborhood of another simplex may {\it not be} a 
convex polytope. Let us
give three examples of this crucial fact.

\example{Example 2.10}
Let $P$ be the $3$-simplex. Then, the result of cutting off an edge $AB$ 
and gluing in its place a transverse edge (that is the result of a ``combinatorial $2$-flip'')
is not the combinatorial type of a $3$-polytope.

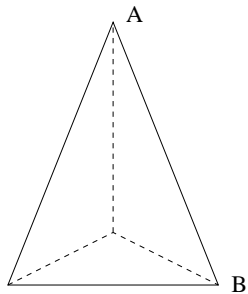
\hfil\scaledpicture 2.3in by 2.8in (fig3 scaled 500) \hfil
\endexample

\example{Example 2.11}
More generally, let $P$ be a simple convex polytope and $F_P$ a simplicial  face of dimension $q$, with $q>2$. Then, we cannot perform a flip along a strict face of
$F_P$.
\endexample

\example{Example 2.12}
Consider the following polytope (``hexagonal book'').

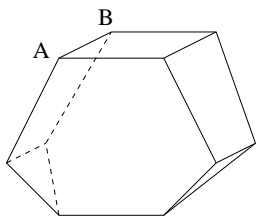
\hfil\scaledpicture 2.4in by 2.0in (fig4 scaled 500) \hfil
\vskip -.5cm

Then, the $2$-flip along the edge $AB$ does not exist.
\endexample

We finish with Section with the following result.

\proclaim{Proposition 2.13}
Let $P$ be a simple convex polytope and let $Q$ be obtained from $P$ by a flip of type $(a,b)$. Let $W$ be the elementary cobordism
between $P$ and $Q$. Assume that $P$ has $d$ facets. Then $W$ has $d+2$ facets if $a\not =1$ and $d+3$ facets if $a=1$.
\endproclaim

\demo{Proof}
In the special case where $a=b=1$, then $P=Q$ is the segment and $W$ is the pentagon.

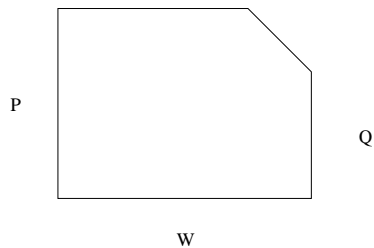
\hfil\scaledpicture 2.9in by 1.9in (fig5 scaled 600) \hfil
\vskip -.5cm

Thus $d$ is equal to $2$ and $W$ has $d+3$ facets.

Assume that $a$ and $b$ are different from one. Then $P$ and $Q$ have the same number $d$ of facets and there is a $1:1$ correspondance
between the facets of $P$ and that of $Q$ : according to Proposition 2.9, each facet of $P$ is transformed through a flip
(case (i) or (ii)) or just shifted (case (iii)) to a facet of $Q$. There are $d$ facets of $W$ which realize the previous trivial and
elementary cobordisms. Adding to this number $2$ to take account of $P$ and $Q$ gives that $W$ has $d+2$ facets.

Assume that $a=1$ and $b\not =1$. Then, as before, the $d$ facets of $P$ correspond to $d$ facets of $W$ realizing cobordisms with
$d$ facets of $Q$. But this time $Q$ has $d+1$ facets and this extra facet belongs to an extra facet of $W$ which does not intersect
$P$. Adding the two facets $P$ and $Q$ gives thus $d+3$ facets for $W$.

Finally, when $b=1$ and $a\not =1$, then the polytope $Q$ has $d-1$ facets; interverting the r\^ole of $P$ and $Q$ in the previous case
yields that $W$ has $(d-1)+3=d+2$ facets.
$\square$
\enddemo
\head
{\bf 3. Elementary surgeries}
\endhead

In this Section, we translate the notions of cobordisms and flips of simple polytopes at the level of the links.

We will make use several times of the following result:

\proclaim{Theorem of Extension of Equivariant Isotopies}
Let $M$ and $V$ be smooth compact manifolds endowed with a smooth torus action. Let $f\ : \ V\times [0,1]\to M$ be an equivariant isotopy. Then $f$ can be 
extended to an equivariant diffeotopy $F \ : \ M\times [0,1]\to M$ such that ${(F_t)}_{\vert V}\equiv f_t$ for $0\leq t \leq 1$.
\endproclaim

A proof of this fact {\it in the non equivariant case} can be found in \cite{Hi}, Chapter 8. Now, we may assume
that the diffeotopy extending an equivariant isotopy is also equivariant (see \cite{Br}, Chapter VI.3), so that this Theorem holds in the equivariant setting.

Let $A\in\Cal A$ and let $F$ be a {\it simplicial} face of $P_A$ of codimension $b$. As explained in Section 1, it
gives rise to an invariant submanifold $X_F$ of $X_A$ (see \lastnum[0]) with trivial invariant tubular neighborhood. 

By Corollary 1.5, as $F$ is simplicial of codimension $b$, then $X_F$ is equivariantly
diffeomorphic to $\Bbb S^{2a-1}\times (\Bbb S^1)^p$ (where $a=n-p-b$).

 But now, we can perform on $X_A$ an equivariant surgery as follows: choose a closed invariant tubular neighborhood
$$
\nu \ :\ X_F\times \overline{\Bbb D^{2b}} \longrightarrow \overline {W_F}
$$
where $W_F\subset X_A$ is an open (invariant) neighborhood of $X_F$. Then fix an equivariant identification
$$
\xi \ : \ \Bbb S^{2a-1}\times (\Bbb S^1)^p \longrightarrow X_F\ .
$$

Finally, set 
$$
\phi\equiv \nu\circ (\xi, \text{Id})\ : \ \Bbb S^{2a-1}\times (\Bbb S^1)^p \times \overline{\Bbb D^{2b}}\longrightarrow
\Bbb S^{2a-1}\times (\Bbb S^1)^p \times \overline{\Bbb D^{2b}}\ .
$$
We call $\phi$ a {\it standard product neighborhood} of $X_F$.

Then, remove $W_F$, and glue $\overline{\Bbb D^{2a}}\times (\Bbb S^1)^p\times \Bbb S^{2b-1}$ by $\phi$ along the boundary. 
We obtain thus a topological manifold $Y$. Since the natural torus actions on $\overline{\Bbb D^{2a}}\times (\Bbb S^1)^p\times \Bbb S^{2b-1}$ and on 
$\Bbb S^{2a-1}\times (\Bbb S^1)^p\times \overline{\Bbb D^{2b}}$ coincide on their common boundary, this topological manifold
supports a continuous action of
$(\Bbb S^1)^n$ which extends the natural torus action on $X_A\setminus W_F$. Using invariant collars for the boundary of $X_A\setminus W_F$ and for
the boundary of $\overline{\Bbb D^{2a}}\times (\Bbb S^1)^p\times \Bbb S^{2b-1}$, we may smooth $Y$ as well as the action in such a way that
the natural inclusions of $X_A\setminus W_F$ and $\overline{\Bbb D^{2a}}\times (\Bbb S^1)^p\times \Bbb S^{2b-1}$
in it are equivariant embeddings. As a consequence of the Theorem of Extension of Equivariant Isotopies, it can be proven that,
up to equivariant diffeomorphism, there are no other differentiable
structure and smooth action on $Y$ satisfying this property (see \cite{Hi}, Chapter 8 for the non equivariant case). 
The manifold $Y$ endowed with such a differentiable structure and such a smooth torus action is the result of 
our surgery.

Here is a combinatorial description of this surgery. Recall that $P_A$ identifies with the quotient of $X_A$ by the natural torus action. The neighborhood $W_F$ 
corresponds then to a neighborhood of $F$ in $P_A$. Consider now a simplex $\Delta$ of same dimension as $P_A$ and a face $F'$ of $\Delta$ of same dimension as $F$.
By Corollary 1.4, the link $X_{\Delta}$ corresponding to $\Delta$ is equivariantly diffeomorphic to $\Bbb S^{2n-2p-1}\times (\Bbb S^1)^p$ and a neighborhood $W_{F'}$
of $X_{F'}$ (coming from a neighborhood of $F'$ in
$\Delta$) is equivariantly diffeomorphic to $W_F$. The complement $X_{\Delta}\setminus W_{F'}$ is equivariantly diffeomorphic to
$$
(\Bbb S^{2n-2p-1}\setminus \Bbb S^{2a-1}\times \Bbb D^{2b})\times (\Bbb S^1)^p=\overline{\Bbb D^{2a}}\times \Bbb S^{2b-1} \times (\Bbb S^1)^p
$$
The surgery consists of removing $W_F$ in $X_A$ and $W_{F'}\ed W_{F}$ in $X_{\Delta}$ and of gluing the resulting manifolds along their boundary:
$$
X_A\setminus W_F \cup_{\psi} X_{\Delta}\setminus W_{F'}
\tag \numerote
$$ 
The map $\psi$ may be written as $\phi\circ (\phi')^{-1}$ for $\phi$ (respectively $\phi'$) a standard product neighborhood of $X_F$ in $X_A$ 
(respectively of $X_{F'}$ in $X_{\Delta}$).
 
We conclude from this description and from Corollary 2.8 that, at the level of the associate polytope, this 
surgery coincides exactly to a flip.

\definition{Definition 3.1}
Let $A\in\Cal A$. Let $(a,b)$ be a couple of positive integers satisfying $a+b=n-p$. Let $F$ be a 
simplicial face of $P_A$ of codimension $b$. We call {\it elementary surgery of type $(a,b)$ along $X_F$} the following equivariant
transformation of $X_A$: 
$$
(X_A\setminus \Bbb S^{2a-1}\times (\Bbb S^1)^p\times \Bbb D^{2b})\cup_{\phi} (\overline{\Bbb D^{2a}}\times (\Bbb S^1)^p\times \Bbb S^{2b-1})\ .
$$

Here $\Bbb S^{2a-1}\times (\Bbb S^1)^p\times \Bbb D^{2b}$ is embedded in $X_A$ by means of a standard product neighborhood $\phi$
and the gluing is made along the common boundary by the restriction of $\phi$ to this boundary.

In the particular case where $a=1$, we restrict the definition of elementary surgery to the case where $X_A$ is equivariantly
diffeomorphic to $X_B\times\Bbb S^1$ and where the surgery is made as follows
$$
(X_B\setminus (\Bbb S^1)^p\times \Bbb D^{2b})\times\Bbb S^1\cup_{\phi} ((\Bbb S^1)^p\times \Bbb S^{2b-1})\times\overline{\Bbb D^{2}}\ .
$$
\enddefinition

These surgeries depend {\it a priori} on the choice of $\phi$. But, in fact

\proclaim{Lemma 3.2}
The result of an elementary surgery is independent of the choice of $\phi$, that is: given two standard product neighborhoods $\phi$ and $\phi'$ ,
the manifolds 
$$
X_{\phi}=(X_A\setminus \Bbb S^{2a-1}\times (\Bbb S^1)^p\times \Bbb D^{2b})\cup_{\phi} (\overline{\Bbb D^{2a}}\times (\Bbb S^1)^p\times \Bbb S^{2b-1})\ .
$$
and
$$
X_{\phi'}=(X_A\setminus \Bbb S^{2a-1}\times (\Bbb S^1)^p\times \Bbb D^{2b})\cup_{\phi'} (\overline{\Bbb D^{2a}}\times (\Bbb S^1)^p\times \Bbb S^{2b-1})\ .
$$
are equivariantly diffeomorphic.
\endproclaim

\demo{Proof}
It is enough to prove that $\phi$ and $\phi'$ are equivariantly isotopic. As in the non equivariant case, the uniqueness of gluing
for isotopic diffeomorphisms is a direct consequence of the Theorem of Extension of Isotopies.

Now, any two invariant tubular neighborhoods of $X_F$ are equivariantly isotopic \cite{Br}, Chapter VI.2. Thus, we may assume that 
$$
\phi(\Bbb S^{2a-1}\times (\Bbb S^1)^p\times \overline{\Bbb D^{2b}})=\phi'(\Bbb S^{2a-1}\times (\Bbb S^1)^p\times \overline{\Bbb D^{2b}})
$$
and that the map $f=\phi'\circ \phi^{-1}$ is of the form
$$
(z,\exp{it},w) \in \Bbb S^{2a-1}\times (\Bbb S^1)^p \times \overline{\Bbb D^{2b}} \longmapsto (f_1(z,\exp{it}), f_2(z,\exp{it}),A(z,\exp{it})\cdot w)
$$
where $A$ is a smooth invariant map from $\Bbb S^{2a-1}\times (\Bbb S^1)^p$ to the group of matrices SO$_{2b}$. Moreover, the equivariance of $f$ implies that each
matrix $A(z,\exp{it})$ is of the form
$$
\pmatrix
\exp{i\theta_1} & &0 \\
&\ddots \\
0&  & \exp{i\theta_b}
\endpmatrix
$$

We may thus easily equivariantly isotope $f$ to
$$
(z,\exp{it},w) \in \Bbb S^{2a-1}\times (\Bbb S^1)^p \times \overline{\Bbb D^{2b}} \longmapsto (f_1(z,\exp{it}), f_2(z,\exp{it}),w)
$$
and it is enough to prove that the equivariant diffeomorphism $\tilde f=(f_1,f_2)$ of $\Bbb S^{2a-1}\times (\Bbb S^1)^p$ is equivariantly isotopic to the identity.

Still by equivariance, we have 
$$
\tilde f(z,\exp{it})=\exp{it}\cdot \tilde f(z,1)
$$
so we may equivariantly isotope $\tilde f$ to a map of the form
$$
(z,\exp{it})\in \Bbb S^{2a-1}\times (\Bbb S^1)^p \longmapsto (h(z),\exp{it})\in \Bbb S^{2a-1}\times (\Bbb S^1)^p
$$
where $h$ is an equivariant diffeomorphism of $\Bbb S^{2a-1}$. Finally, using Lemma 3.3 (stated and proved below), $h$ and thus $f$ are equivariantly isotopic
to the identity. This is enough to show the result.
$\square$
\enddemo

\proclaim{Lemma 3.3} 
Let $h$ be an equivariant diffeomorphism of the sphere $\Bbb S^{2a-1}$. Then $f$ is equivariantly isotopic to the identity.
\endproclaim

\demo{Proof}
We proceed by induction on $a$. For $a=1$, the map $h$ is a translation so the result is clear. Assume the result for some $a\geq 1$ and let $h$ be an equivariant 
diffeomorphism of $\Bbb S^{2a+1}$. 

By equivariance, the submanifold
$$
X=\{z\in \Bbb S^{2a+1} \quad\vert\quad z_{a+1}=0\}\ed \Bbb S^{2a-1}
$$
is invariant by $h$.

We shall construct two invariant tubular neighborhoods of $X$. First, consider, for $0<\epsilon<1$, 
$$
X_{\epsilon}=\{z\in \Bbb S^{2a+1} \quad\vert\quad \vert z_{a+1}\vert ^2 \leq \epsilon \}\ed \Bbb S^{2a-1}\times\overline{\Bbb D^2}
$$
and the equivariant bundle map
$$
z\in X_{\epsilon} \buildrel \xi \over \longmapsto \dfrac{1}{\sqrt{1-\vert z_{a+1} \vert ^2}}(z_1,\hdots, z_a,0)\in X
$$

Secondly, let $f$ be the restriction of $h^{-1}$ to $X$. Set $\tilde X_{\epsilon}=f^*X_{\epsilon}$ (pull-back bundle by $f$), and let $\tilde f$ denote the natural
map between $\tilde X_{\epsilon}$ and $X_{\epsilon}$. The map $h\circ \tilde f$ defines the second tubular neighborhood of $X$ in $\Bbb S^{2a+1}$.

By \cite{Br}, Chapter VI.3, there exists an equivariant isotopy of tubular neighborhoods
$$
H\: \ X_{\epsilon}\times [0,1] \longrightarrow \Bbb S^{2a+1}
$$
with $H_0\equiv \text{Id}$ and $H_1(X_{\epsilon})\equiv h\circ\tilde f (\tilde X_{\epsilon})\equiv h(X_{\epsilon})$. In particular, $H_1$ differs from $h$ by an
equivalence of equivariant bundles
$$
\CD
X_{\epsilon} \aro > h^{-1}\circ H_1 >> X_{\epsilon} \\
\aro V \xi VV \aro V \xi VV \\
X \aro >f >> X
\endCD
$$

Since $X\ed \Bbb S^{2a-1}$, by induction, the map $f$ is equivariantly isotopic to the identity and it is easy to lift this isotopy to an isotopy $G$ between
$H_1$ and $h$.

Combining $H$ and $G$, we obtain an equivariant isotopy
$$
F\ : \ [0,1]\times X_{\epsilon} \longrightarrow \Bbb S^{2a+1}
$$
such that $F_0$ is the natural inclusion map and $F_1\equiv h_{\vert X_{\epsilon}}$.

By the Theorem of Extension of Equivariant Isotopies, $F$ extends to an equivariant diffeotopy between
some map $g$ with $g_{\vert X_{\epsilon}}\equiv h$ and the identity. As this construction can be achieved for any choice of $0<\epsilon <1$, we may assume that 
$g\equiv h$ on the whole sphere.
$\square$
\enddemo

We note that the result of such a surgery {\it may or may not} be a link. Indeed, in Examples 2.10, 2.11 and 2.12, we may perform 
elementary surgeries but the quotient space of the new manifold by the action of the real torus cannot be identified with a
simple polytope, therefore the new manifold is not a link.

Consider now the following more subtle case. Let $X_A$ be a link and let $Q$ be the simple convex polytope obtained from $P_A$ by
performing a flip of type $(a,b)$ along some simplicial face $F$. Then, call $Y$ the manifold obtained from $X_A$ by performing
an elementary surgery of type $(a,b)$ along $X_F$. As the surgery is equivariant, the manifold $Y$ is endowed with a smooth action of
the real torus on it. It follows from Corollary 2.8 that the quotient space of $Y$ by this
action can be identified with $Q$. This means that this quotient space is in bijection with $Q$, that the orbit over a point in the
interior of $Q$ is $(\Bbb S^1)^n$, whereas the orbit over a point in the interior of a facet of $Q$ is $(\Bbb S^1)^{n-1}$ and so on.
 We still call associate polytope the resulting polytope.
Finally, each closed face of $Q$ corresponds to an invariant submanifold of $Y$ with trivial invariant tubular neighborhood. 
In fact, every such face $S$ is obtained from a face
$R$ of $P_A$ by a certain flip, as precised in Proposition 2.9. The corresponding invariant submanifold $Y_S$ is thus obtained from $X_R$ by performing the corresponding
elementary surgery. More precisely,
write
$$
Y=(X_A\setminus W_F)\cup_{\psi} (X_{\Delta}\setminus W_{F'})
$$
as in \lastnum[0], then we have
$$
Y_S=(X_R\setminus W_F\cap X_R)\cup_{\psi} (X_{R'}\setminus W_{F'}\cap X_{R'})
$$
for some well-chosen face $R'$ of $\Delta$. Let
$$
\nu \ : \ X_R \times \Bbb D^{2b'}\longrightarrow W_R\subset X_A
$$
be a trivial invariant tubular neighborhood of $X_R$ (we denote the codimension of $X_R$ in $X_A$ by $b'$). We assume that $W_R$ is small enough to have 
$$
\nu^{-1}(W_R\cap W_F)=(X_R\cap W_F) \times \Bbb D^{2b'}
$$
Then the composition
$$
(X_{R'}\cap W_{F'})\times \Bbb D^{2b'} \buildrel (\psi, \text{Id}) \over \longmapsto (X_R\cap W_F)\times \Bbb D^{2b'} \buildrel \psi^{-1} \over
\longmapsto W_{F'}
$$
can be extended to a (trivial) invariant tubular neighborhood 
$$
\nu'\ : \ X_{R'} \times \Bbb D^{2b'}\longrightarrow W_{R'}\subset X_{\Delta}
$$
since $\psi^{-1}\circ\nu$ maps $X_R\cap W_F$ onto $X_{R'}\cap W_{F'}$. Finally, set $\nu_S\equiv \nu\cup_{\psi}\nu'$. Then $\nu_S$ maps
$$
(X_R\setminus W_F)\times \Bbb D^{2b'} \cup_{(\psi, \text{Id})} (X_{R'}\setminus W_{F'})\times \Bbb D^{2b'}=Y_S\times \Bbb D^{2b'}
$$
to $W_R\setminus W_F \cup_{\psi} W_{R'} \setminus W_{F'}$, that is, $\nu_S$ is a trivial invariant
tubular neighborhood of $Y_S$. 

Assume that $Y_S$ is equivariantly diffeomorphic to some 
$\Bbb S^{2a'-1}\times (\Bbb S^1)^{p'}$. Then we may perform an elementary surgery corresponding to this choice of $Y_S$. In particular, we may perform
an elementary surgery corresponding to any choice of a flip of $Q$, as soon as the corresponding invariant submanifold of $Y$ is equivariantly diffeomorphic to some 
$\Bbb S^{2a'-1}\times (\Bbb S^1)^{p'}$. In this case, we say that the flip is {\it good}. 

We may then repeat this process and construct manifolds obtained from a link by a finite number of elementary surgeries corresponding to good flips of the associate
polytope.

Nevertheless, 
{\it it is not clear a priori that $Y$ as well as the manifolds obtained from $Y$ are equivariantly diffeomorphic to a link, 
that is to a transverse intersection of special real quadrics.}

\definition{Definition 3.4}
We call {\it pseudolink} a manifold obtained from a link by a finite number of elementary surgeries corresponding to good flips of the
associate polytopes.
\enddefinition

We will see now that every flip is good.

\proclaim{Proposition 3.5}
Let $X$ be a pseudolink such that its associate polytope $P$ is a $d$-simplex. Then $X$ is, up to product by circles, equivariantly 
diffeomorphic to the unit euclidean sphere $\Bbb S^{2d+1}$ of $\C^{d+1}$ endowed with the natural action of $(\Bbb S^1)^{d+1}$ on it.
\endproclaim

\demo{Proof}
The proof is by induction on $d$. If $d=0$, then $X$ is obviously a product of circles and the Proposition is satisfied.

Assume now that the Proposition is true for simplices of dimension at most $d$ and consider $X$ a pseudolink whose associate polytope
$P$ is a $(d+1)$-simplex. Then $P$ can be seen as a pyramid with base a $d$-simplex $P'$ and can be decomposed into a closed neighborhood
of $P'$ glued along the common boundary with a closed neighborhood of a $0$-simplex $v$ (a point). This means that
$X$ is equivariantly diffeomorphic to the gluing of an invariant closed neighborhood of $X_P'$ with an invariant closed neighborhood of $X_v$ by the identity along the 
common boundary. We may assume that these neighborhoods are tubular and thus trivial. Using the induction hypothesis and standard product neighborhoods,
we may write
$$
X\ed \Bbb S^{2d+1}\times (\Bbb S^1)^p\times \overline{\Bbb D^2} \cup_{\phi} \overline{\Bbb D^{2b}}\times (\Bbb S^1)^p\times \Bbb S^1
$$
for some $p\geq 0$ and some equivariant diffeomorphism $\phi$ of $\Bbb S^{2d+1}\times (\Bbb S^1)^{p+1}$.
 Using Lemma 3.3, we may
assume that $\phi$ is the identity. Therefore, 
$X$  is, up to product by circles, equivariantly 
diffeomorphic to the unit euclidean sphere $\Bbb S^{2d+3}$ of $\C^{d+2}$ endowed with the natural action of $(\Bbb S^1)^{d+2}$ on it.
$\square$
\enddemo

\proclaim{Corollary 3.6}
Every flip of the associate polytope of a pseudolink is good.
\endproclaim

We finish this Section with a Proposition which will be useful in the sequel.

\proclaim{Proposition 3.7}
Let $A\in\Cal A_k$ and $B\in\Cal A_l$. Assume that $X_B$ is obtained from $X_A$ by performing an elementary surgery of type $(a,b)$
corresponding to a flip. Then,

\noindent (i) If $1<a<n$ or $a=b=1$, then $k=l$.

\noindent (ii) If $a=1$ and $b\not =1$, then $k=l+1$.

\noindent (iii) If $a=n$ and $a\not = 1$, then $k=l-1$.
\endproclaim

\demo{Proof}
As the links $X_A$ and $X_B$ have same dimension, as well as $P_A$ and $P_B$, the numbers $n$ and $p$ are the same for both links.
This implies that $k$ (respectively $l$) is equal to $n$ minus the number of facets of $P_A$ (respectively $P_B$) (see Lemma 0.11).
Now, the results follow easily from the fact that a flip of type $(a,b)$ does not create nor destroy any facet if $1<a<n$ or 
$a=b=1$ (see the figure in the proof of Proposition 2.13),
creates a facet if $a=1$ and $b\not =1$ and destroys a facet if $a=n$ and $a\not = 1$ (see Proposition 2.9).
$\square$
\enddemo

\head
{\bf 4. The Rigidity Theorem}
\endhead

We are now in position to prove:

\proclaim{Rigidity Theorem 4.1}

\noindent (i) Every pseudolink is a link.

\noindent (ii) Let $A\in \Cal A_k$ and $B\in\Cal A_k$ for some $k$. Then $X_A\ed X_B$ if and only if $P_A=P_B$.
\endproclaim

\remark{Remark 4.2}
Let $F_0$ denote the product of complex projective lines $\Bbb P^1\times\Bbb P^1$ and let $F_1$ denote the Hirzebruch surface
obtained by adding a section at the infinite to the line bundle of Chern class $1$ over $\Bbb P^1$. Both are projective toric varieties
and thus admit a smooth, hamiltonian action of $(\Bbb S^1)^2$ with quotient space a convex polygon. 
In both cases, the polygon is a $4$-gon
(see \cite{Fu}), so the two  quotient spaces are combinatorially equivalent as convex polygons. Nevertheless, the two manifolds are
not even topologically the same (see \cite{M-K}): $F_0$ is diffeomorphic to a product $\Bbb S^2\times\Bbb S^2$, whereas $F_1$ is the
only non-trivial $\Bbb S^2$-bundle over $\Bbb S^2$. This example shows that the Rigidity result stated above is not obvious at all 
and is very particular to our situation.
\endremark
\medskip

\remark{Remark 4.3}
Let $p=0$ and $n\geq 2$. Then, $X_A$ is the unit euclidean sphere $\Bbb S^{2n-1}$ of $\C^n$. We may perform an equivariant surgery as follows:
$$
\eqalign{
(&X_A\setminus \Bbb S^{1}\times \Bbb D^{2n-2})\cup (\overline{\Bbb D^{2}}\times \Bbb S^{2n-3})\cr
=&(\overline{\Bbb D^{2}}\times \Bbb S^{2n-3})\cup (\overline{\Bbb D^{2}}\times \Bbb S^{2n-3})=\Bbb S^2\times\Bbb S^{2n-3}
}
$$
This surgery looks like an elementary surgery of type $(1,n)$. In particular, it is easy to check that the quotient space of
$\Bbb S^2\times\Bbb S^{2n-1}$ by the induced torus action can be identified with 
the prism with base a $(n-2)$-simplex, that is the simple convex
polytope obtained from the $(n-1)$-simplex $P_A$ by a flip of type $(1,n)$. Nevertheless, this is not an elementary surgery
by Definition 3.1 ($X_A$ is simply-connected) and the resulting manifold is not a link by Rigidity Theorem 4.1 but a
quotient of a link by an action of $\Bbb S^1$. The simply-connected link corresponding to the prism with base a $(n-2)$-simplex is
$$
\eqalign{
&(\Bbb S^{2n-1}\setminus \Bbb S^1\times \Bbb D^{2n-2})\times\Bbb S^1\cup (\Bbb S^1\times \Bbb S^{2n-3})\times\overline{\Bbb D^{2}}\cr
=&(\overline{\Bbb D^{2}}\times \Bbb S^1)\times\Bbb S^{2n-3}\cup (\Bbb S^1\times \overline{\Bbb D^{2}})\times\Bbb S^{2n-3}
=\Bbb S^3\times\Bbb S^{2n-3}
}
$$
\endremark
\medskip

\demo{Proof}
Let $P$ be a convex simple polytope. Call {\it length} of $P$ the minimal number of flips necessary to pass from the simplex 
(of same dimension as $P$) to $P$. This number exists by Lemma 2.3.

The proof is by induction on the length of the associate polytope. More precisely, the induction hypothesis (at order $l$) 
is that statements (i) and (ii) are true for links and pseudolinks with associate polytopes of length less than or equal to $l$.
This hypothesis is satisfied at order $0$ by Propositions 1.2 and 3.5.

Assume the hypothesis at order $l$, and consider $X$ a pseudolink with associate polytope $P$ of length $l+1$. Then, if $P$
undergoes some well-chosen flip, we obtain a simple convex polytope $Q$ with length $l$. As usually, let $(a,b)$ denote the type of
the flip and $F$ the simplicial face along which the flip is made. Remark that this implies that $P$ is obtained from $Q$ by performing a flip of type $(b,a)$ along some simplicial face $F'$. Perform an elementary surgery of type $(a,b)$ along the submanifold
of $X$ corresponding to $F$. We recover a pseudolink $Y$ whose associate polytope is $Q$. 
By induction, $Y$ is a link $X_A$ for
$A$ belonging to some $\Cal A_k$. 

Define $k'$ as $k$ if $1<a<n$ or $a=b=1$, as $k+1$ if $a=1$ and $b\not = 1$, and as $k-1$ otherwise. In this last case,
notice that $k-1$ is positive: $X$ is obtained from $X_A$ by an elementary surgery of type $(1,n)$, so, by Definition 3.1, the link
$X_A$ is not simply-connected. By Realization Theorem 0.13, there exists $B\in \Cal A_{k'}$ such that
$P_B$ is combinatorially equivalent to $P$. Perform an elementary surgery of type $(a,b)$ along the submanifold
of $X_B$ corresponding to $F$. By induction, the result of this surgery is a link $X_{A'}$. Due to the choice of $k'$, we have
$A'\in \Cal A_k$ by Proposition 3.7. Therefore, the second statement of the induction hypothesis implies that $X_{A'}\ed X_A$.

The conclusion of what preceeds is that both $X_B$ and $X$ are obtained from the same link $X_{A'}\ed X_A$ by performing an elementary
surgery of type $(b,a)$ along the same invariant submanifold (the submanifold corresponding to $F'$ in $Q$). Therefore, $X_B$ and $X$
are equivariantly diffeomorphic and $X$ is a link. This proves the first statement for associate polytopes of length $l+1$. 
Moreover, if you consider now any link $X_C$ with
$P_C=P$ and $C\in \Cal A_{k'}$, then the same proof implies that $X_B\ed X_C$. As these considerations do not depend on the value of
$k'$, this proves one implication of statement (ii). But the converse is easy: two equivariantly diffeomorphic links have the
same combinatorics of orbits, that is have combinatorially equivalent associate polytopes. The induction hypothesis is valid for length
$l+1$. This finishes the proof.
$\square$
\enddemo  

\proclaim{Corollary 4.4}
Let $A\in \Cal A_k$ and $B\in\Cal A_0$. Then $X_A\ed X_B\times (\Bbb S^1)^k$ if and only if $P_A=P_B$.
\endproclaim

\demo{Proof}
By Lemma 0.9, there exists $A'\in \Cal A_0$ such that the link $X_A$ is equivariantly diffeomorphic to $X_{A'}\times (\Bbb S^1)^k$.
In particular, this implies that $P_{A'}=P_A$. Now apply Rigidity Theorem 4.1.
$\square$
\enddemo

\proclaim{Corollary 4.5}
Let $\Phi\ :\ [0,1]\to \Cal A\cap M_{np}(\Bbb R)$ be a continuous path of admissible matrices of same dimensions. Set $A_t=\Phi(t)$. Then $X_{A_0}$ is equivariantly
diffeomorphic to $X_{A_1}$.
\endproclaim

\demo{Proof}
Let $I\subset \{1,\hdots, n\}$ such that $0$ belongs to the convex hull of $(((A_0)_i)_{i\in I})$. Then $0$ belongs to the convex hull of $(((A_t)_i)_{i\in I})$
for all $t$, otherwise there would be a time $t_0$ at which the weak hyperbolicity condition would be broken and the path $\Phi$ would not be a path of admissible
matrices. As a consequence of Lemma 0.12 and \lastnum[-2], the associate polytopes $K_{A_t}$ have all the same combinatorial type. Moreover this implies that all the
$X_{A_t}$ belong to the same $\Cal A_k$. We may thus conclude from Rigidity Theorem 4.1 that $X_{A_0}$ and $X_{A_1}$ are equivariantly diffeomorphic.
$\square$
\enddemo
 
\proclaim{Corollary 4.6}
Let $A\in\Cal A$ and $B\in \Cal A$ and $C\in \Cal A$. Then $X_C\ed X_A \times X_B$ (up to product by circles) if and only if $P_C=P_A\times P_B$.
\endproclaim

\demo{Proof}
It is an immediate consequence of Example 0.6 and Rigidity Theorem 4.1, noting that, in Example 0.6, we have $P_C=P_A\times P_B$.
$\square$
\enddemo

The second statement of the Rigidity  Theorem 4.1 is definitely false 
if we replace equivariant diffeomorphism by diffeomorphism. A counterexample is given in
\cite{LdM2}, p.242. We will see other interesting counterexamples in Section 6 (see Example 6.2).

We may now rely the two previous Sections in the following Theorem. 
As a direct consequence of the description of flips given in Section 2,
of the description of elementary surgeries given in Section 3 and of Rigidity Theorem 4.1, we have

\proclaim{Theorem 4.7}
Let $A\in\Cal A$ and let $B\in\Cal A$ with same dimensions $n$ and $p$. 
Assume that  $P_B$ is obtained from $P_A$ by performing
a flip of type $(a,b)$ along some simplicial face $F$. Then, $X_B$ is obtained (up to equivariant diffeomorphism) from $X_A$ by performing an elementary 
surgery of type $(a,b)$ along some $X_F$.
\endproclaim

As noted above, the converse of the Theorem is false. Indeed, in Examples 2.10, 2.11 and 2.12, we may perform elementary surgeries which will not correspond to flips. In other words,
{\it the class of links (up to equivariant diffeomorphism) is not stable under
elementary surgeries.}

\proclaim{Corollary 4.8}
Let $A\in\Cal A$. Then $X_A$ is obtained (up to equivariant diffeomorphism) from $\Bbb S^{2n-2p-1}\times (\Bbb S^1)^p$ by performing
a finite number of elementary surgeries.
\endproclaim

\demo{Proof}
Let $W$ be the simple polytope obtained from the product $P_A\times [0,1]$ by cutting off a neighborhood of a vertex of $P_A\times
\{1\}$ by a hyperplane (cf Lemma 2.3). 
Then $W$ is a cobordism between $P_A$ and the simplex of dimension $n-p-1$. If it is trivial, then $P_A$ is the $(n-p-1)$-simplex,
otherwise it can be decomposed
into a finite number of elementary cobordisms. Now apply Theorem 4.7 for each elementary cobordism and conclude in both cases with Corollary 1.4.
$\square$
\enddemo

\proclaim{Corollary 4.9}
Let $A\in\Cal A$ and let $B\in\Cal A$ with same dimensions. 
Assume that $X_B$ is obtained from $X_A$ by an elementary surgery.
Then there exists an equivariant cobordism between $X_A\times (\Bbb S^1)^2$ 
and $X_B\times(\Bbb S^1)^2$.
\endproclaim

\demo{Proof}
Let $k\in\Bbb N$ such that $A\in\Cal A_k$. Let $(a,b)$ be the type of the elementary surgery transforming $X_A$ into $X_B$. Let
$W$ be the corresponding elementary cobordism between $P_A$ and $P_B$. We define an integer $l$ as follows: if $a=1$, then $k>0$
by Definition 3.1 and we take $l=k-1$; otherwise $l=k$. By use of the Realization Theorem 0.13, 
there exists a link $X_C$ such that $P_C=W$ and $C\in\Cal A_l$. By Lemma 0.11 and Proposition 2.13, we know that $P_C$ has 
$n-l+2$ facets. As it has dimension $n-p$, then $C$ is a configuration of $n+2$ points in $\R^{p+1}$, so $X_C$ has
dimension $2n-p+2$. Using the fact that $P_A$ and $P_B$ are disjoint facets of $P_C$ and that $X_A$ and $X_B$ have dimension $2n-p-1$, 
we may embed by Proposition 1.1 the link $X_A\times\Bbb S^1$ (respectively $X_B\times\Bbb S^1$) as a smooth submanifold  
of $X_C$ of codimension $2$ with trivial normal bundle. The manifold obtained from $X_C$ by removing an open trivial tubular neighborhood
of each of these submanifolds is an equivariant cobordism between $X_A\times (\Bbb S^1)^2$ 
and $X_B\times(\Bbb S^1)^2$.
$\square$
\enddemo

\head
{\bf 5. Wall-crossing}
\endhead

We will now use the previous results to resolve the wall-crossing problem (compare with \cite{Bo}, \S 4). Let us start with an example to make the next explanations clearer.

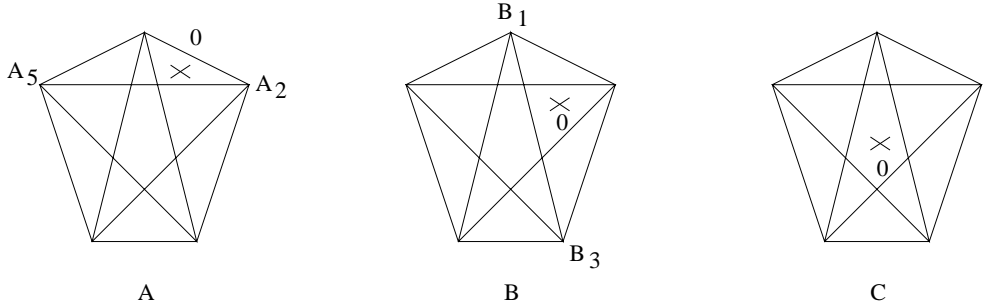
\hfil\scaledpicture 9.3in by 2.9in (fig6 scaled 500) \hfil
\vskip -.5cm

\example{Example 5.1}
Consider the links related to the three admissible configurations represented in the previous picture (the vertices of each 
configuration are numbered clockwise).

Here $n$ is equal to $5$ and $p$ to $2$. Note that $B$ and $C$ are translations of $A$ in $\Bbb R^2$. Nevertheless, the corresponding links are very different. From 
\cite{LdM1} (see Example 0.5) or \cite{McG}, we can conclude that
$$
\eqalign{
X_A&\ed \Bbb S^5\times \Bbb S^1\times \Bbb S^1 \cr
X_B&\ed \Bbb S^3\times \Bbb S^3\times \Bbb S^1 \cr
X_C&\ed \#(5)\Bbb S^3\times \Bbb S^4 
}
$$
where $\#(5)\Bbb S^3\times \Bbb S^4$ denotes the connected sum of five copies of $\Bbb S^3\times \Bbb S^4$. By Corollary 4.5, as long as we move smoothly the configuration
$A$ without breaking the weak hyperbolic condition, i.e. {\it without crossing a wall}, the link $X_A$ remains unchanged. But to go from $A$ to $B$ we have to cross
the wall $A_2A_5$, and to go from $B$ to $C$ we have to cross the wall $B_1B_3$; finally notice that we cannot pass directly from $A$ to $C$ with a single wall-crossing.
The best we can do is to perform two wall-crossings.
\endexample
 
\definition{Definition 5.2}
Let $A\in\Cal A$. A {\it wall} of $A$ is an hyperplane of $\Bbb R^p$ passing through $p$ vectors of $A$ and no more than $p$ (the data of the hyperplane is thus equivalent
to the data of the $p$ vectors) and which does not support a facet of $\Cal H(A)$.
\enddefinition

>From the definition, the intersection of the set $\{A_1,\hdots, A_n\}$ with each open half-space defined by the wall is not vacuous.

\definition{Definition 5.3}
Let $A\in \Cal A$ and $B\in\Cal A$ of same dimensions $n$ and $p$. Let $W$ be a wall of $A$. We say that $B$ is obtained from $A$ by {\it crossing the wall} $W$ if

\noindent (i) The configuration $B$ is a translate of $A$ by some vector $v$ of $\Bbb R^p$.

\noindent (ii) The configuration $A+tv$ is admissible for every $t$ in $[0,1]$ {\it except for one value} $t_0\in ]0,1[$.

\noindent (iii) At $t_0$, the point $0\in\Bbb R^p$ belongs to the translate of $W$ by $t_0v$ and does not belong to any other wall.
\enddefinition

In other words, $0$ ``moves'' continuously in the direction $-v$ and crosses the wall $W$, hence the terminology.

Let $A\in\Cal A$ and let $W$ be a wall of $A$. Then $W$ parts $\Bbb R^p$ into two open half-spaces containing the $n-p$ vectors of $A$ not belonging to $W$. More
precisely, one of the two open half-spaces, let us denote it by $W^+$, contains $0$ and $a$ vectors of $A$, whereas the other (that we call $W^-$) contains $b$
vectors of $A$. We say that the wall $W$ is of {\it type} $(a,b)$. We have $a+b=n-p$ and $1\leq a\leq n-p-1$ and $1\leq b\leq n-p-1$.

Now, let $B$ be obtained from $A$ by crossing $W$. If, by abuse of notations, we still call $W^+$ and $W^-$ the open half-spaces of $\Bbb R^p$ separated by
the translate of $W$, then $W^+$ still contains $a$ vectors of $B$ (which are exactly the translates of the $a$ vectors of $A$ lying in $W^+$) and $W^-$ contains $b$
vectors of $B$, but now $0$ lies in $W^-$. In particular, before the wall-crossing, $0$ belongs to the convex hull of the set consisting of the $p$ vectors of the wall
$W$ and any vector of $W^+$; after crossing the wall, $0$ belongs to the convex hull of the set consisting of the $p$ vectors of the wall
$W$ and any vector of $W^-$.

\proclaim{Wall-crossing Theorem 5.4}
Let $A\in \Cal A$ and $B\in\Cal A$ of same dimensions $n$ and $p$. Assume that $p>0$. Then, the following propositions are equivalent:

\noindent (i) The convex polytope $P_B$ is obtained from $P_A$ by a flip of type $(a,b)$ along the simplicial face $F_J$.

\noindent (ii) There exists $X_{B'}\ed X_B$ and $X_{A'}\ed X_A$ such that $X_{B'}$ is obtained from $X_{A'}$ by a single wall-crossing of $A'$, 
which is of type $(a,b)$.
\endproclaim

In the particular case where $p=0$, the notion of wall is meaningless. This explains the restriction $p>0$ in the statement of Wall-crossing Theorem 5.4.

Combining this result with Theorem 4.7 yields immediately

\proclaim{Corollary 5.5}
Under the same hypotheses, $X_B$ is obtained from $X_A$ by an elementary surgery of type $(a,b)$ along $X_{F_J}$.
\endproclaim

In other words, the class of links (up to equivariant diffeomorphism) is not stable under elementary surgeries but is 
{\it stable under elementary surgeries coming from wall-crossings}.

\demo{Proof of Wall-crossing Theorem 5.4}
The argument is purely convex. Assume (i). Then we can form the simple convex polytope $P_C$ with $P_A$ and $P_B$ as separated facets and with one single extra vertex 
of index $(a,b)$. Let $k\in\Bbb N$ such that $A\in\Cal A_k$. We define an integer $l$ as in the proof of Corollary 4.9:
if $a=1$, then $k>0$ (the assumption $p>0$ excludes the case $a=b=1$) and we take $l=k-1$; otherwise $l=k$.
Note that $P_C$ has dimension $n-p$ and has $n+2-l$ facets by Proposition 2.13. By Realization Theorem 0.13, there exists a link $X_C$ corresponding to $P_C$
with $C\in\Cal A_l$. We know that $C$ is a
configuration of $n+2$ vectors of $\Bbb R^{p+1}$. We set $C=(C_0,\hdots, C_{n+1})$. 
We may assume that $C_+=C\setminus\{C_0\}$ satisfies $X_{C_+}\ed X_A\times \Bbb S^1$
and that $C_-=C\setminus\{C_{n+1}\}$ satisfies $X_{C_-}\ed X_B\times \Bbb S^1$ (see Corollary 4.9). 
Moreover, as $P_A\cap P_B$ is vacuous (as a face of $P_C$), then $C\setminus\{C_0, C_{n+1}\}$
is not admissible. We say that $\{C_0,C_{n+1}\}$ is {\it indispensable}. In particular, this means that there exists an hyperplane 
of $\Bbb R^{p+1}$ passing through $0$
strictly separating $\{C_0,C_{n+1}\}$ from $\overline C=C\setminus\{C_0, C_{n+1}\}$. Scaling each vector of $\overline C$ by a strictly 
positive real number if necessary, 
we may assume that $\overline C$ lies in an affine hyperplane $H$ of $\Bbb R^{p+1}$ without changing the equivariant diffeomorphism type 
of $X_C$ (see Corollary 4.5).

Under this assumption, the convex hull of $C_+$ is a pyramid with base $\overline C$ and apex $C_{n+1}$ and containing $0$. In particular, $C_{n+1}$ is indispensable. 
This implies that, if we project
$0$ onto the hyperplane $H$ by letting
$$
\bar 0=H\cap (0C_{n+1})
$$
where $(0C_{n+1})$ denotes the line passing through the origin and through the point $C_{n+1}$; then, identifying $H$ with $\Bbb R^p$ and $\bar 0$ with the zero of 
$\Bbb R^p$ yields an admissible configuration $A'$
of $n$ vectors in $\R^p$ satisfying $X_{A'}\ed X_A$ (cf Lemma 0.9).

Performing the same transformation on the convex hull of $C_-$ viewed as a cone over $\overline C$ with apex $C_0$, 
we obtain an admissible configuration $B'$ of $n$ vectors in $\R^p$ satisfying $X_{B'}\ed X_B$ and such that
$B'$ is obtained from $A'$ by translation. 

The picture below should illustrate this construction. Taking $\bar 0$ as $O_1$ (respectively $O_2$) gives the configuration $A'$
(respectively $B'$).

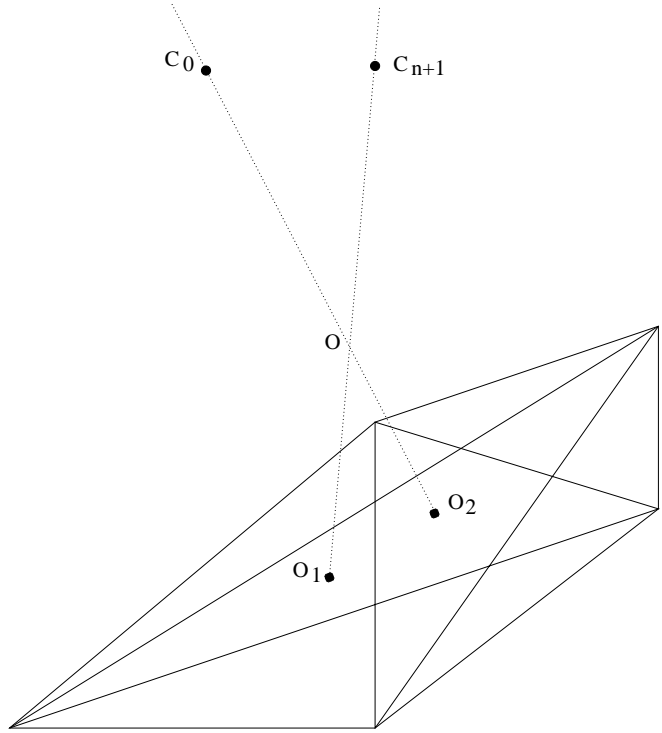
\hfil\scaledpicture 8.9in by 9.9in (fig11 scaled 350) \hfil

>From the construction, there is a translation sending the configuration $A'$ to $B'$. 
Let us now prove that this translation induces exactly one wall-crossing and characterize it.

\proclaim{Lemma 5.6}
Let $I\subset \{1,\hdots,n\}$ of cardinal $p$. Assume that $\{(A'_i)_{i\in I}\}$ defines a wall $W$ of $A'$. Then 
$W$ is crossed when changing from $A'$ to $B'$
if and only if $0$ is in the convex hull of $\{C_0,C_{n+1}\}\cup \{(C_i)_{i\in I}\}$.
\endproclaim

\demo{Proof of Lemma 5.6}
The proof is direct. Let $W$ be a wall of $A'$ defined by $I$. The hyperplane passing through $W$ and through $0$, let us call it $H_1$, separates $\Bbb R^{p+1}$
into two open half-spaces. Clearly, $W$ is crossed when changing from $A'$ to $B'$ if and only if $C_0$ and $C_{n+1}$ does not belong to the same open half-space.
If it is the case, then $H_1$ cuts the segment $[C_0,C_{n+1}]$ in one point $C_{t_0}$ and $0$ belongs to the convex hull of $\{C_{t_0}\}\cup
\{(C_i)_{i\in I}\}$. Therefore, $0$ is in $\Delta$, the convex hull of $\{C_0,C_{n+1}\}\cup \{(C_i)_{i\in I}\}$.

Conversely, assume that $C_0$ and $C_{n+1}$ belongs to the same open half-space defined by $H_1$. Then, the intersection of $\Delta$ and $H_1$ is included in $W$. Thus,
it does not contain $0$.
$\square$
\enddemo

Now, by Lemma 0.12 and by \lastnum[-2], a set of $p+2$ vertices of $C$ including $C_0$ and $C_{n+1}$ and containing $0$ in its convex hull corresponds to a vertex of $P_C$ 
which neither belongs
to $P_A$ nor to $P_B$. As the flip transforming $P_A$ into $P_B$ is elementary, there exists only one such simplex, and thus $B'$ is obtained from $A'$ by a single
wall-crossing along the wall $W_J$ corresponding to the extra vertex of $P_C$. Let us determine the type of the wall.

Let $I$ be the set of indices defining $W$. 
As before, let $W^+$ (respectively $W^-$) be the open half-space containing $\bar 0$ (respectively not containing $\bar 0$) before performing the wall-crossing. 
A point $A'_i$
belongs to $W^+$ if and only if the convex hull of $\{A'_i\}\cup \{A'_j\ \vert \ j\in I\}$ in $\Bbb R^p$ contains $\bar 0$. Since $0$ belongs to the segment
$[\bar 0, C_{n+1}]$, this is the case if and only if 
the convex hull of $\{C_{n+1}\}\cup \{C_i\}\cup \{C_j\ \vert \ j\in I\}$ contains $0$ in $\Bbb R^{p+1}$. Through \lastnum[-2], this determines a vertex $v$ of 
$P_A\subset P_C$. Moreover, since $0$ belongs to $\{C_0,C_{n+1}\}\cup \{C_j\ \vert \ j\in I\}$ by Lemma 5.6 and to 
$\{C_0,C_{n+1}\}\cup \{C_i\}\cup \{C_j\ \vert \ j\in I\}$, we know, still by \lastnum[-2], that there is an edge from $v$ 
to the extra vertex of $P_C$ (that is the vertex of $P_C\setminus (P_A\sqcup P_B)$). As this
vertex has index $(a,b)$, the wall $W$ separates $A'$ into $a$ vectors belonging to $W^+$ and $b$ vectors belonging to $W^-$.

Conversely, assume (ii). Let us define a new admissible configuration as follows: let 
$$
C_i=\pmatrix
{A'}_i \\
-1
\endpmatrix \in \Bbb R^{p+1}
\leqno 1\leq i\leq n
$$ 
and let $\bar 0=(0,-1)\in\Bbb R^p\times\Bbb R$. Consider the hyperplane $H=\Bbb R^p\times \{1\}\subset \Bbb R^{p+1}$. Let $C_0$ be the intersection of $H$ with the line
$(0\bar 0)$. We may now move $\bar 0$ inside $\Bbb R^p\times\{-1\}$ {\it without moving the points} $C_i$ to realize the wall-crossing from $A'$ to $B'$. Define
$C_{n+1}$ as the intersection of $H$ with $0\bar 0$ after the translation of $\bar 0$. Then $C$ is obviously an admissible configuration. We obtain exactly the same picture 
as before.

Moreover, $C\setminus \{C_{n+1}\}$ is an admissible configuration which is a pyramid with base $\overline C=(C_1,\hdots, C_n)$ and apex $C_0$, thus 
$$
X_{C\setminus\{C_{n+1}\}}=X_C\cap \{z_{n+1}=0\}\ed X_{A'}\times\Bbb S^1
$$
In the same way,
$$
X_{C\setminus\{C_{0}\}}=X_C\cap \{z_{0}=0\}\ed X_{B'}\times\Bbb S^1
$$

>From the construction, we obviously have $X_{\overline C}=\emptyset$. Therefore $P_C$ is a cobordism between $P_{A'}$ and $P_{B'}$. But as above, using Lemmas 0.12  and 5.6
and \lastnum[-2], it is straightforward to check that $P_C$ has a single extra vertex which is of index $(a,b)$ and that $P_C$ is an elementary cobordism between $P_A$ and 
$P_B$ along some simplicial face $F_J$.
$\square$
\enddemo

\proclaim{Corollary 5.7}
Let $A\in\Cal A$. Then there exists $A'\in\Cal A$ such that 

\noindent (i) The link $X_A$ is equivariantly diffeomorphic to $X_{A'}$.

\noindent (ii) The configuration $A'$ is obtained by wall-crossings
from a configuration $A''$ satisfying $X_{A''}\ed \Bbb S^{2n-2p-1}\times (\Bbb S^1)^p$.
\endproclaim

\demo{Proof}
Let $A'$ be a generic perturbation of $A$, that is a small perturbation of $A$ whose convex hull is simplicial. In this situation, an hyperplane of $\Bbb R^p$ contains at 
most $p$ vertices of $A'$. By Corollary 4.5,
we may assume that $X_{A'}\ed X_A$. For simplicity, assume that the convex hull of $(A'_1,\hdots A'_p)$ is a facet of $\Cal H(A'_1,\hdots,
A'_n)$. Consider the region $\Cal R$ of $\Bbb R^p$ defined as follows: $\Cal R$ is the union of the simplices whose vertices are
constituted by $p-1$ points among $(A'_1,\hdots A'_p)$ and two points among $(A'_{p+1},\hdots A'_n)$. 

The shaded region on the picture below is an example of such a $\Cal R$.

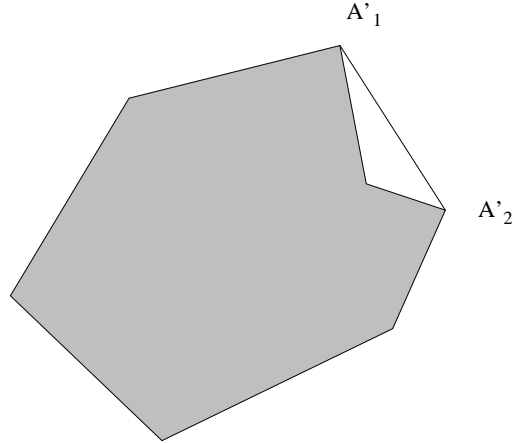
\hfil\scaledpicture 4.8in by 4.2in (fig12 scaled 500) \hfil

Notice that a point of
$\Cal H(A'_1,\hdots,
A'_n)$ which is sufficiently close to the center of $\Cal H(A'_1,\hdots A'_p)$ does not belong to $\Cal R$. Define $A''$ as an admissible
configuration obtained as a translate of $A'$
such that $0$ does not belong to the corresponding translate of $\Cal R$. In particular, $A''$ is obtained from $A'$ by wall
crossings. Then $A''_1$, ..., $A''_p$ are indispensable points of
$A''$, so by Lemma 0.9, we have that $A''\in\Cal A_k$ for $k\geq p$. This implies that $P_{A''}$ has dimension $n-p-1$ and has
at most $n-p$ facets. Therefore $k=p$ and $P_A$ is the $(n-p-1)$-simplex, so by Corollary 1.4 we have  
$X_{A''}\ed \Bbb S^{2n-2p-1}\times (\Bbb S^1)^p$.
$\square$
\enddemo

\remark{Remark 5.8}
Generically, we may take $A'=A$.
\endremark

\head
{\bf 6. Elementary surgery of type $\boldkey ( \boldkey 1 \boldkey , \boldkey n \boldkey )$}
\endhead

Let $X_A$ be a link. Assume that $P_A$ is obtained from the simplex (of same dimension) by performing uniquely flips of type $(1,n)$.
Then in this case, we may describe explicitely the diffeomorphism type of the link. First, note:

\proclaim{Lemma 6.1}
Let $A\in\Cal A_k$ with $k>1$. Let $X_B$ be obtained from $X_A$ by performing an elementary surgery 
of type $(1,n)$ along some invariant submanifold corresponding to a vertex. Then the
diffeomorphism type of $X_B$ is independent on the choice of the vertex on which the flip occurs.
\endproclaim

\demo{Proof}
Let $v$ and $v'$ be two vertices of $P_A$. We want to prove that, if $X_B$ and $X_{B'}$ denotes the links obtained from $X_A$ by 
performing an elementary surgery of type $(1,n)$ along $X_v$ (respectively $X_{v'}$), then these two links are diffeomorphic.
It is enough to show this in the case where $v$ and $v'$ belong to the same edge $E$. Let us describe $X_E$. By Corollary 1.5, the
link $X_E$ is diffeomorphic to $\Bbb S^3\times (\Bbb S^1)^p$. The real torus $(\Bbb S^1)^{p+2}=\Bbb S^1\times \Bbb S^1
\times T$ acts on $X_E$ in the following manner:
decompose $\Bbb S^3$ as the union of two solid tori $(\Bbb S^1\times \Bbb D^2)\times (\Bbb D^2\times\Bbb S^1)$. Then 
$\Bbb S^1\times\Bbb S^1$ acts on each solid torus in the natural way (that is the first factor by translations on $\Bbb S^1$ and the second factor 
tangentially to each circle
on $\Bbb D^2$) and this describes the induced action on $\Bbb S^3$; finally, $T$ acts by translations on $(\Bbb S^1)^p$. Therefore, $X_v$
is exactly given as $(\Bbb S^1\times \{0\})\times (\Bbb S^1)^p$, that is as the core circle of the first solid torus product with
$(\Bbb S^1)^p$; and $X_{v'}$ is exactly given as $(\{0\}\times \Bbb S^1)\times (\Bbb S^1)^p$, 
that is as the core circle of the second solid torus product with
$(\Bbb S^1)^p$. There exists an isotopy in $\Bbb S^3$ which sends $\Bbb S^1\times \{0\}$ to $\{0\}\times \Bbb S^1$ and this isotopy
can be extended by the identity on $(\Bbb S^1)^p$ to obtain an isotopy in $X_E$ sending $X_v$ to $X_{v'}$. Moreover, as it is
the identity on $(\Bbb S^1)^p$, it maps the circle which will be fulled by a 2-disk in the surgery giving $X_B$ to 
the circle which will be fulled by a 2-disk in the surgery giving $X_{B'}$. Therefore the two elementary surgeries give the
same result that is, $X_B$ is diffeomorphic to $X_{B'}$.
$\square$
\enddemo

Of course, in the previous Lemma, the class of $X_B$ {\it modulo equivariant diffeomorphisms} depends on the vertex on which the surgery
occurs: generally, the corresponding flips give different combinatorial types so, by Rigidity Theorem 4.1, different equivariant smooth
classes of links. Here is such an example.

\example{Example 6.2}
Consider the following polyhedron (the hexagonal book)

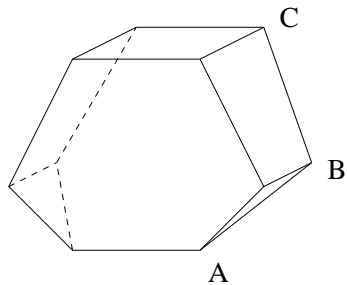
\hfil\scaledpicture 2.7in by 2.2in (fig9 scaled 600) \hfil
\vskip -.5cm

Let $X_A$ be the corresponding link with $A\in\Cal A_1$. Then, we may perform an elementary surgery of type $(1,3)$ on $X_A$ in
three manners, corresponding to the three vertices $A$, $B$ and $C$ indicated on the picture. By Lemma 6.1, the resulting
manifolds are all diffeomorphic but, by Rigidity Theorem 4.1, any two of them are not equivariantly diffeomorphic. In particular,
this gives an example of a manifold which admits three different ``structures of link''.
\endexample
\medskip
We may now describe explicitely the links corresponding to polytopes obtained from the simplex (of same dimension) by cutting off
vertices.

\proclaim{Theorem 6.3 (see [McG])}
Let $X_A$ be a simply-connected link such that $P_A$ is obtained from the $q$-simplex (of same dimension) by $l$ flips of type
$(1,n)$ (we assume that $l>0$). 
Then $X_A$ is diffeomorphic to the following connected sum of products of spheres:
$$
X_A\simeq \sc _{j=1}^l j \pmatrix l+1 \\ j+1 \endpmatrix \Bbb S^{2+j}\times\Bbb S^{2q+l-j-1} 
$$
\endproclaim

The proof of this Theorem is done for polygons in \cite{McG} (Theorem 3.4) but the proof of this generalization is the same.
Notice that this Theorem shows that, for any dimension of the associate polytope and for any value of $p$, there exist
infinite families which are connected sums of products of spheres as in Example 0.5.

Going back to Example 6.2, we see that the manifold 
$$
\sc (10) \Bbb S^3\times \Bbb S^8 \sc (20) \Bbb S^4\times\Bbb S^7 \sc (19) \Bbb S^5\times\Bbb S^6
$$
admits three different actions of $(\Bbb S^1)^8$ with quotient a convex polyhedron. 

This Example can be easily generalized as follows.

\example{Example 6.4}
Consider the $l$-gonal book $P_l$ 
for $l\geq 3$. It is obtained from the tetrahedron by $(l-3)$ flips of type $(1,3)$. By Theorem 6.3, it thus
gives rise to a $2$-connected link diffeomorphic to
$$
X_l=\sc _{j=1}^{l-3} j \pmatrix l-2 \\ j+1 \endpmatrix \Bbb S^{2+j}\times\Bbb S^{2+l-j}
$$
Consider a $l$-gonal facet of $P_l$. Number its vertices as indicated in the following picture.

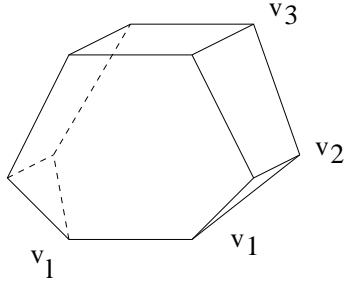
\hfil\scaledpicture 2.8in by 2.2in (fig10 scaled 600) \hfil
\vskip -.5cm

The simple convex polyhedra obtained from $X_{l-1}$ by cutting off a vertex $v_i$ are all combinatorially different when
$i$ ranges from $1$ to $[l/2]$ (where $[-]$ denotes the integer part). One of these polyhedra being the $l$-gonal book, 
we have by Lemma 6.1 that the corresponding links are all diffeomorphic to $X_l$.

In other words, the manifold $X_l$ admits at least $[l/2]$ structures of link. 
Therefore, the number of structures of link that $X_l$ has tends to infinity when $l$ tends to infinity.
 Notice that the dimension of $X_l$ is
$l+4$.
\endexample

\vfill
\eject
\specialhead
{\nom Part II: The cohomology ring of a link}
\endspecialhead
\bigskip

Thanks to Theorems 0.13 and 4.1, there is exactly one 2-connected link (up to equivariant diffeomorphism) associated to any simple
convex polytope (recall that we always consider a convex polytope only up to combinatorial equivalence). 
In this part, we give an explicit formula for the cohomology ring of a 2-connected
link in terms of its associate polytope. We use this formula to show that the cohomology of a link can have arbitrary amount of
torsion.

\head
{\bf 7. Notations and statement of the results}
\endhead

We denote by $P$ a simple convex polytope and by $X$ the associated 
2-connected link, that is we drop the subscript $A$ referring to
the choice of a matrix.

Given a finite simplicial complex $\Gamma$, {\it we make no distinction between $\Gamma$ and the poset of faces of $\Gamma$ ordered by inclusion.} In particular,
let $E$ be a set and $F$ a poset whose elements are subsets of $E$ ordered by inclusion. If every nonempty subset $J$ of $I\in F$ also belongs to $F$, then we consider
$F$ as a simplicial complex whose $k$-faces  are the elements of $F$ of cardinal $k+1$.   

Furthermore, we note:
\smallskip
\item{$\bullet$} $d$ the dimension of $P$;
\smallskip

\item{$\bullet$} $n$ the number of facets of $P$;
\smallskip

\item{$\bullet$} $\partial P$ the boundary of $P$. We consider it as a cell complex;
\smallskip

\item{$\bullet$} $P_b$ the barycentric subdivision of $\partial P$. In the same way, the barycentric subdivision of a simplicial complex $\Gamma $ will be denoted 
$\Gamma _b $. If a set $I$ numbers a simplex $\sigma$ of $\Gamma$, then we number the center of $\sigma$ in $\Gamma_b$ by the same set $I$, that is we identify
a simplex of $\Gamma$ and its center in $\Gamma_b$;
\smallskip

\item{$\bullet$} $\Cal F $ the set of the facets of $P$;
\smallskip

\item{$\bullet$} $\Cal I $ a subset of $\Cal F $;
\smallskip

\item{$\bullet$} $\vert \Cal I\vert$ the cardinal of $\Cal I$;
\smallskip

\item{$\bullet$} $\bar{\Cal I }$ the complement of $\Cal I $ in $\Cal F $;
\smallskip

\item{$\bullet$} $F_{\Cal I }$ the intersection of the facets of $P$ that are in $\Cal I $. It 
is either empty or a face of $P$;
\smallskip

\item{$\bullet$} $\Delta $ the poset of nonempty subsets $\Cal I $ of $\Cal F$ such that 
$F_{\bar{\Cal I }} =\emptyset $ ordered by inclusion. It is a simplicial complex;
\smallskip

\item{$\bullet$} $P _{\Cal I }$ the union of the facets of $P$ that are in $\Cal I $;
\smallskip

\item{$\bullet$} $K_{\Cal I }$ the poset of nonempty subsets $I$ of $\Cal I$ such that $F_I$ is a (nonempty) face of $P$ ordered by inclusion. It is a simplicial complex. 
We will often consider its barycentric subdivision $(K_{\Cal I})_b$ as a simplicial subcomplex of $P_b$ by identifying a subset $I$ to the center of the face $F_I$ in the 
barycentric subdivision of $\partial P$;
\smallskip

\item{$\bullet$} $\tilde{\Cal I }$ the poset of proper subsets $I$ of 
$\bar{\Cal I }$ such that $F_{\bar{\Cal I } \backslash I}$ is not empty ordered by {\it reverse inclusion}. It is also a simplicial complex;
\smallskip

\item{$\bullet$} $\hat I$ the complement of a subset $I$ in $\bar{\Cal I}$;
\smallskip

\item{$\bullet$} $P ^* $ the dual polytope of $P$;
\smallskip

\item{$\bullet$} $\delta_i^j$ the Kronecker symbol;
\smallskip

\item{$\bullet$} $H_i(A, \Bbb Z)$ (respectively $\tilde H_i(A, \Bbb Z)$) the $i$-th homology group (respectively reduced homology group) of a manifold or a simplicial
complex $A$ with coefficients in $\Bbb Z$. By convention, we set $\tilde H_{-1}(\emptyset, \Bbb Z)=\Bbb Z$;
\smallskip
 
\item{$\bullet$} $H^i(A, \Bbb Z)$ (respectively $\tilde H^i(A, \Bbb Z)$) the $i$-th cohomology group (respectively reduced cohomology group) of a manifold or a simplicial
complex $A$ with coefficients in $\Bbb Z$;
\smallskip

\item{$\bullet$} the simplex whose vertices are the elements of a finite set $E$ will be denoted 
$\Delta _E $ and its boundary $S^E $ (in some context, $\Delta _E $ will be noted 
$\sigma _E $);
\smallskip

\definition{Definition 7.1}
  For a nonempty face $F$ of $P$, the vector space underlying the affine space in 
which $F$ has nonempty interior will be called the (vector) space of $F$. By abuse of notation, {\it we will still denote by $F$ the space of $F$}. No confusion
should arise from this abuse.
\enddefinition

\definition{Definition 7.2}
  A proper face of $P$ will be called an $\Cal I $-face (respectively an $\bar{\Cal I }$-face) if 
every facet of $P$ containing it is in $\Cal I $ (respectively in $\bar{\Cal I }$).
\enddefinition

We prove now some preliminary results on simple polytopes.

\proclaim{Lemma 7.3}
 Let $P$ be a simple polytope and let $\Cal I\subset \Cal F$. Then, a 
nonempty intersection of elements of $\Cal I $ is an $\Cal I $-face.
\endproclaim

\demo{Proof}
This comes directly from the fact that the neighborhood of a face in a simple polytope is the product of
this face by a simplex. Hence, for every face $F$ of $P$, there is a unique subset $\Cal I$ such that $F_{\Cal I}=F$ and $F$ is an $\Cal I$-face.
$\square$
\enddemo

This Lemma is false for non-simple polytopes. In the following picture, the polytope is a pyramid with rectangular base and apex $v$, whereas the set $\Cal I$ consists
of two faces whose intersection is $v$. Nevertheless, $v$ is not an $\Cal I$-face.

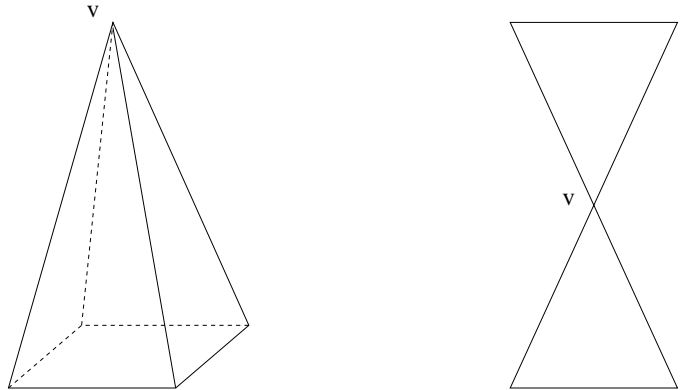
\hfil\scaledpicture 8.0in by 4.6in (fig13 scaled 400) \hfil

We then have

\proclaim{Lemma 7.4}
  Let $P$ be a simple polytope. Consider a subset $\Cal I $ of $\Cal F$. Then,

\noindent (i) The complex $(K_{\Cal I })_b$ is homotopy equivalent to 
${P }_{\Cal I }$.

\noindent (ii) The set ${P }_{\Cal I }$ has the same homotopy type as its interior in $\partial P$.
\endproclaim

\demo{Proof}
  The barycentric subdivision of $\partial P$ is a simplicial complex whose vertices 
are all the (nonempty) faces of $P$. By Lemma 7.3, the complex $(K_{\Cal I })_b$ is isomorphic to the 
subcomplex of this subdivision associated to $\Cal I $-faces. Each point $M$ of $P_{\Cal I }$ 
belongs to a {\it unique} minimal simplex of $P_b$
and this simplex has at least one vertex belonging to $(K_{\Cal I})_b$ (the center of the minimal face which contains it). Take the barycentric
coordinates of $M$ in this simplex. We may then construct a 
retraction of ${P }_{\Cal I }$ on $(K_{\Cal I })_b$ by cancelling the bad barycentric 
coordinates (i.e. coordinates associated to vertices which do not belong to $(K_{\Cal I })_b$).

To prove (ii), just remark that the previous construction yields also a retraction of the interior of $P_{\Cal I}$ onto $(K_{\Cal I })_b$.
$\square$
\enddemo

This Lemma is in fact a variation of the following well known fact:

\proclaim{Lemma 7.5}
  Let $\Delta $ be a simplicial complex, $\Gamma $ a subcomplex. Then the 
"mirror complex" of $\Gamma $ in $\Delta $, i.e. the complex of the faces of 
$\Delta _b$ that are disjoint from $\Gamma _b $ is homotopy 
equivalent to (and even a deformation retract of) $\Delta \setminus \Gamma $.
\endproclaim

We may now state

\proclaim{Cohomology Theorem 7.6}
  For any $i$, we have an isomorphism:
$$
H^i (X,\Bbb Z ) \simeq 
\bigoplus_{\Cal I \subset \Cal F } \tilde{H}_{d+|\bar{\Cal I }|-i-1} (P _{\Cal I} ,\Bbb Z )
$$

  We note $\psi ([c])$ the inverse image by this isomorphism of a class $[c]$ in any factor 
of the second member.

  Moreover, consider two classes $[c] \in \tilde H_k (P _{\Cal I } ,\Bbb Z )$ and 
$[c'] \in \tilde H_{k'} (P _{\Cal J },\Bbb Z )$. Note $[c] \cap [c']$ their intersection class in 
$\tilde H_{k+k'-d+1} (P _{\Cal I \cap \Cal J}, \Bbb Z )$. Then, up to sign, the cup product 
of their images by $\psi $ is given by~:
\setbox1=\hbox{$\psi([c] \cap [c'])$}
\setbox2=\hbox to \wd1{\hfill $0$ \hfill}
$$
\psi([c]) \smile \psi([c']) = \left\{ 
\eqalign{
\box1 & \quad\text{ if } \Cal I \cup \Cal J = \Cal F \cr
\box2 & \quad\text{ else} 
}\right . 
$$
\endproclaim

\remark{Remark 7.7}
The following formula for the homology groups of $X$ in terms of $P^*$ also holds:
$$
\tilde H_i(X,\Bbb Z)\simeq \bigoplus_{\Cal I\subset \Cal F} \tilde H_{i-|\Cal I|-1}(P_{\Cal I}^*,\Bbb Z)
$$
where $\Cal F$ is identified with the set of vertices of $P^*$ and where $P_{\Cal I}^*$ denotes the maximal simplicial subcomplex
of $P^*$ with vertex set $\Cal I$. In some cases, this formula is easier to use to compute the homology groups. We will prove this formula at the same time as the formula
of Cohomology Theorem 7.6.
\endremark

\remark{Remark 7.8}
If $\Cal I $ and $\Cal J $ are complementary in $\Cal F $ and we take classes 
$[c] \in \tilde{H}_k (P _{\Cal I }, \Bbb Z )$ and 
$[c'] \in \tilde{H}_{k'} (P _{\Cal J },\Bbb Z )$ with $k+k' = d-2$, then their intersection 
class in $\tilde{H}_{-1} (\emptyset , \Bbb Z ) \simeq \Bbb Z $ is their linking number. In particular, {\it Poincar\'e duality on $X$ is given by Alexander duality on
$\partial P$}.
\endremark

\example{Example 7.9}
Let $P$ be the cube. Number its facets in the following way: $1$, $2$ and $3$ denote three faces adjacent to a vertex $v$, 
whereas $1'$ (respectively $2'$, $3'$) is the opposite
face to $1$ (respectively $2$, $3$).

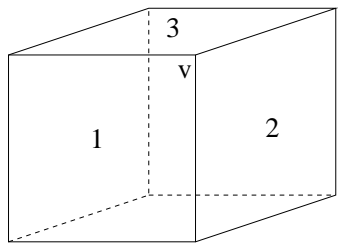
\hfil\scaledpicture 3.5in by 2.5in (fig14 scaled 450) \hfil

The sets $P_{\{1,2,1',2'\}}$, $P_{\{1,3,1',3'\}}$ and $P_{\{2,3,2',3'\}}$ have the homotopy type of a circle. Let us denote by $[c_{12}]$ (respectively
$[c_{13}]$ and $[c_{23}]$) a generator of $\tilde H_1(P_{\{1,2,1',2'\}},\Bbb Z)$ (respectively $\tilde H_1(P_{\{1,3,1',3'\}},\Bbb Z)$ and
$\tilde H_1(P_{\{2,3,2',3'\}},\Bbb Z)$).

The sets $P_{\{1,1'\}}$, $P_{\{2,2'\}}$ and $P_{\{3,3'\}}$ have the homotopy type of a pair of points. Let us denote by $[c_1]$ (respectively $[c_2]$ and $[c_3]$)
a generator of $\tilde H_0(P_{\{1,1'\}},\Bbb Z)$ (respectively $\tilde H_0(P_{\{2,2'\}},\Bbb Z)$ and
$\tilde H_0(P_{\{3,3'\}},\Bbb Z)$).

Finally, let us denote by $[c]$ a generator of the top-dimensional cohomology group of the associated link $X$.

Cohomology Theorem 7.6 gives the cohomology groups of $X$.
$$
\vbox{
\def\tv{\vrule height 12pt depth 5pt width .5pt}
\offinterlineskip
\halign{\tv \hfill \kern .7em {$#$} \kern .7em \hfill \tv &&\hfill \kern .7em {$#$}\kern .7em\hfill\tv\cr
\noalign{\hrule}
i & \tilde H^i (X,\Bbb Z)\simeq \cr
\noalign{\hrule}
1,2,4,5,7,8 & \{0\}  \cr  
\noalign{\hrule}
3 & \Bbb Z\cdot \psi([c_{12}])\oplus \Bbb Z\cdot \psi([c_{13}])\oplus \Bbb Z\cdot \psi([c_{23}])\cr
\noalign{\hrule}
6 & \Bbb Z\cdot \psi([c_1])\oplus \Bbb Z\cdot \psi([c_2])\oplus \Bbb Z\cdot \psi([c_3])\cr
\noalign{\hrule}
9 & \Bbb Z\cdot [c]\cr
\noalign{\hrule}
}}
$$
and the only non-zero cup products are, up to sign,
$$
\eqalign{
&\psi([c_{12}])\smile\psi([c_3])=\psi([c_{13}])\smile\psi([c_2])=\psi([c_{23}])\smile\psi([c_1])=[c] \cr
&\psi([c_{12}])\smile\psi([c_{13}])=\psi([c_1])\cr
&\psi([c_{12}])\smile\psi([c_{23}])=\psi([c_2])\cr
&\psi([c_{13}])\smile\psi([c_{23}])=\psi([c_3])
}
$$
>From Corollary 4.6 and Example 0.6, we know that $X$ is a product of spheres $\Bbb S^3\times\Bbb S^3\times\Bbb S^3$. We recover here its cohomology ring.
\endexample

\demo{Proof of the first part of Theorem 7.6 and of Remark 7.7}
By Lemma 0.7, the link $X$ has the same homotopy type as the complement $\Cal S$ of a coordinate 
subspace arrangement $\Cal L$ of $\Bbb C^n$ (see \lastnum[-12]
and \lastnum[-11]; as for the case of 
$X$ and $P$, we drop the subscript referring to a matrix $A$). Notice that $\Cal L$ is only defined up to a permutation of the coordinates of $\Bbb C^n$. Now, 
fix a numbering of the facets of $P$ by integers from $1$ to $n$. Then, by \lastnum[-3], this indeterminacy on $\Cal L$ is cancelled.

We make use of the formulas given in \cite{DL} which describe the
cohomology ring of a coordinate subspace arrangement.
  Let us first recall De Longueville's notations and results adapted to our case.

Let $\Delta$ be the simplicial complex defined at the beginning of this Section. Let $(e_1,\hdots,e_n)$ be the canonical basis of
$\C^n$. We may associate to $\Delta$ the following coordinate 
subspace arrangement
$$
\Cal A_{\Delta}=\{\text{Vect}_{\C}(e_i)_{i\in \Cal I} \quad\vert\quad \Cal I\subset\Delta\}
\tag\numerote
$$

Using \lastnum[-4], it is straightforward to check that

\proclaim{Lemma 7.10}
We have $\Cal A_{\Delta}=\Cal L$.
\endproclaim

Finally, let $\sigma$ be a face of $\Delta$; we define
$$
link_{\Delta} \sigma=\{\tau\in\Delta \quad\vert\quad \sigma\cap\tau=\emptyset, \ \sigma\cup\tau\in\Delta\}
$$

Geometrically, $link_{\Delta}\sigma$ is the boundary of the subcomplex of $\Delta$ formed by the simplices to which $\sigma$ belongs.

\remark{Remark 7.11}
Let $\sigma_{\Cal I}$ be a face of $\Delta$ indexed by $\Cal I\subset \Cal F$. Then, we have
$$
link_{\Delta} \sigma_{\Cal I}=\{ I\subset \bar{\Cal I} \quad\vert\quad F_{\hat I}=F_{\bar {\Cal I}\setminus I}=\emptyset\}\ .
$$ 
Therefore, the set $\tilde {\Cal I}$ defined at the beginning of
this Section is exactly the set of nonempty subsets of $\bar{\Cal I}$ which are not in $link_{\Delta}\sigma_{\Cal I}$.
\endremark
\medskip
 
With these notations, the Goresky-Mac Pherson formula of \cite{G-McP} states that the reduced cohomology group 
$\tilde{H}^i (\Cal S,\Bbb Z )$ is isomorphic to the sum of the groups 
$\tilde{H}_{2|\bar{\Cal I }|-i-2} (link_{\Delta } \sigma _{\Cal I},\Bbb Z )$, the sum being taken 
over all the elements $\sigma $ in $\Delta $. As $\Cal S$ and $X$ are homotopy equivalent, 
the same result is also true for $X$.

  On the other hand, consider two elements $\sigma _1 $ and $\sigma _2 $ of the complex $\Delta$ and two classes $[c_1]$ and $[c_2]$ 
of $\tilde{H}_{2|\sigma |-i-2} (link_{\Delta } \sigma _1 , \Bbb Z )$ and 
$\tilde{H}_{2|\sigma |-i-2} (link_{\Delta } \sigma _2 , \Bbb Z )$ represented by $c_1$ and $c_2$. Noting also $\psi ([c])$ 
the cohomology class associated to some class $[c]$, De Longueville shows in \cite{DL} that, up to sign:
\setbox1=\hbox{$\psi ([\langle i_2 - i_1 \rangle  * c_1 * c_2])$}
\setbox2=\hbox to \wd1 {\hfill $0$ \hfill}
$$
\psi([c_1] ) \smile \psi([c_2] ) = \left\{ \eqalign{
\box1 & \quad\text{ if }\sigma _1 \cup \sigma _2 = \Cal F\cr
\box2 & \quad\text{ else} 
}\right . 
$$
where $i_1 $ and $i_2 $ are elements out of $\sigma_1 $ and $\sigma_2 $ respectively, and where $*$ denotes the join of two cycles.

  To prove the Theorem, we will establish isomorphisms between the groups which compose 
the cohomology of $X$, then study the behaviour of the product in the polytopal case.

\proclaim{Lemma 7.12}
  For any $\Cal I $, the group
$\tilde{H}_{d+|\bar{\Cal I }|-i-1} (P _{\Cal I },\Bbb Z )$ is 

\item{$\bullet$} isomorphic to 
$\tilde{H}_{2|\bar{\Cal I }|-i-2} (link_{\Delta } \sigma _{\Cal I },\Bbb Z )$ if 
$\Cal I $ is in $\Delta $;

\item{$\bullet$} zero if $\Cal I $ is not in $\Delta $ and not $\Cal F $;

\item{$\bullet$} zero if $\Cal I = \Cal F $ and $i \neq 0$;

\item{$\bullet$} isomorphic to $\Bbb Z $ if $\Cal I = \Cal F $ and $i=0$.
\endproclaim

\demo{Proof of Lemma 7.12}
  Let us begin with the simple special case: $\Cal I = \Cal F $. In this case, 
$\Cal F $ is not in $\Delta $ and $P _{\Cal I }$ is all $\partial P$.

  We then have 
$\tilde{H}_{d+|\bar{\Cal I }|-i-1} (P _{\Cal I },\Bbb Z ) = 
\tilde{H}_{d-i-1} (\Bbb S^{d-1} ,\Bbb Z )$ which is zero, except if $d-i-1 = d-1$, i.e. $i=0$ 
in which case this group is isomorphic to $\Bbb Z $.

  Consider now that $\Cal I $ is not in $\Delta $ and not $\Cal F $. Then the facets of 
$\bar{\Cal I }$ exist and intersect. The set $P _{\bar{\Cal I }}$ is therefore starshaped in 
$\partial P$ and then so is $P _{\Cal I }$ ($\partial P$ is considered 
as a sphere). Hence, $P _{\Cal I }$ is contractible and all its reduced homology 
groups vanish.
 
  We will establish that, in the other cases, $link_{\Delta } \sigma _{\Cal I }$ and $P _{\Cal I }$ have 
complements in some spheres that are homotopy equivalent. The isomorphism will follow 
from Alexander duality applied twice.

  First, except in the special case thereup, $link_{\Delta } \sigma _{\Cal I }$ is a 
subcomplex of $S^{\bar{\Cal I }}$, which is a sphere of dimension $n-|\Cal I |-2$. By Lemma 7.5, its 
complement in this sphere is homotopy equivalent to its mirror complex. Here this mirror complex is the 
subcomplex of the barycentric subdivision of $S^{\bar {\Cal I }}$, whose vertices are the 
ones corresponding to elements of $\tilde{\Cal I}$ (see Remark 7.11), i.e. is isomorphic to the simplicial complex 
${\tilde{\Cal I}}_b$. Hence, by Alexander duality,
$\tilde H_{2|\bar{\Cal I }| - i - 2}(link_{\Delta } \sigma _{\Cal I }, \Bbb Z )$ is isomorphic to 
$\tilde H^{i - |\bar{\Cal I }| - 1} ({\tilde{\Cal I }}_b, \Bbb Z )$ and thus to $\tilde H^{i - |\bar{\Cal I }| - 1} ({\tilde{\Cal I }}, \Bbb Z )$.

  On the other side, $P _{\Cal I}$ is the complement in $\partial P$ of 
$P _{\bar{\Cal I }}$ (of its interior precisely but, by Lemma 7.4, they are 
homotopically equivalent). Still by Lemma 7.4, $P _{\bar{\Cal I }}$ is homotopically equivalent to $(K_{\bar{\Cal I }})_b$. Then, by 
Alexander duality, we get an isomorphism between 
$\tilde H^{i - |\bar{\Cal I }| - 1} (K_{\bar{\Cal I }}, \Bbb Z )$ and 
$\tilde H_{|\bar{\Cal I }| + d - i - 1} (P _{\Cal I }, \Bbb Z )$.

  But we claim that the complexes $\tilde{\Cal I }$ and $K_{\bar{\Cal I }}$ are 
isomorphic. In fact, by definition of $\tilde{\Cal I}$, the map 
$I \to \hat I$ sends $\tilde{\Cal I}$ to the set of $\bar{\Cal I }$-faces, 
reversing inclusion.
$\square$
\enddemo

  It is now easy to complete the proof of the first part of Theorem 7.6. Finally, noting that $P^*_{\bar{\Cal I}}$ is isomorphic to $K_{\bar{\Cal I}}$ for 
$\Cal I\in \Delta$, we deduce from the proof of Lemma 7.12 that $\tilde H^{2|\bar{\Cal I }| - i - 2}(link_{\Delta } \sigma _{\Cal I }, \Bbb Z )$ is isomorphic to
$\tilde H_{i - |\bar{\Cal I }| - 1} (P^*_{\bar{\Cal I}}, \Bbb Z )$. This leads to the formula of Remark 7.7.
$\square$
\enddemo

\remark{Notation 7.13}
  For a class $[c]$ in $\tilde{H}_k (link_{\Delta } \sigma _{\Cal I },\Bbb Z )$, its image in 
$\tilde H_{k + |\bar{\Cal I }| + d + 1}(P_{\Cal I},\Bbb Z)$ by the forementioned isomorphism will be denoted 
$\phi ([c])$.
\endremark
\medskip

  In order to prove the second part of Theorem 7.6, we have to explicitely establish 
the correspondance between the groups. As we need to explicitely compute Alexander duals, we have 
to deal with orientations.

\head
{\bf 8. Orientation}
\endhead

  We talk here about Alexander duality on spheres of the form $S^{\Cal I }$ for subsets 
$\Cal I $ of $\Cal F $ and on the sphere $\partial P$. These spheres have then to be 
oriented (in fact, this is not really necessary as long as we work up to sign, but even 
then suitable choices a bit simplify matters). Let us start with the orientation of $\partial P$. We consider $P$ as being realized
in $\R^d$. We orient $\R^d$ and thus obtain an orientation of $P$.
\medskip
\noindent $\underline{\text{Orientation of a facet and of a boundary:}}$
recall that if we consider an oriented po\-lytope, there is a canonical orientation of 
its boundary by stating that for any facet $F$ of this polytope, a basis consisting of the 
normal outward pointing vector followed by a positively oriented basis of the space of 
the facet is a positively oriented basis of the space of the polytope.
\medskip
\noindent $\underline{\text{Orientation of a face of }P\text{ :}}$
  consider a $k$-tuple $(H_1 ,\hdots ,H_k )$ of facets of $P $ with nonempty intersection. 
Then $F_{(H_1 ,\hdots ,H_k )}$ denote the intersection of these facets endowed with the 
following orientation: taking a basis $(v_1 ,\hdots ,v_k ,\Cal B )$ of the space of $P $, 
where $v_i $ denotes the normal outward pointing vector of $H_i$ and $\Cal B $ is a basis of 
the space of our face, we state that both basis have the same orientation. Remark that 
even a $0$-dimensional face has two "orientations". 

\remark{Remark 8.1}
To orient a face of $P$ is equivalent to
order the set of facets containing it. In particular, given an orientation of a 
convex polytope, {\it there is no canonical orientation of the faces which are not facets}.
\endremark 
\definition{Definition 8.2}
  A $d$-tuple $(H_1 ,\hdots , H_d )$ of facets of $P $ with nonempty intersection will be 
called direct if $(v_1 ,\hdots , v_d )$ is a positively oriented basis. It will be called undirect else.
\enddefinition

\remark{Notation 8.3}
  For a $k$-tuple $I = (H_1 ,\hdots ,H_k )$ and a $k'$-tuple $J = (H'_1 ,\hdots ,H'_{k'} )$ 
disjoint from $I$ of facets of $P $ such that $F_I $ and $F_J $ have nonempty 
intersection, the face associated to the $(k+k')$-tuple $(H_1 ,\hdots ,H_k ,H'_1 ,\hdots ,H'_{k'})$ 
will be denoted $F_{I+J}$.
\endremark

\medskip
\noindent $\underline{\text{Orientation of an intersection:}}$
consider a $n$-dimensional oriented vector space $E$ and two oriented subspaces $F$ and 
$F'$, of respective strictly positive dimension $d$ and $d'$ and whose sum is $E$. Then the vector space 
$F \cap F'$ is oriented with the convention that if $\Cal B = (v_1 ,\hdots ,v_{d+d'-n})$ is a 
basis of $F \cap F'$, if $(w_1 ,\hdots ,w_{n-d'},v_1 ,\hdots ,v_{d+d'-n})$ is a positive basis of 
$F$ and $(v_1 ,\hdots ,v_{d+d'-n},w'_1 ,\hdots ,w'_{n-d})$ a positive basis of $F'$, then the 
basis $\Cal B $ of $F \cap F'$ and the basis 
$(w_1 ,\hdots ,w_{n-d'},v_1 ,\hdots ,v_{d+d'-n},w'_1 ,\hdots ,w'_{n-d})$ have the same sign. In the special case where $F\cap F'$ is reduced to $\{0\}$, then we
state that $F\cap F'$ is positively oriented if $(w'_1 ,\hdots ,w'_{n-d},w_1 ,\hdots ,w_{n-d'})$ is a positive basis of $\Bbb R^d$. This convention is taken to
guarantee the statement of Lemma 8.5 (see below) in this special case.

\remark{Remark 8.4}
With this definition, the orientations of $F\cap F'$ and $F'\cap F$ may be different. 
\endremark
\medskip

The previous convention is a generalization of the convention of orientation of a face, since we have:

\proclaim{Lemma 8.5}
  With the orientation conventions thereup, $F_{I+J}$ is equal to 
$F_I \cap F_J $ as oriented face.
\endproclaim

\demo{Proof}
We use Notation 8.3. Let $v_i$ (respectively $v'_i$) denote the normal outward pointing vector of $H_i$ (respectively ${H_i}'$). 
We may assume that $F_I$ and $F_J$ are orthogonal. Let $\Cal B$
be a basis of $F_I\cap F_J$.
Then $(v_1,\hdots, v_k,v'_1,\hdots,v'_{k'},\Cal B)$ is a positive basis of $\Bbb R^d$ if and only if $(v'_1,\hdots,v'_{k'},\Cal B)$ is a positive basis of $F_I$ whereas
$(v'_1,\hdots,v'_{k'}, \Cal B, v_1,\hdots v_k)$ is a positive basis of $\Bbb R^d$ if and only if
$(\Cal B,v_1,\hdots,v_k)$ is a positive basis of $F_J$. The claim follows then easily. 
$\square$
\enddemo

\proclaim{Lemma 8.6}
Let $P$ be an oriented polytope. Let $F$ be a face of $P$. Fix an orientation of $F$.
  With the orientation conventions thereup, the oriented boundary of $F$ is given 
by:
$$
\partial F = 
\sum_{H \in \Cal F , F  \cap H \neq F,\emptyset } F\cap H
$$
where $F$ is considered as an oriented polytope and $H$ is endowed with the canonical 
orientation of $\partial P$.
\endproclaim

\demo{Proof}
We may find $I = (H_1 ,\hdots ,H_k )$ such that $F_I=F$ as oriented face. Now, set $\Cal F=\{H_1,\hdots, H_n\}$. 
For $k<i\leq n$, the oriented face $F_{I+\{i\}}$ is a facet of $F_I$ (if non-empty) which
is easily seen to be positively oriented with respect to the convention about the orientation of a facet. Therefore,
$$
\partial F=\sum_{k<i\leq n} F_{I+\{i\}}
$$
The result follows now from Lemma 8.5.
$\square$
\enddemo

Now, we orient the spheres $S^{\Cal I}$. We consider an order on $\Cal F $ and for any subset $\Cal I $ of $\Cal F $ we orient 
$\Delta _{\Cal I }$ compatibly with the restriction to $\Cal I $ of the order on $\Cal F $ as explained below. Then, $S^{\Cal I}$ is oriented as boundary of
$\Delta _{\Cal I }$.
\medskip

\noindent $\underline{\text{Orientation of a simplex:}}$
consider a finite set $E$ having at least two elements and a total order $\leq$. We can associate to this order an orientation of $\Delta _E $ by stating that if 
$e_0 \leq \hdots \leq e_{|E|-1}$ are the ordered elements of $E$, then the basis 
$\overrightarrow{e_0 e_1 } ,\hdots , \overrightarrow{e_0 e_{|E|-1}}$
is a positively oriented basis of the space of $\Delta _E $. The order and the 
orientation are then called compatible.

\remark{Convention 8.7}
In the sequel, a subset $\Cal I$ of $\Cal F$ will always be considered as an ordered set, with the order induced from the order of $\Cal F$.
In particular, the simplex $\sigma_{\Cal I}$ is thus an oriented simplex.
\endremark

\remark{Notation 8.8}
Let $E$ and $F$ be disjoint finite sets with orders $\leq _E $ on $E$ and 
$\leq _F $ on $F$. Then $EF$ denotes the set $E\cup F$ endowed with the following order:
any element of $E$ is less than any element of $F$ and 
the restriction of the order to $E$ (respectively to $F$) is $\leq _E $ (respectively $\leq _F $).
\endremark
\medskip

We finish this Section with another convention of orientation that will be needed.
\medskip

\noindent $\underline{\text{Orientation of a join:}}$
consider two oriented simplices $\Delta _E $ and $\Delta _F $ on disjoint finite sets 
$E$ and $F$, whose orientations are compatible with the orders $\leq _E $ on $E$ and 
$\leq _F $ on $F$. We orient their join $\Delta _E * \Delta _F $ compatibly with the 
order on $EF$. We easily check that 
this orientation does not depend on the chosen orders. Indeed, we have $\Delta_E * \Delta_F=\Delta_{EF}$.
\medskip
To sum up, {\it given a total order on $\Cal F$, then, with the conventions thereup, an orientation is fixed on any face of $P$ as well as on any sphere 
$S^{\Cal I}$ for $\Cal I\subset \Cal F$}.

\head
{\bf 9. Alexander duals up to a sign}
\endhead

To compute Alexander duals, we make use of \cite{Al}, t. 3, ch. XIII. We first recall this construction in our context.
Let $P$ be an oriented simple convex polytope. Let $K$ be $\partial P$ seen as a cell complex. Let $m$ be its dimension. 
Given an oriented 
cell $\sigma$ of $K$, 
its star dual $\sigma^*$ is defined as the maximal 
subcomplex of the barycentric subdivision $K_b$ of $K$ whose vertices are the centers of
the faces of $K$ containing $\sigma$ (see \cite{Al}, t. 1, p.143--144). An orientation is fixed on $\sigma^*$ by demanding that the
intersection number of $\sigma$ with $\sigma^*$ is $+1$ (\cite{Al}, t. 3, p.11--17). 
We denote by $K^*$ the complex of the star duals of the faces
of $K$. It is an abstract simplicial complex whose $k$-simplices are the star duals of dimension $k$, that is the star duals of $(m-k)$-simplices of $K$. Indeed, $K^*$
is $\partial (P^*)$.  Let $K_0$ be a closed subcomplex of $K$ and let $K_0^*$ be the 
subcomplex of the star duals of the faces of $K_0$. 

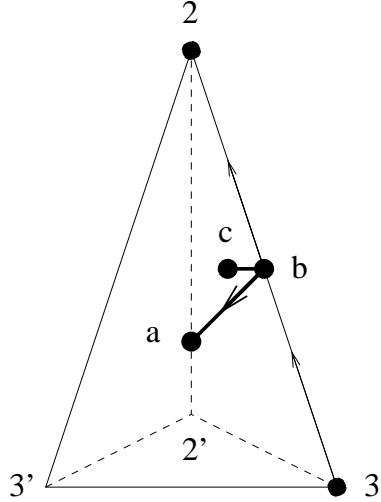
\hfil\scaledpicture 2.6in by 3.4in (fig15 scaled 700) \hfil
\vglue -.5mm
In the previous picture, let $\sigma$ denote the oriented edge $32$. We assume that the orientation of the tetrahedron $(233'2')$ is given by the standard orientation
of $\Bbb R^3$. Then the star dual of $\sigma$ is the sum of the oriented sum $cb+ba$ of the barycentric subdivision of $(233'2')$.

Let $k$ be a positive integer and let $[c]\in \tilde H_k(K_0,\Bbb Z)$ be a homology class represented by the cycle $c$. 
In $K$, the cycle $c$ is the boundary of a $(k+1)$-chain $d$. Decompose $d$ as 
$$
d=\sum a_i \sigma_i
\eqno {a_i\in\Bbb Z}
$$
where $\sigma_i$ are cells of $K$. We can assume that $a_i$ is zero if $\sigma_i$ is in $K_0$. 
Else, the boundaries of $d$ and of $d'$ where the sum appearing in $d$ is restricted to 
$K\setminus K_0$ differ from a boundary in $K_0$, hence both 
represent $[c]$.
 
Consider the star dual of $d$, that is the $(m-k-1)$-cochain
$$
d^*=\sum a_i \sigma_i^*
$$
Then $d^*$ is a coboundary in $K^*$ but only a cocycle in $K^*\setminus K_0^*$. The cohomology class of $d^*$ in
$\tilde H^{m-k-1}(K^*\setminus K^*_0,\Bbb Z)\simeq \tilde H^{m-k-1}(K\setminus K_0,\Bbb Z)$ is the Alexander dual of $[c]$.

Let us give an example. We use the previous picture. Let $K_0$ denote the subcomplex of $(233'2')$ constituted by the two edges $22'$ and $33'$. Then $(233'2')^*$ is a 
tetrahedron whose facets are $2$, $3$, $3'$ and $2'$ whereas $(233'2')^*\setminus K_0^*$ is this tetrahedron minus the four (open) facets and minus the two (open)
edges $22'$ and $33'$. The class of the $0$-cycle
$2-3$ is a generator of $\tilde H_0(K_0,\Bbb Z)$. In $(233'2')$, it is the boundary of the oriented edge $32$. The Alexander dual of $K_0$ in $(233'2')^*\setminus K_0^*$
is the oriented edge $32$ as shown in the following picture. It is a cocycle whose class generates $\tilde H^1((233'2')^*\setminus K_0^*,\Bbb Z)$. 
The picture represents the tetrahedron $(233'2')^*$. The subcomplex $(233'2')^*\setminus K_0^*$ is constituted 
by the bold edges. Finally, the orientation of the edge $32$ is given by the arrow. As before, the orientation of $(233'2')^*$ comes from the standard orientation 
of $\Bbb R^3$.

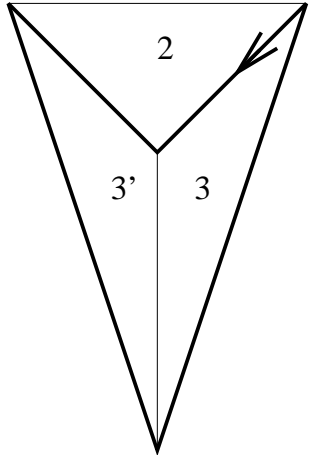
\hfil\scaledpicture 2.1in by 3.1in (fig16 scaled 700) \hfil

\remark{Remark 9.1}
The barycentric subdivision of $K^*\setminus K^*_0$ identifies naturally with the mirror complex of $K_0$. Via this identification, $d^*$ is a cochain of this mirror
complex. In the example given above, $d^*$ is then the cochain $cb+ba$ drawn in $(ABCD)_b$. Nevertheless, $d^*$ is {\it generally not a cocycle of the mirror complex}
of $K_0$, since the barycentric subdivision of a cocycle of a complex does not generally remain a cocycle in the barycentric subdivision of this complex.
\endremark

\demo{End of the proof of Theorem 7.6}  
Let $\Cal I $ be a proper subset of $\Cal F $, and $k$ an integer. 
Consider a class $[c]$ in $\tilde{H}_{k} (link_{\Delta } \sigma _{\Cal I} , \Bbb Z )$ 
represented by a simplicial $k$-cycle $c$. Then, in $S^{\bar{\Cal I }}$, the cycle $c$ is the boundary of 
some simplicial $(k+1)$-chain 
$$
d = \sum_{I \subset \bar{\Cal I }} a_I \sigma _I 
$$
As before, we assume that $a_I$ is zero if $I\in link_{\Delta } \sigma _{\Cal I}$, hence by Remark 7.11 we may keep only the sets $I$ 
which belong to $\tilde {\Cal I}$. Recall that $\sigma_I$ is oriented from the order of the poset $\Cal F$ (see Section 8).

\remark\nofrills{NB: \ } 
 in the sums we will use, we will only consider subsets that have a prescribed 
cardinal (for instance $k+2$ thereup). We will omit this precision in the 
sequel.
\endremark
\medskip

  For a simplex $\sigma _I $ in $S^{\bar{\Cal I }}$, we denote by $\sigma ^* _I $ its 
star dual in $(S^{\bar{\Cal I }})^* $. Then, the Alexander dual of $[c]$ in 
$\tilde H^{|\bar{\Cal I }|-k-3} ((S^{\bar{\Cal I }})^* \backslash (link _{\Delta } \sigma _{\Cal I })^*)$ is 
given by the class of 
$$
\sum_{I \in \tilde{\Cal I }} a_I (\sigma _I^* )\ .
$$

Indeed, $(S^{\bar{\Cal I }})^* \backslash (link _{\Delta } \sigma _{\Cal I })^*$ is isomorphic to $\tilde {\Cal I}$ and the previous cochain is a cocycle in 
$\tilde{\Cal I}$.

  Let us now place in $P _b $. Recall that $\tilde{\Cal I}$ is identified with $K_{\bar{\Cal I}}$ via the map 
$$
I\in \tilde{\Cal I}  \longmapsto \hat I \in K_{\bar{\Cal I}}
$$ 

Denote $F_{\hat I}^*$ the image of $\sigma_I^*$ via this map.
We now have to compute the Alexander dual in $\partial P$ of 
the cohomology class of
$$
\sum_{I \in \tilde{\Cal I }} a_I (F _{\hat I}^*)
$$
Obviously, the simplicial complex $K_{\bar{\Cal I}}$ is isomorphic to $\partial (P^*)\setminus P^*_{\Cal I}$. Via this identification, 
$F _{\hat I}^*$ is the star dual in $(\partial P)^*=\partial (P^*)$ of $F_{\hat I}$ in $P_{\Cal I}$.

Consider  
$$
\sum_{I \in \tilde{\Cal I }} a_I \partial F_{\langle \hat I\rangle }
\tag\numerote
$$
where the angles mean that the set $\hat I$ is ordered in a way which may be different from the natural order
induced by $\Cal F$; as explained in Section 8, the face $F_{\langle \hat I\rangle }$ is thus oriented. In the special case where $\hat I$
is a singleton, then $F_{\langle \hat I\rangle }$ may be $F_{\hat I}$ or $-F_{\hat I}$, that is the facet $\hat I$ 
with the orientation reversed.
It follows from the construction of Alexander duals recalled above that, if the order on each 
subset of $\tilde{\Cal I}$ is suitably chosen, 
then the former expression is a cycle whose class in $\tilde H_{k-|\bar{\Cal I}|+d+1} (P _{\Cal I }, \Bbb Z )$ 
is the searched Alexander dual.

Let us enlighten all this discussion with an example. Let $P$ be the cube numbered as in Example 7.9. Let $\Cal I=\{1,1'\}$. Then,
$$
link_{\Delta} \sigma_{\Cal I}=\{2, 2', 3, 3', 22',33'\}
$$

Let $c$ be the $0$-cycle $2-3$ in $link_{\Delta} \sigma_{\Cal I}$. We are exactly in the situation drawn in the two previous
pictures. We thus have that the Alexander dual of $[c]$ in $\tilde{\Cal I}$ is the oriented 1-cocycle $32$. Via the map recalled above, it corresponds to
the oriented 
$1$-cocycle $2'3'$ in $K_{\bar{\Cal I}}$. This cocycle is the star dual of the edge $2'3'$ of $P$. The boundary of this edge, that is $12'3'-1'2'3'$ is a $0$-cycle
whose class is a generator of $\tilde H_0(P_{\Cal I},\Bbb Z)$ as shown in the following picture.

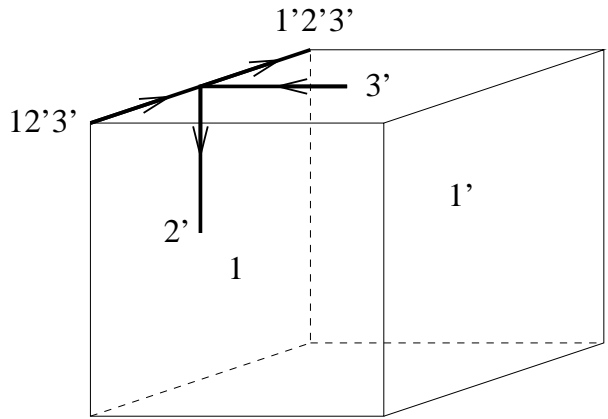
\hfil\scaledpicture 4.1in by 2.8in (fig17 scaled 700) \hfil

  The expression \lastnum[0] can be rewritten in another form using Lemma 8.6. This gives the following formula:
$$
\phi ([c]) = \left [
\sum_{H \in \Cal I , I \in \tilde{\Cal I }, H \cap F_{\hat I} \neq \emptyset} 
a_I F_{\langle \hat I\rangle}\cap H\right ]
$$

  Let us now prove that the cup product operation on the cohomology of $X$ corresponds (up 
to sign) to the operation of intersection on our homology classes (in the nonzero case), 
i.e. $\phi ([\langle i_{\bar{\Cal J }} - i_{\bar{\Cal I }}\rangle  * c * c'] ) = \pm \phi ([c]) \cap \phi ([c'])$, where
$i_{\bar{\Cal J }}$ (respectively $i_{\bar{\Cal I }}$) is an element of $\bar{\Cal J }$ (respectively $\bar{\Cal I }$).

  Consider two subsets $\Cal I $ and $\Cal J $ of $\Cal F $. If $\Cal I \cup \Cal J $ is not equal to 
$\Cal F $, then the cup product of classes associated to homology elements of $P _{\Cal I }$ 
and $P _{\Cal J }$ is zero as it corresponds to the case $\sigma \cup \sigma ' \neq [n]$ 
in \cite{DL}, Theorem~1.1. In the sequel, we assume that $\Cal I \cup \Cal J = \Cal F $.

  If we take $\Cal I $ equal to $\Cal F $, then only $\tilde{H}_{d-1} (P _{\Cal I }, \Bbb Z )$ 
is nonzero and a class $[c]$ in it is a multiple of the top-class of $\partial P$. 
Moreover, $\psi ([c])$ is in $H^0 (X,\Bbb Z )$, hence is a multiple (the same up to sign) 
of the unity of the cohomology ring of $X$. Therefore, both the intersection with $[c]$ and 
the cup product with $\psi ([c])$ are, up to sign, multiplication by this integer. This 
proves the formula in the particular case $\Cal I = \Cal F $ ($\Cal J = \Cal F $ is identical).

  From now on, we assume that $\Cal I $ and $\Cal J $ are distinct from $\Cal F $ (in particular 
they are nonempty as well). As we are working up to sign, we can assume: 
\medskip

\noindent $\underline{\text{Hypothesis :}}$ for the order on the facets, any element of $\bar{\Cal I }$ is less than 
any element of $\bar{\Cal J }$, i.e. $\bar{\Cal I}\cup\bar{\Cal J}=\bar{\Cal I}\bar{\Cal J}$ as ordered sets.
\medskip

  We thus consider an element $[c_{\Cal I} ]$ of 
$\tilde{H}_{k} (link_{\Delta } \sigma _{\Cal I } , \Bbb Z )$ and an element $[c_{\Cal J} ]$ of
$\tilde{H}_{k'} (link_{\Delta } \sigma _{\Cal J } , \Bbb Z )$. Let $[c]$ be 
$[c_{\Cal I} * c_{\Cal J} * \langle i_{\bar{\Cal J}} -i_{\bar{\Cal I}} \rangle  ]$ in 
$\tilde{H}_{k+k'} (link_{\Delta } \sigma _{\Cal I \cap \Cal J } , \Bbb Z )$. We have to see 
that $\phi ([c])$ is, up to sign, the intersection of $\phi ([c_{\Cal I}] )$ with $\phi ([c_{\Cal J}] )$.
Let $d_{\Cal I}$ (respectively $d_{\Cal J}$) be a $(k+1)$-chain of $S^{\bar {\Cal I}}$ (respectively a $(k+1)$-chain of $S^{\bar {\Cal J}}$) 
whose boundary has $[c_{\Cal I}]$ (respectively $[c_{\Cal J}]$) for class.

First, we find a chain in $S^{\bar {\Cal I}\cup\bar{\Cal J}}$ whose boundary has $[c]$ for class. 

\proclaim{Lemma 9.2}
  Consider two disjoint nonempty finite sets $A$ and $B$. Consider a $k$-chain $d_A $ and 
a $k'$-chain $d_B $ in subcomplexes $K_A $ and $K_B $ of $\Delta _A $ and $\Delta _B $. 
Then, up to sign, $\partial (d_A * d_B )$ is homologous to 
$\partial d_A  * \partial d_B * \langle i_A - i_B \rangle $ in 
$K_A * \Delta _B \bigcup \Delta _A * K_B $, where $i_A $ and $i_B $ denote elements of 
$A$ and $B$ respectively.
\endproclaim

\demo{Proof}
  In fact, we show that $\partial (d_A * d_B )$ is homologous to 
$\partial d_A * \langle i_B - i_A \rangle  * \partial d_B $ which is clearly equal up to sign to 
$\partial d_A  * \partial d_B * \langle i_A - i_B \rangle $.

We have $\partial (d_A * d_B ) = \partial d_A * d_B + (-1)^{k+1} d_A * \partial d_B $. We 
then just have to see that $\partial d_A * d_B $ and 
$\partial d_A * \langle i_B \rangle  * \partial d_B $ differ from a boundary and that 
$(-1)^{k+1} d_A * \partial d_B $ and $\partial d_A * \langle -i_A \rangle  * \partial d_B $ do too.

  The boundary of $\langle i_B \rangle  * \partial d_B $ is $\partial d_B $. Hence, $d_B $ and 
$\langle i_B \rangle  * \partial d_B $ differ from a cycle and this cycle is a boundary in 
$\Delta _B $ as it is not $0$-dimensional. This gives immediately that 
$\partial d_A * d_B $ and $\partial d_A * \langle i_B \rangle  * \partial d_B $
differ from a boundary in $K_A * \Delta _B $.

  We have $\partial d_A * (-\langle i_A \rangle ) = (-1)^{k+1} \langle i_A \rangle  * \partial d_A $ and, as above, 
$\partial d_A * \langle -i_A \rangle  * \partial d_B $ and $(-1)^{k+1} d_A * \partial d_B $ differ from a 
boundary in $\Delta _A * K_B $.

  This proves the lemma.
$\square$
\enddemo

  In our context, 
the lemma shows that we can take $d_{\Cal I} * d_{\Cal J} $ as chain having the desired boundary.

  We then can compute $\phi ([c])$. Suppose we have 
$$
d_{\Cal I} = \displaystyle{\sum_{I \in \tilde{\Cal I }} a_I \sigma _I } \quad\text{and}\quad 
d_{\Cal J} = \displaystyle{\sum_{J \in \tilde{\Cal J }} b_J \sigma _J }
$$ 
Then, as $\bar{\Cal I }$ and $\bar{\Cal J }$ are disjoint, we have thanks to the chosen order:
$$
d_{\Cal I} * d_{\Cal J} = \sum_{I \in \tilde{\Cal I }, {J \in \tilde{\Cal J }}} 
a_I b_J \sigma _{I \cup J}
$$
In fact, as noted at the beginning of this Section, we may replace $d_{\Cal I}*d_{\Cal J}$ by a homologous chain in 
$S^{\bar {\Cal I}\cup\bar{\Cal J}}\setminus link_{\Delta}(\sigma_{\Cal I\cap \Cal J})$, by keeping in the former equation only the 
couples $(I,J)$ such that $I\cup J$ is in $\widetilde{\Cal I\cap\Cal J}$.

The following lemma ensures us that the intersection of two cycles is again a cycle.

\proclaim{Lemma 9.3}
  Up to a sign that is independent of $I$ and $J$ (but which depends on their cardinal), we get 
$F_{\langle \hat I \cup \hat J\rangle } = F_{\langle \hat I\rangle  + \langle \hat J\rangle }$ (when these intersections are nonempty).
\endproclaim

\demo{Proof}
  Both members of the equality represent 
$F_{\hat I\cup\hat J}$, hence are equal up to sign.
To compute this sign, we have to understand which is the orientation of $F_{\langle \hat I\rangle }$ knowing 
$I$.

 The polytope $P$, the simplicial sphere $S^{\bar{\Cal I}}$ and $\sigma_I$ are oriented. As explained
above, this induces an orientation of the star dual $\sigma_I^*$. Via the isomorphism of complexes given by
$I \to \hat I$, an orientation is fixed on the image $F^*_{\hat I}$ of $\sigma_I^*$. We order $\langle\hat I\rangle$ so that the star dual of 
$F_{\langle \hat I\rangle }$ is $F^*_{\hat I}$.

  Note $\epsilon_{\langle \hat I\rangle }$ being $+1$ if $F_{\hat I}$ is oriented like 
$F_{\langle \hat I\rangle }$ and $-1$ else. We want to show that 
$\epsilon_{\langle \hat I\rangle } \cdot \epsilon_{\langle \hat J\rangle } \cdot 
\epsilon_{\langle \hat I \cup \hat J\rangle}$ neither depends on $I$ nor on $J$.  We will in fact prove more than stated in the lemma, since we will give the exact sign of this product. This will be useful in the next Section. 

The (unoriented) 
star dual of $\sigma_I$ in $(S^{\bar{\Cal I}})_b\setminus (link_{\Delta}\sigma_{\Cal I})_b\simeq \tilde{\Cal I}_b$ is in fact the order complex $C_I$ on the
sets $I'$ such that $I\subset I' \subsetneq \bar{\Cal I}$. Under the identification between the mirror complex of $link_{\Delta}
\sigma_{\Cal I}$ and $(K_{\bar{\Cal I}})_b$, the complex $C_I$ may also be seen as the (unoriented) star dual of
$F_{\langle\hat I\rangle}$. It is easy to check that, due to the simpleness of $P$, it is besides isomorphic to the barycentric 
subdivision of the simplex $\Delta_{\hat I}$ with vertex set $\hat I$. To simplify the
proof, we will, by abuse of notation, call $C_I$ these three complexes.

As a consequence of the identification between $C_I$ and $\Delta_{\hat I}$, an ordering of $\langle \hat I\rangle$ 
induces an orientation of $C_I$. The star dual orientation is the one for which
the intersection number $\sigma_I\times C_I$ is $1$. 
On the other hand, when we see $C_I$ as a subcomplex of $\partial P_b$, 
the star dual orientation is the one for which the intersection number $F_{\langle\hat I\rangle}\times C_I$ is $1$. 
In particular, for any orientation of $C_I$, these two intersection numbers are the same. 

Put on $C_I$ the orientation given by the
natural order of $\hat I$ as subset of $\Cal F$. 
Remark that the sets $\hat I I$ and $\bar{\Cal I}$ are the same
up to a permutation. Let $\epsilon_{\hat I I}$ denote the sign of this permutation.

We claim that, with the orientation we have fixed on $C_I$, we have:
$$
\sigma_I\times C_I=F_{\langle\hat I\rangle}\times C_I=(-1)^{(\vert \hat I\vert -1)(\vert I\vert -1)}\cdot\epsilon_{\hat I I}
\tag\numerote
$$
This can be shown as follows. We still identify $C_I$ with $(\Delta_{\hat I})_b$. 
Let $\hat I=\{\hat\imath_0<\hdots < \hat\imath_l\}$. Consider the positively oriented simplex $\sigma=J_0<\hdots <J_l$ 
of $(\Delta_{\hat I})_b$ defined by $J_s=\hat\imath_0\hdots \hat\imath_s$.

On the other hand, consider the oriented barycentric subdivision $(S^{\bar{\Cal I}})_b$. 
Consider $(\sigma_I)_b$. Let $I=\{i_0<\hdots <i_{k+1}\}$ and consider the positively oriented simplex
$\sigma'=I_0<\hdots <I_{k+1}$ defined by $I_s=i_s \hdots i_{k+1}$.

In the sphere $(S^{\bar{\Cal I}})_b$, the simplex $\sigma$ corresponds in fact via the map $I\mapsto \hat I$ to the simplex
$\bar{\Cal I}\setminus J_0 <\hdots <\bar{\Cal I}\setminus J_l$. Notice that 
$$
\bar{\Cal I}\setminus J_l=\bar{\Cal I}\setminus \hat I=I=I_0
$$
so we may consider the simplex
$$
\bar{\Cal I}\setminus J_0 <\hdots <\bar{\Cal I}\setminus J_l=I_0<\hdots <I_{k+1}
$$
It is easy to check that this simplex induces
the orientation $\hat I I$ on $S^{\bar{\Cal I}}$.

Consider now the ``reversed simplex''
$$
I_{k+1}<\hdots< I_0=\bar{\Cal I}\setminus J_l <\hdots< \bar{\Cal I}\setminus J_0
$$
It induces on $S^{\bar{\Cal I}}$ an orientation which differs from the previous one and is equal to 
$$
\epsilon=(-1)^{(k+l+1)(k+l+2)/2}\cdot\epsilon_{\hat I I}\ .
$$
In the same way, the simplex $I_{k+1}<\hdots< I_0$ of $(\sigma_I)_b$ is no more positive but with sign
$$
\epsilon'=(-1)^{(k+1)(k+2)/2}
$$
and the simplex $\bar{\Cal I}\setminus J_l <\hdots< \bar{\Cal I}\setminus J_0$ is no more positive
but with sign
$$
\epsilon''=(-1)^{l(l+1)/2}
$$
By \cite{Al}, t. 3, p.11--17, the intersection number $\sigma_I\times C_I$ is given by the product 
$\epsilon\cdot \epsilon'\cdot \epsilon''$.
 A direct computation shows now the claim.

Of course, putting the natural orders on $\hat J$ and $\widehat {I\cup J}$, the same argument implies that
$$
\eqalign{
&\sigma _J\times C_J = F_{\langle \hat J\rangle }\times C_J = (-1)^{(\vert \hat J\vert-1) (\vert J\vert -1)}\cdot\epsilon_{\hat{J}J}
\cr 
&\sigma_{I \cup J}\times C_{I \cup J}= (-1)^{(\vert \hat I\cup \hat J\vert-1) (\vert I\cup J\vert -1)}\cdot
\epsilon_{\widehat{I \cup J}I\cup J} = 
F_{\langle \widehat{I \cup J}\rangle }\times C_{I \cup J}
}
\tag\numerote
$$

Let us consider now the situation on $\partial P$. By definition, we have
$$
F_{\langle \hat I\rangle } \times F_{\hat I }^* = 1 \quad\text{and}\quad F_{\hat I }\times F_{\hat I }^* = \epsilon_{\langle \hat I\rangle }
$$

  Let now $(H_0 ,..., H_l )$ be a collection of hyperplanes supporting facets of $P$ such that $F_{\hat{I}}$ is equal
to $F_{(H_0 ,..., H_l )}$. The two previous intersection numbers can be interpreted as follows.
Let $\langle B \rangle $ and $B^* $ 
be respective positive basis of $F_{\langle \hat I\rangle }$ and $F_{\hat I }^* $ (more 
exactly of its part lying on $H_0$). Let $B^+_{F_{\hat{I}}}$ be a positive 
basis of the face $F_{\hat{I}}$. Finally, let $v_i$ denote an outward pointing normal vector to $H_i$ for $i$ between $0$ and
$l$. Then 
the basis $(v_0, \langle B \rangle , B^* )$ of $\R^d$ is direct, whereas 
$(v_0, B^+_{F_{\hat{I}}} , B^* )$ is a basis of $\R^d$ whose sign is $\epsilon_{\langle \hat I\rangle }$. On the other hand,
we have that $(v_1,\hdots,v_l)$ is a direct basis of $C_I$, therefore the sign of the permutation transforming $(v_1,\hdots,v_l)$
into $B^*$ is equal to the intersection number $F_{\langle\hat I\rangle}\times C_I$.

With our conventions, to say that $B^+_{F_{\hat{I}}}$ is positive means exactly that the basis $(v_0,\hdots,v_l, B^+_{F_{\hat{I}}})$
is direct. The sign of $(v_0, B^+_{F_{\hat{I}}} , B^* )$ is also given as the sign of the transformation sending it to
$(v_0,\hdots,v_l, B^+_{F_{\hat{I}}})$, or to the product of the sign of the transformation sending it to 
$(v_0,B^*,B^+_{F_{\hat{I}}})$ by the sign of the transformation sending $(v_0,B^*,B^+_{F_{\hat{I}}})$ to $(v_0,\hdots,v_l, B^+_{F_{\hat{I}}})$. By what preceeds, this last sign is equal to the intersection number $F_{\langle\hat I\rangle}\times C_I$.

As a consequence of all this and of \lastnum[-1], we obtain the following identity
$$
\eqalign{
\epsilon_{\langle \hat I\rangle } &= (-1)^{(\vert \hat I\vert-1) (\vert I\vert -1)}\cdot
\epsilon_{\hat{I}I} \cdot (-1)^{dim F_{\hat{I}} \cdot l} \cr
&= 
(-1)^{(\vert \hat I\vert-1) (\vert I\vert -1)}\cdot\epsilon_{\hat{I}I} \cdot (-1)^{d - |\hat{I}| \cdot (|\hat{I}| - 1)}\cr
&= (-1)^{(\vert \hat I\vert-1) (\vert I\vert -1)}\cdot\epsilon_{\hat{I}I} \cdot (-1)^{d}
}
\tag\numerote
$$

  The previous equality is naturally also true for $J$ and $I \cup J$. As a consequence of \lastnum[-2], \lastnum[-1], \lastnum[0] and of the hypothesis made above, we have, after computation,
$$
\epsilon_{\langle \hat I\rangle } \epsilon_{\langle \hat J\rangle } 
\epsilon_{\langle \hat I \cup \hat J\rangle } = (-1)^{1+\vert \hat I\vert \cdot\vert J\vert +\vert \hat J\vert \cdot\vert I \vert}\cdot
\epsilon_{\hat{I}I} \epsilon_{\hat{J}J} \epsilon_{\hat{I}\hat{J}IJ} \cdot 
(-1)^d
$$
  Now, the product $\epsilon_{\hat{I}I} \epsilon_{\hat{J}J}$ sends the ordered set 
$\bar{\Cal I }\cup \bar{\Cal J }$ to $\hat{I}I \hat{J}J$ and $\epsilon_{\hat{I}\hat{J}IJ}$ 
sends this same ordered set to $\hat{I}\hat{J}IJ$. Their product is the sign of the 
permutation which permutes $I$ and $\hat{J}$, hence is equal to 
$(-1)^{|I| \cdot |\hat{J}|}$.

  This finally gives 
$$
\epsilon_{\langle \hat I\rangle } \epsilon_{\langle \hat J\rangle } 
\epsilon_{\langle \hat I \cup \hat J\rangle } = (-1)^{1+\vert \hat I\vert \cdot\vert J\vert +\vert \hat J\vert \cdot\vert I \vert+
d + |I| \cdot |\hat{J}|}=(-1)^{d+1+\vert \hat I\vert \cdot\vert J\vert},
$$
a sign which is independent of $I$ and $J$.
$\square$
\enddemo

  Thanks to this lemma, we can claim that, up to sign
$$
\phi ([c]) = \left [\sum_{I \in \tilde{\Cal I }, J \in \tilde{\Cal J }} a_I b_J \partial 
F_{\langle \hat I\rangle  + \langle \hat J\rangle }
\right ]
$$

By Lemma 8.6, this gives us then, up to sign
$$
\phi ([c]) = \left [\sum_{I \in \tilde{\Cal I }, J \in \tilde{\Cal J }, H \in \Cal I \cap \Cal J} 
a_I b_J F_{\langle \hat I\rangle  + \langle \hat J\rangle}  \cap H
\right ]
$$

  On the other side, we have to compute the intersection of $\phi ([c_{\Cal I}] )$ and 
$\phi ([c_{\Cal J}] )$. Let us write them
$$
\displaylines{
\phi ([c_{\Cal I}] ) = \left [
\sum_{H \in \Cal I ; I \in \tilde{\Cal I} ; F_{\langle \hat I\rangle} \cap H\neq \emptyset } 
a_I F_{\langle \hat I\rangle}\cap H 
\right ]\cr
\phi ([c_{\Cal J}] ) = \left [
\sum_{H' \in \Cal J ; J \in \tilde{\Cal J} ; F_{\langle \hat J\rangle}\cap H' \neq \emptyset } 
b_{J} F_{\langle \hat J\rangle}\cap H'
\right ]
}
$$
These two classes are naturally realised in the boundaries of $\Cal F _{\Cal I }$ and $\Cal F _{\Cal J }$ 
but do not then meet transversely. We can nevertheless "push" them in the interior of 
these sets so that they do.

\definition{Definition 9.4}
  Consider a simple polytope $P$ and for each of its facet $H$ an affine function 
$l_H $ on the space of $P $ which is zero on $H$ and positive on $P \backslash H$. 
For $\epsilon > 0$, call $H_{\epsilon } = l_H ^{-1} (\epsilon ) \cap P $ and for a face 
$F$ of $P $, note $F_{\epsilon} = \displaystyle{\bigcap_{H \supset F} H_{\epsilon }}$.
\enddefinition

\proclaim{Lemma 9.5}
  Consider now two faces $F$ and $F'$ of a simple polytope $P $ that are not contained in a 
common facet and have nonempty intersection. Then, if $\epsilon $ is small enough, 
$\partial F_{\epsilon }$ and $\partial F'_{\epsilon }$ meet transversely and their 
intersection is $\partial (F \cap F')_{\epsilon}$. Moreover, this also works when we deal 
with oriented faces.
\endproclaim

  This lemma is clear.

  We now can compute the homology class of the intersection of our two cycles. For this, 
consider for every facet of $P $ an affine function satisfying the properties Definition 9.4.

  Take $\epsilon > 0$ small enough. Define $\phi_{\epsilon }([c_{\Cal I}] )$ as follows : 
for an element $I$ of $\tilde {\Cal I }$ and a facet $H$ of $\Cal I $ meeting $F_I$, call 
$(F_I\cap H)_{H, \epsilon }$ the set $(F_I\cap H )_{\epsilon }$ when we consider $H$ as a simple polytope and restrict the affine 
functions of the facets meeting 
$H$ to the facets of $H$. Just remark now that
$$
\phi ([c_{\Cal I}] ) = \left [
\sum_{H \in \Cal I ; I \in \tilde{\Cal I} ; F_{\langle \hat I\rangle}\cap H \neq \emptyset } 
a_I (F_{\langle \hat I\rangle}\cap H )_{H, \epsilon }
\right ]
$$
since the cycle in the brackets thereup is homotopic to $\sum_{H \in \Cal I ; I \in \tilde{\Cal I} ; F_{\langle \hat I\rangle} \cap H\neq \emptyset } 
a_I F_{\langle \hat I\rangle}\cap H $.

Of course, the same is true for $\phi([c_{\Cal J}])$. But these cycles meet transversely and, 
thanks to Lemma 9.5, their intersection can be written:
$$
\phi ([c_{\Cal I}] ) \cap \phi([c_{\Cal J}] ) =
\left [ 
\sum_{H \in \Cal I ; I \in \tilde{\Cal I } ; J \in \tilde{\Cal J } F_{\langle \hat I\rangle} \cap F_{\langle \hat J\rangle}\cap H \neq 
\emptyset } a_I b_J (F_{\langle \hat I\rangle} \cap F_{\langle \hat I\rangle}\cap H)_{H, \epsilon }
\right ]
$$
And this last expression is then $\phi([c_{\Cal I} \cap c_{\Cal J}] )$.

  We get finally, up to sign:
$$
\phi ([c_{\Cal I}] ) \cap \phi ([c_{\Cal J}] ) =\left [ 
\sum_{I \in \tilde{\Cal I }, J \in \tilde{\Cal J }, H \in \Cal I \cap \Cal J} 
a_I b_J F_{\langle \hat I\rangle  + \langle \hat J\rangle}  \cap H 
\right ]
= \phi ([c])
$$
  This completes the demonstration of the Theorem.
$\square$
\enddemo

\head 
{\bf 10. Computation of the sign}
\endhead
  
In the previous Section, the product of two generators of the cohomology of $X$ was 
computed up to sign. Here do we compute it precisely. This gives:

\proclaim{Sign Theorem 10.1}
  Consider $[c] \in \tilde H_k (P _{\Cal I } ,\Bbb Z )$ and $[c'] \in \tilde H_{k'} (P _{\Cal J },\Bbb Z )$ as in
the statement of Cohomology Theorem 7.6. Set  $K'=|\bar{\Cal J}|-d+k'-1$. Denote by
$\epsilon_{\bar{\Cal I}\bar{\Cal J}}$ the sign of the permutation transforming $\bar{\Cal I}\bar{\Cal J}$ into
$\bar{\Cal I }\cup\bar{\Cal J}$. Then, 
$$
\psi ([c]) \smile \psi ([c']) = \epsilon \psi ([c] \cap [c'])
$$ 
with 
\setbox1=\hbox{$\epsilon_{\bar{\Cal I}\bar{\Cal J}}\cdot (-1)^{(d+1+n + K'|\bar{\Cal I}|)}$}
\setbox2=\hbox to \wd1{$\hfill 1 \hfill$}
$$
\epsilon = \left\{ \eqalign{
& \box1\text{ if neither }\bar{\Cal I }\text{ nor }\bar{\Cal J } \text{ is empty.} \cr
& \box2 \text{ if at least one is.}
}\right .
$$
\endproclaim

\demo{Proof}
  In the special case where $\Cal I=\Cal F$, the class $[c]$ is a multiple of the top class of $\partial P$
and $\psi([c])$ is a multiple of the unity of the cohomology ring of $X$. The intersection of $[c]$ with any class
$[c']$ and the cup product  of $\psi([c])$ with $\psi([c'])$ are just multiplications by integers and $\epsilon$ is $1$.

  For the general case, let $[c] \in \tilde H_k (P _{\Cal I } ,\Bbb Z )$ and $[c'] \in \tilde H_{k'} (P _{\Cal J },\Bbb Z )$.
Due to Lemma 7.12, they correspond to classes $[c_1] \in\tilde{H}_K (link_{\Delta } \sigma _{\Cal I} , \Bbb Z )$ and 
$[c_2] \in\tilde{H}_{K'} (link_{\Delta } \sigma _{\Cal J} , \Bbb Z )$ with
$$
K=|\bar{\Cal I}|-d+k-1 \quad\text{and}\quad  K'=|\bar{\Cal J}|-d+k'-1
$$ 

Let us recall now de Longueville's results. 
The cup product of these two classes is the class of 
$(-1)^{n+K(K'+1)+1} \langle i_{\bar{\Cal J }} - i_{\bar{\Cal I} }\rangle  * c_1 * c_2 $ in 
$\tilde{H}_{K+K'+2} (link_{\Delta } (\sigma _{\Cal I} \cap \sigma _{\Cal J} ), \Bbb Z )$.
Due to the proof of Lemma 9.2, if we take the class associated to 
the boundary of $d_{\Cal I} * d_{\Cal J} $ instead of $\langle i_{\bar{\Cal J} } - i_{\bar{\Cal I} }\rangle  * c_1 * c_2$, 
the sign is $(-1)^{n+KK'}$.

When passing to the classes in $\partial P$, a sign comes: it is explicitely described in the proof
of Lemma 9.3 and 
is equal to 
$$
(-1)^{d+1+(K'+2)(|\bar{\Cal I}|-K-2)}=(-1)^{d+1+K'(|\bar{\Cal I}|-K)}
$$ 
under the hypothesis that in the order of $\Cal F$, the elements of $\bar{\Cal I}$ are lower than the elements of $\bar{\Cal J}$. 
There exists a permutation which 
reorders $\bar{\Cal I }\cup\bar{\Cal J}$ 
such as this assumption holds and we thus have to multiply the result by $\epsilon_{\bar{\Cal I}\bar{\Cal J}}$, the sign of this permutation.

Putting all these results together gives the formula of Sign Theorem 10.1.
$\square$
\enddemo

\head
{\bf 11. Applications to the topology of the links}
\endhead

In this Section we make use of the previous results on the cohomology ring of a 2-connected link $X$ to investigate how
complicated can the topology of a link be. We will see that the complexity increases when the dimension $d$ of the associate
polytope $P$ increases and that the topology of a link may finally be ``arbitrarily complicated''.

For $d=0$, the unique 2-connected link is a point, for $d=1$ it is $\Bbb S^3$ (this is the case $p=0$ and $n=2$). For the polygons,
the situation is not so easy and the links are products of odd-dimensional spheres or connected sums of products of spheres: this case
was completely described in \cite{McG} (cf Theorem 6.3). In higher dimensions, the only known case is the special case where
$p=2$ \cite{LdM1}, \cite{LdM2} where the same type of manifolds is obtained (cf Example 0.5). On the other hand, the generalization of
MacGavran's results stated as Theorem 6.3 shows that, for any value of $d$, there is an infinite number of examples where the link
is a connected sum of products of spheres. This leads naturally to the following question, whose positive answer was stated as a
conjecture in \cite{Me1}

\proclaim{Question A}
Is it true that any 2-connected link may be decomposed into a product of odd-dimensional spheres and connected sums of products
of spheres?
\endproclaim

A weaker version of this question is

\proclaim{Question A'}
At least, is it true that the cohomology ring of a 2-connected link is isomorphic to the cohomology ring of  
 a product of odd-dimensional spheres and connected sums of products
of spheres?
\endproclaim

This supposes to resolve first the (easier?)

\proclaim{Question A''}
Is it true that the homology of a 2-connected link is always without any torsion?
\endproclaim

An immediate application of Cohomology Theorem 7.6 is that the answer is yes if $d$ is lower than $4$.

\proclaim{Corollary 11.1}
  If the polytope $P $ has dimension at most $4$, then the homology of the associated 
manifold is torsion free.
\endproclaim

\demo{Proof}
  In this case, every homology group of the form $\tilde{H}_k (P _{\Cal I },\Bbb Z)$ is
torsion free, as $P _{\Cal I }$ lies in $\partial P$ which is a sphere of 
dimension $\leq 3$ (see \cite{Al}, t. 3, Chapter XIII, paragraph 4.12). So is a direct sum of such groups as are the cohomology groups of $X$ by Cohomology Theorem 7.6.
$\square$
\enddemo

We emphasize that this result obtained as an easy consequence of Cohomology Theorem 7.6 should not be easily deduced from the classical form of the Goresky-Mac Pherson
formula (for example in the version of \cite{DL}) applied to the complement of subspace arrangement $\Cal S$, since the dimension of the complex $\Delta$ on which the
homology computations have to be done can be much greater than $3$. Therefore, this Corollary illustrates all the interest in having a formula in terms of subsets of the
associate polytope.

We will now prove that, even in dimension $3$, the answer to questions A and A' is negative. To see this, we will first compute how the cohomology of a link $X$ changes when
performing an elementary surgery of type $(1,n)$ on $X\times \Bbb S^1$, that is when performing a vertex cutting on $P$. Recall that, by Lemma 6.1, the diffeomorphism
type of the new link $X'$ is independent of the choice of the vertex to be cut off.

\proclaim{Proposition 11.2}
  Let $X$ and $X'$ as above. Assume that $d\geq 2$. Then: 
$$
\eqalign{
&H^0 (X', \Bbb Z ) \simeq H^{n+d+1}(X', \Bbb Z ) \simeq \Bbb Z \cr
&H^1 (X', \Bbb Z ) \simeq H^2 (X', \Bbb Z ) \simeq H^{n+d-1} (X', \Bbb Z ) \simeq 
H^{n+d} (X', \Bbb Z ) \simeq 0\cr
&H^i (X', \Bbb Z ) \simeq H^i (X, \Bbb Z ) \oplus H^{i-1} (X, \Bbb Z ) \oplus 
\Bbb Z ^{\left (\smallmatrix n-d \\ i-2d+1 \endsmallmatrix\right )} \oplus \Bbb Z^{\left (\smallmatrix n-d \\ i-2 \endsmallmatrix\right )} \text{ if }3 \leq i \leq n+d-2}
$$
where $\pmatrix l \\ k \endpmatrix$ is 
zero if $k<0$ or $k>l$.

  Moreover, the product is given by the following rules considering two cohomology 
classes $[c]$ and $[c']$ of $X'$: 

\smallskip
\noindent $\underline{\text{Rule 1:}}$  if $[c]$ or $[c']$ is in $H^0 (X', \Bbb Z )$ or 
$H^{n+d+1} (X', \Bbb Z )$, then the product is the obvious one.
\smallskip
  
Assume this is not the case. Then note $S_{i, j}$ for $3 \leq i \leq n+d-2$ and $1 \leq j \leq 4$, 
the sums thereup when they exist, that is 
$$
H^i(X',\Bbb Z)=S_{i,1}\oplus S_{i,2}\oplus S_{i,3}\oplus S_{i,4}\ .
$$
 
For $j=1$ or $j=2$, decompose $S_{i,j}$ as $\oplus_{\Cal I\subset\Cal F} S_{\Cal I, j}$ as in Cohomology Theorem 7.6. 
Finally denote by 
$S_j$, for $1 \leq j \leq 4$ the sums of $S_{i,j}$ when $i$ varies. We assume that $[c]$ is 
in $S_{\Cal I ,j}$ and $[c']$ in $S_{\Cal J ,j'}$.

\smallskip
\noindent $\underline{\text{Rule 2 :}}$ if $\{ j,j' \} \neq \{ 1 \} , \{ 1,2 \} , \{ 3,4 \} $ then 
$[c] \smile [c'] = 0$.
\smallskip

  Call $\varphi _1 $ and $\varphi _2$ the applications of $H^i (X, \Bbb Z )$ in 
$S_{i,1}$ and $S_{i+1 , 2}$.

\smallskip
\noindent $\underline{\text{Rule 3 :}}$ if $j = j' = 1$, then we can assume that $[c] = \varphi _1 ([c_1])$ and $[c']= \varphi _1 ([c'_1])$. Then 
$[c] \smile [c'] = - \varphi _1 ([c_1] \smile [c'_1])$.

\smallskip
\noindent $\underline{\text{Rule 4 :}}$ if $j = 1$ and $j' = 2$, then we can assume that 
$[c] = \varphi _1 ([c_1])$ and 
$[c']= \varphi _2 ([c'_2])$. 
Then 
$[c] \smile [c'] = - \varphi _2 ([c_1] \smile [c'_2])$.

\smallskip
\noindent $\underline{\text{Rule 5 :}}$ the cup product from $S_3 \times S_4 $ to 
$H^{n+d+1}(X, \Bbb Z ) \simeq \Bbb Z $ is a unimodular bilinear form, which is diagonal 
in the canonical basis (when these basis are suitably ordered). Note that the product 
vanishes when dimensions do not correspond.
\endproclaim

In particular, if the cohomology of $X$ has no torsion, then so has the cohomology of $X'$.

\remark{Remark 11.3}
  The isomorphisms are not completely canonical. Some judicious choices have 
to be made to obtain the desired rules about the cup product.
\endremark
\medskip

\demo{Proof}
  Let $v$ be the cut vertex, $\Cal F _v $ the set of the facets of $P $ that contain 
$v$ and $F$ the "new" facet (we will not distinguish a facet of $P $ -even in 
$\Cal F _v $- from the "same" facet of $P '$).

\remark{Notation 11.4}
  For a subset $\Cal I $ of $\Cal F $, we will denote $\Cal I _2 $ the subset of the 
facets of $P '$ having the same elements as $\Cal I $ and $\Cal I _1 $ the subset of 
the facets of $P '$ where we add $F$ to the ones of $\Cal I $.
\endremark
\medskip

Let $\Cal I\subset\Cal F$ such that the intersection of $\Cal I$ with $\Cal F _v $ is proper and nonempty; then $v$ belongs to the topological
boundary of $P_{\Cal I}$ and both $P '_{\Cal I _1 }$ and 
$P '_{\Cal I _2 }$ are homotopy equivalent to $P _{\Cal I }$. Therefore, the three sets 
have the same reduced homology groups.

  Consider now a subset $\Cal I $ of $\Cal F $ that contains $\Cal F _v $. Then 
$P '_ {\Cal I _1 }$ is homotopy equivalent to $P _{\Cal I }$, hence has the same reduced 
homology groups and $P '_{\Cal I _2 }$ is homotopy equivalent to $P _{\Cal I }$ minus a 
point. Therefore, if $\Cal I \neq \Cal F $, then the reduced homology groups of 
$P '_ {\Cal I _2 }$ are isomorphic to the ones of $P _{\Cal I }$ except 
$\tilde{H}_{d-2} (P '_{\Cal I _2 }, \Bbb Z )$ which is isomorphic to 
$\tilde{H}_{d-2} (P _{\Cal I }, \Bbb Z ) \oplus \Bbb Z $. And if $\Cal I = \Cal F $, then 
$P '_{\Cal I_2}$ is contractible, hence has no reduced homology.

  Consider now a subset $\Cal I $ of $\Cal F $ that is disjoint from $\Cal F _v $. Then 
$P '_{\Cal I _2 }$ is homotopy equivalent to $P _{\Cal I }$, hence has the same reduced 
homology groups and $P '_{\Cal I _1 }$ is homotopy equivalent to the disjoint union of 
$P _{\Cal I }$ with a point. Therefore, if $\Cal I \neq \emptyset $, then the reduced 
homology groups of $P '_ {\Cal I _1 }$ are isomorphic to the ones of $P _{\Cal I }$ except 
$\tilde{H}_0 (P '_{\Cal I _2 }, \Bbb Z )$ which is isomorphic to 
$\tilde{H}_0 (P _{\Cal I }, \Bbb Z ) \oplus \Bbb Z $. And if $\Cal I = \emptyset $, then 
$P '_{\{ F \} } = F$ is contractible and has no reduced homology.

  Let $i$ be an integer. 
Then, the above results allow us to compute $H^i (X', \Bbb Z )$. 
This gives:
$$
\eqalign{
H^i (X', \Bbb Z ) &\simeq \bigoplus_{\Cal I \subset \Cal F } 
\tilde{H}_{d + |\bar {\Cal I} _1 | - i - 1} (P' _{\Cal I _1 }, \Bbb Z ) 
\bigoplus_{\Cal I \subset \Cal F } \tilde{H}_{d + |\bar {\Cal I} _2 | - i - 1} 
(P' _{\Cal I _2 }, \Bbb Z ) \cr
&\simeq\bigoplus_{\Cal I \subset \Cal F } \tilde{H}_{d + |\bar{\Cal I }| - i - 1} 
(P' _{\Cal I _1 }, \Bbb Z ) 
\bigoplus_{\Cal I \subset \Cal F } \tilde{H}_{d + |\bar{\Cal I }| - i} 
(P' _{\Cal I _2 }, \Bbb Z )
}
$$
which is isomorphic to
$$
\displaylines{
\bigoplus_{\Cal I \subset \Cal F ,\ \Cal I \cap \Cal F _v \neq \emptyset } 
\tilde{H}_{d + |\bar{\Cal I }| - i - 1} (P _{\Cal I }, \Bbb Z ) 
\bigoplus_{\Cal I \subset \Cal F ,\ \Cal I \cap \Cal F _v = \emptyset ,\ \Cal I \neq \emptyset } 
\left ( \tilde{H}_{d + |\bar{\Cal I }| - i - 1} (P _{\Cal I }, \Bbb Z )
\oplus 
\Bbb Z ^{\delta_{i+1} ^{d + |\bar{\Cal I }| }}
\right )\hfill\cr
\hfill\bigoplus_{\Cal I \subset \Cal F ,\ \Cal I \not\supset \Cal F _v } 
\tilde{H}_{d + |\bar{\Cal I }| - i} (P _{\Cal I }, \Bbb Z ) 
\bigoplus_{\Cal I \subset \Cal F ,\ \Cal I \supset \Cal F _v ,\ \Cal I \neq \Cal F } 
\left ( \tilde{H}_{d + |\bar{\Cal I }| - i} (P _{\Cal I }, \Bbb Z ) 
\oplus \Bbb Z ^{\delta_{d-2}^{d + |\bar{\Cal I }| - i}} \right )
}
$$
\noindent and finally to
$$
\displaylines{
\bigoplus_{\Cal I \subset \Cal F ,\ \Cal I \neq \emptyset } 
\tilde{H}_{d + |\bar{\Cal I }| - i - 1} (P _{\Cal I }, \Bbb Z )
\bigoplus_{\Cal I \subset \Cal F ,\ \Cal I \neq \Cal F } 
\tilde{H}_{d + |\bar{\Cal I }| - i} (P _{\Cal I }, \Bbb Z )\hfill\cr
\hfill\bigoplus_{\Cal I \subset \Cal F ,\ \Cal I \cap \Cal F _v = \emptyset ,\ \Cal I \neq \emptyset } 
\Bbb Z ^{\delta_{i - d + 1} ^{|\bar{\Cal I }|}}\bigoplus_{\Cal I \subset \Cal F ,\ \Cal I \supset \Cal F _v ,\ \Cal I \neq \Cal F } 
\Bbb Z ^{\delta_{i-2}^{|\bar{\Cal I }|}}
}
$$
\medskip

  The sum $\bigoplus_{\Cal I \subset \Cal F ,\ \Cal I \neq \emptyset } 
\tilde{H}_{d + |\bar{\Cal I }| - i - 1} (P _{\Cal I }, \Bbb Z )$ is isomorphic to 
$H^i (X, \Bbb Z )$, except if 
$d+n-i-1 = -1$, i.e. $i = d+n$.
\medskip

  Also, the sum $\bigoplus_{\Cal I \subset \Cal F ,\ \Cal I \neq \Cal F } 
\tilde{H}_{d + |\bar{\Cal I }| - i} (P _{\Cal I }, \Bbb Z )$ is isomorphic to 
$H^{i-1} (X, \Bbb Z )$, except if 
$d-i = d-1$, i.e. $i = 1$.
\medskip

  On the other side, 
$$
\sum_{\Cal I \subset \Cal F ,\ \Cal I \cap \Cal F _v = \emptyset ,\ \Cal I \neq \emptyset } 
\delta_{i - d + 1} ^{|\bar{\Cal I }|}
$$ 
is the number of nonempty subsets of 
$\Cal F \backslash \Cal F _v $ having $n-i+d-1$ elements. It is $\left (\smallmatrix n-d \\ n-i+d-1\endsmallmatrix\right )$ except if 
$n-i+d-1 = 0$ i.e. $i = n+d-1$, in which case this sum is zero.
\medskip

  We also have that
$$
\sum_{\Cal F _v \subset \Cal I \subset \Cal F ,\ \Cal I \neq \Cal F } 
\delta_{i - 2} ^{|\bar{\Cal I }|} = 
\sum_{\bar{\Cal I } \subset \Cal F ,\ \bar{\Cal I } \cap \Cal F _v = \emptyset ,\ 
\bar{\Cal I} \neq \emptyset } 
\delta_{i - 2} ^{|\bar{\Cal I }|}
$$ 
is the number of nonempty 
subsets of $\Cal F \backslash \Cal F _v $ having $i-2$ elements. It is $\left ( \smallmatrix n-d \\ i-2 \endsmallmatrix\right )$ except 
if $i-2 = 0$ i.e. $i = 2$, in which case this sum is zero.
\medskip

Putting all these results together and remarking that $(n-d) - (n-i+d-1) = i-2d+1$, we get the isomorphisms of the Proposition.
\medskip

  The proof of the first part of Proposition 11.2 is completed. Let us now 
describe the cup product.

  Rule 1 is obvious.

  To continue, we have to define clearly our sums $S_j $ because they derive from 
isomorphims which are, as we shall see right now, not canonical.

  Look first at the isomorphism $\tilde{H}_0 (P '_{\Cal I _1 }, \Bbb Z ) \simeq 
\tilde{H}_0 (P _{\Cal I }, \Bbb Z ) \oplus \Bbb Z $ where $\Cal I $ is nonempty and does not 
meet $\Cal F _v $. This isomorphism is canonical when (not reduced) homology is concerned, but the cycles 
that are added (multiples of the singleton $\langle v\rangle $) are not cycles in reduced homology. 
Look now at the isomorphism $\tilde{H}_{d-2} (P ' _{\Cal I _2 }, \Bbb Z ) \simeq 
\tilde{H}_{d-2} (P _{\Cal I }, \Bbb Z ) \oplus \Bbb Z $ where $\Cal I \neq \Cal F $ and contains 
$\Cal F _v $. The projection of $\tilde{H}_{d-2} (P ' _{\Cal I _2 }, \Bbb Z )$ over 
$\tilde{H}_{d-2} (P _{\Cal I }, \Bbb Z )$ is canonical (hence is its kernel which is the 
factor $\Bbb Z $), but the inclusion of $\tilde{H}_{d-2} (P _{\Cal I }, \Bbb Z )$ in 
$\tilde{H}_{d-2} (P ' _{\Cal I _2 }, \Bbb Z )$ is not.

  Consider a nonempty subset $\Cal I $ of $\Cal F $ disjoint from $\Cal F _v $. Choose now any 
reduced homology class  in $\tilde{H}_0 (P '_{\Cal I _1 }, \Bbb Z )$ whose value on the 
connected component $F$ of $P '_{\Cal I _1 }$ is equal to $1$ and call $[c_{\Cal I }]$ this 
class. It is clear that the groups $\Bbb Z \cdot [c_{\Cal I }]$ and 
$\tilde{H}_0 (P _{\Cal I }, \Bbb Z )$ whose inclusion in 
$\tilde{H}_0 (P '_{\Cal I _1 }, \Bbb Z )$ results from the inclusion 
$P _{\Cal I } \subset P '_{\Cal I _1 }$ give the desired isomorphism. Doing this for every $\Cal I$, we thus have
$$
S_3=\bigoplus_{\Cal I \subset \Cal F ,\ \Cal I \cap \Cal F _v = \emptyset ,\ \Cal I \neq \emptyset}
\Bbb Z\cdot [c_{\Cal I}]
$$

  Consider now $\bar{\Cal I }$. It is a proper subset of $\Cal F $ which contains $\Cal F _v $.
The linking operation on $\tilde{H}_0 (P '_{\Cal I _1 }, \Bbb Z ) \times 
\tilde{H}_{d-2} (P ' _{\bar{\Cal I} _1 }, \Bbb Z )$ is well defined and the subgroup of the 
homology classes that are not linked with $[c_{\Cal I }]$ is isomorphic to 
$\tilde{H}_{d-2} (P _{\bar{\Cal I} }, \Bbb Z )$. As a consequence, 
$\tilde{H}_{d-2} (P ' _{\bar{\Cal I} _1 }, \Bbb Z )$ is the direct sum of this subgroup with the 
group generated by the class $[c'_{\bar{\Cal I }}]$ of a sphere that ''turns around $F$'' 
(this group is also the kernel of the projection coming from the inclusion 
$P ' _{\bar{\Cal I} _1 } \subset P _{\bar{\Cal I }}$). We thus obtain
$$
S_4=\bigoplus_{\Cal I \subset \Cal F ,\ \Cal I \cap \Cal F _v = \emptyset ,\ \Cal I \neq \emptyset}
\Bbb Z\cdot [c'_{\bar{\Cal I }}]
$$

  Rule 5 is now clear. More precisely, if we take $[c_{\Cal I }]$ and $[c'_{\bar{\Cal J }}]$ as explained above, the 
cup product of the corresponding cohomology classes is zero if $\Cal I \neq \Cal J $. Indeed, if $\Cal I \neq \Cal J $, then 
$\Cal I \cup \bar{\Cal J } \neq \Cal F $ or $\Cal I \cap \bar{\Cal J } \neq \emptyset $. By Cohomology Theorem 7.6, the cup product
is automatically $0$ in the first case; and in the second case, it lies in $\tilde{H}_{-1}(\Cal I \cap \bar{\Cal J },\Bbb Z)$. As
this group is reduced to zero, the cup product is zero too. On the other hand 
the cup product of the classes associated to $[c_{\Cal I }]$ and $[c'_{\bar{\Cal I }}]$ is, up to 
sign, the top class of $X'$ (more precise choices allow to obtain exactly the top class 
every time). This gives rule 5.

  For Rule 2, remark first that if both $[c]$ and $[c']$ are in $S_j $ with $j \neq 1$, 
then the union of the corresponding subsets of $\Cal F \cup \{F\}$ is not all $\Cal F \cup \{F\}$ (indeed $F$ is not in 
this union if $j$ is $2$ or $4$ and $\Cal F _v $ does not meet the union if $j = 3$). We 
then just have to see that $[c] \smile [c']$ vanishes if $j \leq 2$ and $j' \geq 3$.

  Consider first a class $[c'_{\bar{\Cal I}}]$ in $S_4 $. It is realized by a $(d-2)$-sphere 
which surrounds $F$. Remark that every (reduced) homology class in a $P_{\Cal I }$ can be 
realized by a cycle which is far away from $v$ (except if $\Cal I = \Cal F $ but then the 
corresponding class is in $H^0 (X', \Bbb Z )$ and rule 1 applies). As $F$ and thus the sphere 
realizing $[c'_{\bar{\Cal I}}]$ can be thought of very close to $v$, they do not intersect 
(neither are they linked). Hence, if $[c']$ is in $S_4 $ and $[c]$ is in $S_{j'}$ with 
$j' \leq 2$, then $[c] \smile [c'] = 0$.

  Consider now a class $[c_{\Cal I }]$ in $S_3 $. Let $\Cal J\not=\Cal F$ and let $[a_{\Cal J}]$ be a class of $\tilde H_k(P'_{\Cal J_2},\Bbb Z)$. By arguments
similar to those used in the proof of Rule 5, we have that the intersection class $[c_{\Cal I }]\cap [a_{\Cal J}]$ corresponds to a non-trivial
cohomology class of $X'$ if and only if  $[a_{\Cal J}]$ is a multiple of $[c'_{\bar{\Cal J}}]$. But such a class is not in $S_2$ and thus the cup product
of a class of $S_2$ with a class of $S_3$ is always zero.

  Rules 3 and 4 derive from our Theorems 7.6 and 10.1. Assume that $F$ is the 
greatest element for the order we consider on $\Cal F \cup \{ F \} $.

  For a proper nonempty subset $\Cal I $ of $\Cal F $, and an element 
$[a] \in \tilde H_k (P _{\Cal I }, \Bbb Z )$, recall that $\psi ([a])$ is its image in 
$H^{|\bar{\Cal I}| + d - k - 1} (X, \Bbb Z )$. Let
$\psi _i ([a]) = \varphi _i (\psi ([a]))$ for $i=1,2$. Via our isomorphisms, $[a]$ is identified 
to some classes $[a_j] \in \tilde H_k (P' _{\Cal I _j }, \Bbb Z )$ for $j=1,2$. Noting $\psi '$ the 
application on $X'$ which is equivalent to $\psi $ on $X$, we have 
$\psi _j ([a]) = \psi ' ([a_j] )$ for $j=1,2$.

  Consider now $[a] \in \tilde H_k (P _{\Cal I }, \Bbb Z )$ and 
$[b] \in \tilde H_{k'} (P _{\Cal J }, \Bbb Z )$ with $\Cal I$ and $\Cal J$ proper and nonempty. Assume 
moreover that $\Cal I \cup \Cal J = \Cal F $ (else cup products are zero). Remark that 
$[a_1] \cap [b_j] = ([a] \cap [b])_j$ for $j=1,2$. For a finite set $E$ denote by $K'(E)$ the number
$|E|-d+k'-1$. We then compute:
$$
\eqalign{
\psi _1 ([a]) \smile \psi _1 ([b]) =& \psi '([a_1] ) \smile \psi '([b_1] )\cr
=&\epsilon _{\bar {\Cal I} _1 \bar{\Cal J} _1 } (-1)^{(d+1+n+1+K'(\bar{\Cal J}_1) |\bar{\Cal I}_1|)} \psi' ([a_1] \cap [b_1] )}
$$
and then
$$
\eqalign{
\psi _1 ([a]) \smile \psi _1 ([b])=&- \left( 
\epsilon _{\bar{\Cal I }\bar{\Cal J }} (-1)^{(d+1+n+K'(\bar{\Cal J})|\bar{\Cal I}|)} \psi '(([a] \cap [b])_1 ) \right) \cr
=& - \varphi _1 \left( \epsilon _{\bar{\Cal I }\bar{\Cal J }} (-1)^{(d+1+n+K'(\bar{\Cal J})|\bar{\Cal I}|)} \psi_1([a] \cap [b]) 
\right) \cr
=& - \varphi _1 (\psi ([a]) \smile \psi ([b]))\ .
}
$$ 

Rule 3 results from this. We also have 
$$
\eqalign{\psi _1 ([a]) \smile \psi _2 ([b]) =& \psi '([a_1] ) \smile \psi '([b_2] ) \cr
=&
\epsilon _{\bar{\Cal I} _1 (\bar{\Cal J _2 })} (-1)^{(d+1+n+1+K'(\bar{\Cal J}_2)|\bar{\Cal I}_1|)} 
\psi '([a_1] \cap [b_2] ) \cr
=& (-1)^{1 + \bar{\Cal I}} \left( 
\epsilon _{\bar{\Cal I }(\bar{\Cal J } \cup \{ F \})} (-1)^{(d+1+n+K'(\bar{\Cal J}) |\bar{\Cal I}|)} 
\psi '(([a] \cap [b])_2 ) \right) \cr
=& -\left( 
\epsilon _{\bar{\Cal I }\bar{\Cal J }} (-1)^{(d+1+n+K'(\bar{\Cal J}) |\bar{\Cal I}|)} 
\psi '(([a] \cap [b])_2 ) \right)
}
$$
that is
$$
\psi _1 ([a]) \smile \psi _2 ([b])= -\varphi _2 (\psi ([a]) \smile \psi ([b]))\ .
$$ 

Rule 4 results from this. The Proposition is now proved.
$\square$
\enddemo

\example{Example 11.5}
  Consider the cube as simple polytope. By Corollary 4.6, the associated manifold is the product of three 
$3$-spheres (cf Example 7.9). Cut now a vertex. The resulting simple polytope has dimension $3$ and seven 
facets, hence the associated manifold $X$ has dimension $10$. Note also a 
$\frak S _3 $-symmetry. Let us compute its cohomology ring as an application of Proposition 11.2.

  Number $0$ the "cut face", $1$, $2$, $3$ the adjacent faces to $0$ and $1'$, $2'$, $3'$ the 
"opposite" faces to $1$, $2$, $3$ respectively. 

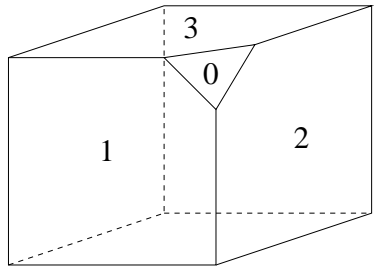
\hfil\scaledpicture 3.5in by 2.5in (fig18 scaled 500) \hfil

The cohomology groups of $X$ are free and 
the Betti numbers are: 
$$
\vbox{
\def\tv{\vrule height 12pt depth 5pt width .5pt}
\offinterlineskip
\halign{\tv \hfill \kern .7em {$#$} \kern .7em \hfill \tv &&\hfill \kern .7em {#}\kern .7em\hfill\tv\cr
\noalign{\hrule}
i &0 ; 10 & 1 ; 9 & 2 ; 8 & 3 ; 7 & 4 ; 6 & 5\cr
\noalign{\hrule}
b_i (X) & 1 & 0 & 0 & 6 & 6 & 2 \cr
\noalign{\hrule}
}}
$$

Denote by $\lambda _i$ for $1 \leq i \leq 3$ the cohomology classes which generate 
$H^3 (\Bbb S^3 \times \Bbb S^3 \times \Bbb S^3 , \Bbb Z )$, and by $\lambda_{ij}$ the cup product $\lambda_i\smile
\lambda_j$. For $l=1,2$ let $\lambda_{i,l}$ (respectively $\lambda_{ij,l}$) be $\varphi_l(\lambda _i)$ (respectively
$\varphi_l (\lambda_{ij})$). The expression $e_{\Cal I}$ for some $\Cal I\subset\{0,1,2,3,1',2',3'\}$ denotes the generator
of a cohomology class of $P_{\Cal I}$ and will be only used when $P_{\Cal I}$ has only one not zero reduced homology group and when
this group is isomorphic to $\Bbb Z$ (e.g. $P_{\Cal I}$ has the homotopy type of a circle). Finally, we denote by $\sigma $ 
a permutation of the set $\{1,2,3\}$. Letting $\sigma$ varies among the permutations of $\{1,2,3\}$, we have:
\smallskip
\item{$\bullet$} $H^3 (X, \Bbb Z )$ is generated by $\lambda _{\sigma (1),1}$ and 
$e_{123\sigma (1)'\sigma(2)')}$;
\smallskip
\item{$\bullet$} $H^4 (X, \Bbb Z )$ is generated by $\lambda _{\sigma (1),2}$ and 
$e_{123\sigma (1)'}$;
\smallskip
\item{$\bullet$} $H^5 (X, \Bbb Z )$ is generated by $e_{123}$ and $e_{01'2'3'}$;
\smallskip
\item{$\bullet$} $H^6 (X, \Bbb Z )$ is generated by $\lambda _{\sigma (1)\sigma (2),1}$ and 
$e_{0\sigma (1)'\sigma (2)')}$;
\smallskip
\item{$\bullet$} $H^7 (X, \Bbb Z )$ is generated by $\lambda _{\sigma (1)\sigma (2),2}$ and 
$e_{0\sigma (1)'}$.
\medskip

  The product of these generators are zero except: 

\noindent i) $\lambda _{\sigma (1),1} \smile \lambda _{\sigma (2),1} = -\lambda _{\sigma (1)\sigma (2),1}$;

\noindent ii) $\lambda _{\sigma (1),1} \smile \lambda _{\sigma (2),2} =- 
\lambda _{\sigma (1) \sigma (2),2}$ and $\lambda _{\sigma (2),1} \smile \lambda _{\sigma (1),2} = 
-\lambda _{\sigma (1)\sigma (2),2}$; 

\noindent iii) The products which give the top class, i.e. 
$-(\lambda _{\sigma (1),1} \smile \lambda _{\sigma (2)\sigma (3),2})$, $e_{\Cal I } \smile e_{\bar{\Cal I }}$ and
$-(\lambda _{\sigma(1)\sigma(2),1} \smile \lambda _{\sigma (3),2})$.
\endexample
\medskip

It is easy to check that, in the previous Example, the cohomology ring of the associated link is isomorphic neither to that of a sphere, 
nor to that of a connected sum of sphere products, nor to that of the product of such manifolds. 
The answer to Questions A and A' is thus negative yet in dimension $3$. Notice that the exact diffeomorphism type of the link
of the previous example is not clear. We may ask

\proclaim\nofrills{Question : \ }
   Describe this manifold more precisely: for instance, can it be decomposed into a connected sum 
of manifolds?
\endproclaim

In dimension $3$, we may in fact characterize precisely which simple polytopes give rise to connected sums of sphere products as links,
and which manifolds appear in this way. We have

\proclaim{Proposition 11.6}
Let $P$ be a simple polyhedron (so $d=3$). Then, the following statements are equivalent:

\noindent (i) The cohomology ring of the associated link $X$ is isomorphic to that of a connected sum of sphere products.

\noindent (ii) The link $X$ is diffeomorphic to a connected sum of sphere products.

\noindent (iii) There exists $l>0$ such that $X$ is diffeomorphic to
$$
\sc _{j=1}^l j \pmatrix l+1 \\ j+1 \endpmatrix \Bbb S^{2+j}\times\Bbb S^{6+l-j-1}\ .
$$

\noindent (iv) There exists $l>0$ such that $P$ is obtained from the $3$-simplex by cutting off $l$ well chosen vertices.
\endproclaim

\demo{Proof}
By application of Theorem 6.3, we know that (iv) implies (iii), and of course (iii) implies (ii) and (ii) implies (i), so it
is sufficient to prove that (i) implies (iv). We assume thus that the cohomology ring of the associated link $X$ is isomorphic to 
that of a connected sum of sphere products.

\definition{Definition 11.7}
  Let $\Cal I $ be a subset of $\Cal F $. We say that $\Cal I $ is a $1$-cycle of facets 
of $P $ if $K_{\Cal I} $ is a cycle (i.e. a connected 
graph all of whose vertices are bivalent).
\enddefinition

  A $1$-cycle of facets can also be viewed as the data of an integer $k \geq 3$ 
and an injective map from $\Bbb Z_k $ into $\Cal I $ such that the images 
of two elements meet if and only if the two elements are equal or consecutive 
in $\Bbb Z_k $, and if moreover the $k$ facets do not have 
a common vertex. The integer $k$ is then called the length of the $1$-cycle of 
facets.

\noindent $\underline{\text{Claim:}}$  consider two disjoint facets $F$ and $F'$ of $P $. Then 
$\Cal F \backslash \{ F , F' \} $ contains a $1$-cycle of facets.

To see this, consider the set $\Cal I _F $ of facets that meet $F$ (except $F$ itself). 
Consider the maps $\phi $ from $\Bbb Z_k $ into $\Cal I _F $ having the 
following properties:

\noindent i) for all $i$ in $\Bbb Z_k $, $\phi (i)$ meets $\phi (i+1)$.

\noindent ii) for all $i$ in $\Bbb Z_k $, consider the segment on $\phi (i)$ 
joining the centers of the edges $\phi (i) \cap \phi (i-1)$ and 
$\phi (i) \cap \phi (i+1)$. We require the polygon obtained by concatenation of 
all these segments to be nontrivial in the homology of 
$P _{\partial } \backslash (F \cup F')$.

There exist such maps: order $\Cal I _F $ such that the bijective order-preserving map from $\Bbb Z_{\vert \Cal I _F \vert}$ to 
$\Cal I _F $ satisfies i). Then this map also satisfies ii), since the polygon obtained from it is homotopic to the boundary of
$F$.
Moreover, let us prove that a minimal subset of $\Cal I _F $ fulfilling these conditions is a 
$1$-cycle of facets.

  First, such a minimal subset cannot contain exactly three globally meeting 
facets, as in this case the polygon considered in the point ii) would be 
contained in a contractible subset (the union of the three faces) of 
$P _{\partial } \backslash (F \cup F')$, which is not allowed.

  Assume now that in this minimal subset $\{C_1 ,..., C_k\} $, the facet $C_1 $ 
 meet $C_j $, for some $j$ such that $2<j<k$. Then $\{C_1 ,..., C_j\} $ and 
$\{C_1 , C_j , C_{j+1} ,..., C_k\} $ satisfy i) and one of them satisfies 
ii), as the polygon of $C_1 ,..., C_k $ is homologically the sum of the 
polygons of these two subsets. Contradiction.

  This completes the proof of the claim.

  We denote by $(*)$ the property, for a simple $3$-dimensional polytope, that 
all its $1$-cycles of facets have length $3$.

  Assume that $P $ does not satisfy $(*)$. Then we can take a $1$-cycle of facets $\Cal I$
of length $k \geq 4$ of $P $. In particular, $I_1 $ and $I_3 $ are disjoint. 
The complement of $P _{\Cal I }$ in $P $ has two connected components $\Cal X $ and 
$\Cal Y $.

  The group $H_1 (P _{\Cal I }, \Bbb Z )$ is isomorphic to $\Bbb Z $, generated 
by the class of the ``polygon'' $T$ whose vertices are the centers of the 
intersections of facets of $\Cal I $.

 Consider now $\Cal J = \{ I_1 ; I_3 \} \cup \left( \Cal F \backslash \Cal I \right) $.
  The group $H_1 (P _{\Cal J }, \Bbb Z )$ is isomorphic to $\Bbb Z $ too, 
generated by the class of a cycle $T '$ which is decomposed as follows: 
for $i = 1$ or $i=3$, let $x_i $ (respectively $y_i $) be in the intersection of $I_i $ with 
$\Cal X $ (respectively $\Cal Y $). Consider a segment in $I_i $ joining $x_i $ to $y_i $ 
and a path in the interior of $\Cal X $ (respectively $\Cal Y $) joining $x_1 $ to $x_3 $ 
(respectively $y_1 $ to $y_3 $). The cycle $T '$ is obtained by the concatenation of 
these four paths.

The next picture represents such a situation. Here $P$ is the cube with the same numbering of facets as in Example 7.9. The 1-cycle of facets is
$\Cal I=\{1,2,1',2'\}$, so $\Cal J=\{1,3,1',3'\}$, whereas $\Cal X=3$ and $\Cal Y=3'$.

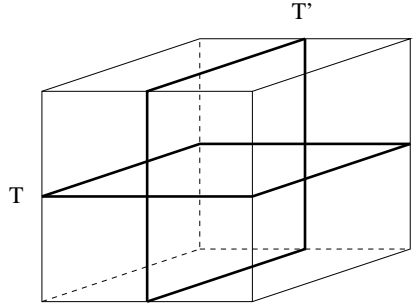
\hfil\scaledpicture 3.9in by 2.9in (fig14bis scaled 500) \hfil

  Now $\Cal I \cup \Cal J = \Cal F $ and $\Cal I \cap \Cal J = \{ I_1 ; I_3 \} $. On $I_3 $ and on $I_1$, 
the intersection of $T$ and $T'$ is exactly one point. In particular, 
the intersection class of these two cycles in $H_0 (I_1 \cup I_3 , \Bbb Z )$ 
cannot be zero. By Cohomology Theorem 7.6, the class $\psi([T])$ (respectively $\psi([T'])$) is non-trivial of dimension $\vert \bar{\Cal I}\vert +1$
(respectively $\vert \bar{\Cal J}\vert +1$). Still by Cohomology Theorem 7.6, the cup product $\psi([T])\smile \psi([T'])$ is a non-trivial cohomology class.

This class does not belong to the top-dimensional cohomology group of $X$, since the top class corresponds to the generator of 
$\tilde H_{-1}(\emptyset,\Bbb Z)$.
This means that the cohomology ring of $X$ is not isomorphic to 
that of a connected sum of sphere products. Contradiction. The polytope $P$ has only $1$-cycles of facets of length $3$.

We now have to show the converse, i.e. if $P $ satisfies $(*)$, then 
$P $ is obtained from the tetrahedron by vertex cutting. Remark that a 
polyhedron which is obtained from the tetrahedron by vertex cutting has (at 
least) two disjoint triangular facets (except if it is the tetrahedron itself).

  Assume that $P $ has a triangular face. Then, if $P $ is not itself the 
tetrahedron, we can perform a flip of type $(3,1)$ along this face so that it 
disappears. The resulting polytope $Q$ satisfies $(*)$ too as we cannot have 
created new $1$-cycles of facets. It has one face less than $P$ and $P$ is obtained from $Q$ by vertex cutting..

  Hence, by induction on the number of facets, we just have to show that a 
polytope having the property $(*)$ has necessarily a triangular face.

  Consider a polytope $P $ fulfilling $(*)$. If $P $ is not a tetrahedron, it 
has two disjoint facets and, according to the claim, a $1$-cycle of facets 
$(F_1 , F_2 , F_3 )$ of length $3$. Now, the plane $H$ passing through the 
centers of the intersections $F_i \cap F_j $ intersects no other facet. The 
intersections $P ^+ $ and $P ^- $ of $P $ with the two half-planes 
delimited by $H$ are simple convex polytopes satisfying $(*)$ and with a triangular face $H\cap P$. If $P ^+ $ is 
$P $ itself, then $P $ has a triangular face. Else $P ^+ $ has strictly 
less faces than $P $ and, by induction, is obtained from the tetrahedron by 
vertex cutting. As it cannot be the tetrahedon (because 
$F_1 \cap F_2 \cap F_3 $ is empty), it has two disjoint triangular facets, and 
in particular one which is disjoint from $H \cap P $. This facet is also a 
triangular facet of $P $, which completes the proof.
$\square$
\enddemo

In higher dimension, the simple polytopes obtained from the simplex (of same dimension) by cutting off vertices still give rise
to links whose cohomology ring is
isomorphic to that of a connected sum of products of spheres by Theorem 6.3. Nevertheless, there are not the only ones and a
nice characterization of all the polytopes having this property seems not to exist. In particular, the results of \cite{LdM2} recalled in Example 0.5
give examples of connected sums of products of spheres which cannot be obtained by Theorem 6.3. We use the notations of Example 0.5.
Set $n=10$ and $n_1=\hdots=n_5=2$. Then, the associated link $X$ is diffeomorphic to $\sc (5) \Bbb S^7 \times \Bbb S^{10}$. Since $X$ is $6$-connected,
it is not diffeomorphic to one of the links obtained by Theorem 6.3: none of them is $3$-connected. Moreover, we may construct other examples. To do that, recall that 
a(n even dimensional) polytope is called neighbourly 
if every subset of cardinal $\frac{d}{2}$ determines a face, and that such a polytope is simplicial (see Section 2 and \cite{Gr}). A 
polytope whose dual is neighbourly is therefore simple and is called a dual neighbourly 
polytope. Here, we will only consider the even dimensional case.

\proclaim{Proposition 11.8}
  Assume that $P$ is dual neighbourly and of even dimension. Then the cohomology ring of $X$ is 
isomorphic to the one of a connected sum of sphere products.
\endproclaim

\demo{Proof}
  We try to compute the reduced homology groups of $P _{\Cal I }$, for $\Cal I $ proper and 
nonempty. Recall that this set is homotopy equivalent to the subcomplex of $P ^* $ 
corresponding to the {\it maximal} subcomplex whose vertices are those related to the facets of 
$\Cal I $. For $k < \frac{d}{2} - 1$, the $k+1$-skeleton of $P _{\Cal I }^* $ is complete by 
definition of neighbourlyness, hence $P _{\Cal I }$ has trivial reduced $k$-(co)homology.

The torsion part of $\tilde H_{\frac{d}{2} - 1} (P _{\Cal I }, \Bbb Z )$ is isomorphic to the torsion part of the 
group $\tilde H^{\frac{d}{2}} (P _{\bar{\Cal I }}, \Bbb Z )$. From Lemma 7.4 and Alexander-Pontrjagin duality (see \cite{Al}, t. 3, p.53), it
is also isomorphic to the torsion part of the 
group $\tilde H_{\frac{d}{2} - 2} (P _{\bar{\Cal I }}, \Bbb Z )$, hence is trivial. 
In the same way, for $k \geq \frac{d}{2}$, the group $\tilde H_k (P _{\Cal I }, \Bbb Z )$ is isomorphic to 
the direct sum of the free part of $\tilde H_{d-k-2} (P _{\bar{\Cal I} }, \Bbb Z )$ and of the torsion 
part of $\tilde H_{d-k-3} (P _{\bar{\Cal I }}, \Bbb Z )$, both being trivial.

  To sum up, the reduced homology groups of $P _{\Cal I }$ vanish except in dimension 
$\frac{d}{2} - 1$ in which case it is free.

  Furthermore, if the homology intersection of two such classes is nonzero, then it must 
lie in the reduced homology group of dimension $-1$ of some subset of $\Cal F $, which must 
be the emptyset. Finally, to conclude, we just have to see that the linking number is a 
unimodular bilinear form on 
$\tilde H_{\frac{d}{2} - 1} (P _{\Cal I }, \Bbb Z ) \times 
\tilde H_{\frac{d}{2} - 1} (P _{\bar{\Cal I }}, \Bbb Z )$, which results from the 
"little Pontrjagin duality" (see \cite{Al}, t. 3, p.91).

  This proves the lemma.
$\square$
\enddemo

\example{Example 11.9}
The (even dimensional) cyclic polytopes (\cite{Gr}, \S 4.7) are examples of neighbourly polytopes. For any $d$ and any $v\geq d+1$, there exists a unique cyclic
polytope $C(d,v)$ of dimension $d$ with $v$ vertices. Let us take $d=4$. Then $C(4,5)$ is the $4$-simplex, while $C(4,6)$ is dual to the product of two triangles.
Using the Dehn-Sommerville equations (\cite{Gr}, Chapter 9), it is easy to check that $C(4,7)$ has $28$ faces of dimension $2$ and that $C(4,8)$ has $40$ such faces. 
Comparing these
numbers with the number of $2$-faces of the $6$-simplex and of the $7$-simplex, this means that, in $C(4,7)$, there exist $7$ subsets $\Cal I$ such that 
$P_{\Cal I}^*$ is not contractible but homotopic to a circle, and, in $C(4,8)$, there exist $16$ such subsets. Using the homology formula of Remark 7.7, Proposition 11.8
and Lemma 0.10, we get easily the following table.
$$
\vbox{
\def\tv{\vrule height 12pt depth 5pt width .5pt}
\offinterlineskip
\halign{\tv \hfill \kern .7em {$#$} \kern .7em \hfill \tv &&\hfill \kern .7em {$#$}\kern .7em\hfill\tv\cr
\noalign{\hrule}
v & 5 & 6 & 7 & 8 \cr
\noalign{\hrule}
X & \Bbb S^9 & \Bbb S^5\times\Bbb S^5& \sc (7) \Bbb S^5\times\Bbb S^6 & \sc (16) \Bbb S^5\times\Bbb S^7 \sc (15) \Bbb S^6\times\Bbb S^6 \cr
\noalign{\hrule}
}}
$$

In the first three cases, the table gives the diffeomorphism type of $X$; in the third case, this follows from the fact that the same
example can be obtained from Example 0.5 (take $n=k=7$ and use Lemma 1.3). On the contrary, it guarantees only the cohomology ring of $X$
in the last case. Notice that this last case can be obtained neither from Theorem 6.3 nor from Example 0.5.
\endexample
\medskip

This leads to the conjecture:

\proclaim{Conjecture}
  If $P $ is dual neighbourly, then $X$ is actually the connected sum of sphere 
products (if not a sphere).
\endproclaim

\remark{Remark 11.10}
One difficult step in proving the conjecture is to prove that, if $P$ is dual neighbourly, then $X$ has the homotopy type of a connected
sum of sphere products. Relating to this is the more general question

\proclaim{Question}
Let $X$ and $X'$ be two links. Assume that they have isomorphic cohomology rings. Are they homotopy equivalent?
\endproclaim

We will go back to this question in Part III.
\endremark
\medskip 

To finish with this part, we have to answer Question A''. Indeed, a link may not only have torsion in (co)homology, but arbitrary
torsion!

\proclaim{Torsion Theorem 11.11}
The (co)homology groups of a $2$-connected link may have arbitrary amount of torsion. More precisely, let $G$ be any abelian finitely 
presented group. Then, there exists a $2$-connected link $X$ such that $H^i(X,\Bbb Z)$ contains $G$
as a free summand (that is $H^i(X,\Bbb Z)=G\oplus \hdots$) for some $2<i<\dim X-2$.
\endproclaim

This is a very surprising result (at least for the authors) since the links are transverse intersections of quadrics with very
special properties ...

\demo{Proof}
Let $G$ be an abelian finitely presented group. Let $K$ be a finite simplicial complex such that $\tilde H_i(K,\Bbb Z)=G$ for some $i>0$.
Let $\{1,\hdots l\}$ be the vertex set of $K$.
Consider the $(l-1)$-simplex and let its set of facets be $\{1,\hdots l\}$. 
For every simplex $I=(i_1,\hdots,i_p)$ of maximal dimension 
of $K$, cut off the face of the $(l-1)$-simplex numbered $\{1,\hdots l\}\setminus I$
by a generic hyperplane. We thus obtain a simple convex polytope $P$. Notice that its number of facets $n$ is the sum of $l$ with
the number $f$ of facets of $K$. Set $\Cal F=\{1,\hdots, l,l+1,\hdots l+f\}$. Finally, consider the associated link $X$. 

The crucial remark is that $link_{\Delta}\sigma_{\{l+1,\hdots,l+f\}}$ is isomorphic to $K$. Indeed, by Remark 7.11, we have
$$
link_{\Delta}\sigma_{\{l+1,\hdots,l+f\}}=\{I\subset \{1,\hdots,l\} \quad\vert\quad F_{\hat I}=\emptyset\}\ .
$$

Now, $F_{\hat I}$ is empty if and only if $\hat I$ numbers a face of the $(l-1)$-simplex which is cut off when passing to $P$, i.e. if and only if $I$ numbers
a simplex of $K$. By 
application of Cohomology Theorem 7.6 and Lemma 7.12, every (reduced) homology group of $K$ will thus appear as a free summand of some 
cohomology group of $X$. So will appear $G$.
$\square$
\enddemo

\remark{Remark 11.12}
The proof of this Theorem is perhaps easier to understand when modified as follows. Starting with a finite simplicial complex $K$ with $l$ vertices, embed it
as a simplicial subcomplex of the $(l-1)$-simplex $\Delta$. Perform a barycentric subdivision of each face of $\Delta\setminus K$. We thus obtain a simplicial
polytope $P^*$ such that $K$ is the {\it maximal} simplicial subcomplex of $P^*$ of vertex set $\{1,\hdots,l\}$. We conclude with Remark 7.7.
\endremark
\medskip

The proof of Torsion Theorem 11.11 is constructive. Here is an example.

\example{Example 11.13 (compare with [Je])}
Consider the minimal triangulation of the projective plane $\Bbb P^2(\Bbb R)$ drawn at the bottom of this example.
The simplices of maximal dimensions are
$$
\{(356), (456),(246),(235),(145),(125),(134),(234),(126),(136)\}
$$

Consider the $5$-simplex and number its facets $\{1,\hdots, 6\}$. Cut off the faces of this simplex numbered
$$
\{(123),(124),(135),(146),(156),(236),(245),(256),(345),(346)\}
$$
by generic hyperplanes. We thus obtain a simple $5$-polytope with $16$ facets giving rise to a $2$-connected link $X$ of dimension $21$.
Set $\Cal F=\{1,\hdots, 16\}$. 
The complex $link_{\Delta}\sigma_{\{7,\hdots,16\}}$ is homotopic
to the projective plane. By Lemma 7.12, this means that $\tilde H_1(P_{\{7,\hdots,16\}},\Bbb Z)$ is isomorphic to $\Bbb Z_2$. Cohomology Theorem 7.6 implies that
$$
\eqalign{
H^9(X,\Bbb Z)\simeq &\bigoplus_{\Cal I\subset \{1,\hdots, 16\}}\tilde H_{|\bar{\Cal I}|-5}(P_{\Cal I},\Bbb Z)\cr
\simeq &\tilde H_1(P_{\{7,\hdots,16\}},\Bbb Z)\oplus \hdots\simeq\tilde H_1(\Bbb P^2(\Bbb R),\Bbb Z)\oplus \hdots
\simeq \Bbb Z_2\oplus \hdots
}
$$

Therefore, not all the homology groups of $X$ are free.
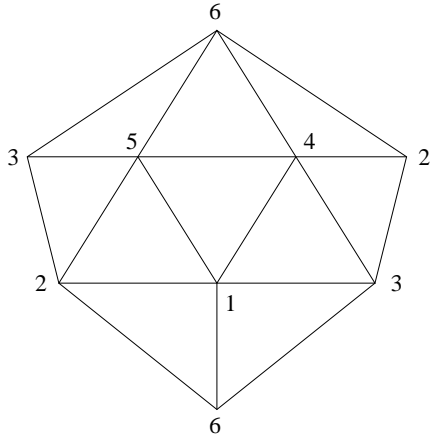
\endexample

\hfil\scaledpicture 6.8in by 6.9in (fig19 scaled 300) \hfil

Notice that, due to Corollary 11.1, the dimension of this counterexample is sharp.

\vfill
\eject

\specialhead
{\nom Part III: Applications to compact complex manifolds}
\endspecialhead
\bigskip

\head
{\bf 12. LV-M manifolds and links}
\endhead

We recall very 
briefly the construction of the LV-M manifold (see \cite{Me1} and \cite{Me2} for more details; this is a generalization
of the construction presented in \cite{LdM-Ve}). Let $m>0$ and $n>2m$ be two
integers. Let $\Lambda=(\Lambda_1,\hdots,\Lambda_n)$ be a set of $n$ vectors of $\C^m$ satisfying the Siegel and the weak hyperbolicity
condition (as vectors of $\R^{2m}$, see Lemma 0.3). Consider the holomorphic foliation $\Cal F$ of the projective space $\Bbb P^{n-1}$
given by the following action
$$
(T,[z])\in\C^m\times\Bbb P^{n-1}\longmapsto [\exp\langle \Lambda_1,T\rangle\cdot z_1,\hdots,\exp\langle \Lambda_n,T\rangle\cdot z_n]\in\Bbb P^{n-1}
\tag\numerote
$$
where the brackets denote the homogeneous coordinates in $\Bbb P^{n-1}$ and where $\langle -,- \rangle$ is the {\it inner} product of
$\C^n$.
Define 
$$
V=\{[z]\in\Bbb P^{n-1}\quad\vert\quad 0\in\Cal H((\Lambda_i)_{i\in I_z})\}
\tag\numerote
$$
where $I_z$ was defined in \lastnum[-20]. We notice that the set $I_z$ is independent of the choice of a representant $z$ of the class
$[z]$. Finally define
$$
\Cal N_{\Lambda}=\{[z]\in\Bbb P^{n-1}\quad\vert\quad \sum_{i=1}^n\Lambda_i\vert z_i\vert ^2=0\}
\tag\numerote
$$
which is a smooth manifold due to the weak hyperbolicity condition (see Lemma 0.3).

\noindent Then it is proven in \cite{Me1} (see also \cite{Me2}) that

\noindent (i) The restriction of $\Cal F$ to $V$ is a regular foliation of dimension $m$.

\noindent (ii) The compact 
smooth submanifold $\Cal N_{\Lambda}$ is a global transverse to $\Cal F$ restricted to $V$, that is cuts every leaf
transversally in an unique point.
\medskip
Therefore, $\Cal N_{\Lambda}$ can be identified with the quotient space of $\Cal F$ restricted to $V$ and thus inherits a complex
structure. We will denote $N_{\Lambda}$ the compact complex manifold obtained in this way. A complex manifold $N_{\Lambda}$ for some
$\Lambda$ will be called a {\it LV-M manifold}. Notice that it has (complex) dimension $n-m-1$. 

The main complex properties of these manifolds are investigated in \cite{Me1}, where\-as a particularly nice connection with projective
toric varieties is explained in \cite{Me2}. We will not need these results, but we will use the following Lemma. Recall that
$\Lambda_i$ is an {\it indispensable point} if $0$ is not in the convex hull of $(\Lambda_j)_{j\not = i}$. 

\proclaim{Lemma 12.1}
Let $N_{\Lambda}$ be a LV-M manifold. Assume that $\Lambda$ has at least $m+1$ indispensable points. 
Then the complex structure of $N_{\Lambda}$ is affine (and even linear), that is may be defined by a holomorphic 
atlas
such that the changes of charts are affine (and even linear) automorphisms
of $\C^{n-m-1}$.
\endproclaim

\demo{Proof}
Assume that $\Lambda_1,\hdots,\Lambda_{m+1}$ are indispensable. By \lastnum[-1], this implies
$$
[z]\in V \Rightarrow z_1\cdot \hdots \cdot z_{m+1}\not = 0
$$
By construction of $N_{\Lambda}$, we just need to construct a foliated atlas of $(V,\Cal F)$ with linear transverse changes of charts.
Look at the map
$$
\displaylines{
(T,w)\in \C^m\times \C^{n-m-1}\buildrel \Phi_z\over\longmapsto [z_1\cdot\exp \langle \Lambda_1, T\rangle, \hdots,
z_{m+1}\cdot\exp \langle \Lambda_{m+1}, T\rangle, \hfill \cr \hfill
w_1\cdot\exp \langle \Lambda_{m+2}, T\rangle, \hdots,
w_{n-m-1}\cdot\exp \langle \Lambda_n, T\rangle] \in V }
$$
for a fixed set $z=(z_1,\hdots, z_{m+1})\in (\Bbb C^*)^{m+1}$. Using the weak hyperbolicity condition, 
it can be shown that the set $(\Lambda_2-\Lambda_1,
\hdots, \Lambda_{m+1}-\Lambda_1)$ has rank $m$. As a consequence, 
$\Phi_z(T,w)=\Phi_{z'}(T',w')$ if and only if 
$$
w'_i=w_i\cdot\exp{\langle \Lambda_{m+1+i}, T-T'\rangle}
\eqno 1\leq i\leq n-m-1
$$ 
and $T-T'$ belongs to a fixed lattice in $\C^m$. 
Therefore $\Phi_z$ is a local homeomorphism and can be used as a local foliated chart for
every point $(z_1,\hdots, z_{m+1}, w)$.
Since the $(m+1)$ first homogeneous coordinates of every point of $V$ are not zero, $V$ can be covered by such charts.
Moreover, the previous computation proves that the changes of charts are uniquely determined by translations along a lattice 
$T\mapsto T+a$ so that the transverse changes of
charts have the form
$$
w\in\Bbb C^{n-m-1} \longmapsto (w_1\cdot\exp\langle \Lambda_{m+2}, a\rangle,\hdots,
 w_{n-m-1}\cdot\exp\langle \Lambda_{n}, a\rangle)
$$
that is are linear.
$\square$
\enddemo

To avoid particular cases in the sequel, we add the special case $m=0$: then there is no action at all and $N$ is by 
definition the projective space
$\Bbb P^{n-1}$.

Let $A\in\Cal A$. The quotient space of $X_A$ by the diagonal action \lastnum[-16] can be identified with
$$
\widetilde X_A=\{[z]\in\Bbb P^{n-1}\quad\vert\quad \sum_{i=1}^n A_i\vert z_i\vert ^2=0\}
\tag\numerote
$$
which is a smooth manifold by Lemma 0.3. In particular, if $X_A$ is not simply-connected, then by Lemma 0.9, it is equivariantly 
diffeomorphic to $X_B\times \Bbb S^1$ for some $B \in \Cal A$. It is then easy to check that $X_B$ and $\widetilde X_A$ are 
equivariantly diffeomorphic.
On the contrary, when $A\in\Cal A_0$, the manifold $\widetilde X_A$ is not a link: for example, 
think about the case where $X_A$ is diffeomorphic to $\Bbb S^3\times\Bbb S^3$ (Example 0.4).

The following Theorem is the motivation for the previous study of the links.

\proclaim{Theorem 12.2}
Let $A\in\Cal A$ of dimensions $p$ and $n$. Then,

\noindent (i) If $p$ is odd, that is if $X_A$ is even-dimensional, then $X_A$ admits a complex structure as a LV-M manifold.

\noindent (ii) If $p$ is even, that is if $X_A$ is odd-dimensional, then $\widetilde X_A$ and $X_A\times\Bbb S^1$ admit a complex
structure as a LV-M manifold.
\endproclaim

\demo{Proof}
Assume that $X_A$ is odd-dimensional, that is that $p$ is even. Setting $m=p/2$ and letting $\Lambda$ denote the image of
$A$ via the standard identification between $\C^m$ and $\R^{2m}$, then $\widetilde X_A$ and $\Cal N_{\Lambda}$ are the same. Therefore,
$\widetilde X_A$ inherits a complex structure. 

If $p$ is odd, define the following matrix
with $n+1$ columns and $p+1$ rows
$$
B=\pmatrix &A &0 \\
&1\hdots 1 &-1
\endpmatrix
$$
This is obviously an admissible configuration and by Lemma 0.9, the links $X_B$ and $X_A\times\Bbb S^1$ are equivariantly
diffeomorphic. As noticed before, this means that 
$\widetilde X_B$ is diffeomorphic to $X_A$ and we are in the previous case.

Finally, if $p$ is even, consider the following matrix with dimensions $n+2$ and $p+2$
$$
C=\pmatrix &A &0 &0 \\
&1 \hdots 1 &-1 &0 \\
&1\hdots 1 &0 &-1
\endpmatrix
$$
Then $X_C$ is equivariantly diffeomorphic to $X_A\times\Bbb S^1\times\Bbb S^1$, and $\widetilde X_C\ed X_A\times\Bbb S^1$ has a complex
structure as a LV-M manifold by what preceeds.  
$\square$
\enddemo

\proclaim{Corollary 12.3}
The product of two links admits a complex structure as a LV-M manifold as soon as it has even dimension.
\endproclaim

\demo{Proof}
Use Example 0.6 and Theorem 12.2, (i).
$\square$
\enddemo

\remark{Remark 12.4}
Let $A\in\Cal A$ and let $A'\in\Cal A$ be obtained from $A$ by a homotopy which does not break the weak hyperbolicity condition.
Then, by Corollary 4.5, the links $X_A$ and $X_{A'}$ are equivariantly diffeomorphic. Nevertheless, the complex structures of
$X_A$ and $X_{A'}$ (if $p$ is odd) or of $\widetilde X_A$ and $\widetilde X_{A'}$ (if $p$ is even) given by Theorem 12.2 are in
general not the same; in this way a link $X_A$ or its diagonal quotient $\widetilde X_A$ comes equipped not only with a complex
structure but with a deformation space of complex structures (see \cite{Me1} where this space is studied).
\endremark

\head
{\bf 13. Holomorphic wall-crossing}
\endhead

Let $N_{\Lambda}$ be a LV-M manifold. Identifying $\R^{2m}$ to $\C^m$ and $\Lambda$ to an element of $\Cal A$, we may talk of a wall $W$ 
of $\Lambda$ (see Definition 5.2) and of a configuration $\Lambda'$ obtained from $\Lambda$ by crossing the wall $W$ (Definition 5.3).
Up to equivariant diffeomorphism, $N_{\Lambda'}$ is obtained from $N_\Lambda$ by performing an equivariant smooth surgery described
in Wall-crossing Theorem 5.4. Nevertheless, $N_\Lambda$ and $N_{\Lambda'}$ being complex manifolds, it is natural to ask which
{\it holomorphic} transformation occurs when performing the wall-crossing. This is what we call the {\it holomorphic wall-crossing
problem}.

\remark{Remark 13.1}
Let $B\in\C^m$ such that $\Lambda'=\Lambda+B$, that is $\Lambda'=(\Lambda_1+B,\hdots,\Lambda_n+B)$. By Definition 5.3, the
configuration $\Lambda+tB$ is admissible for every $t\in [0,1]$, except for one special value $t_0$. 
It follows from \lastnum[-3] and from Corollary 4.5 that $N_{\Lambda}$ and $N_{\Lambda+tB}$ 
are biholomorphic for every $0\leq t<t_0$ and that $N_{\Lambda'}$ and $N_{\Lambda+tB}$ are biholomorphic for every
$t_0<t\leq 1$ (compare with the general case of Remark 12.4). Therefore, the complex structures of the induced links are fixed
before and after crossing the wall.
\endremark
\medskip

In this Section, we will give a complete solution to the holomorphic wall-crossing problem by showing that, in this case, 
the smooth equivariant surgeries occuring during the wall-crossing are in fact holomorphic surgeries. Let us first recall

\definition{Definition 13.2 (see [M-K], p.15)}
Let $M$ be a complex manifold and let $S$ be a holomorphic submanifold of $M$. Let $W$ be a neighborhood of $S$. Finally let $S^*
\subset W^*$ be a pair (holomorphic submanifold, complex manifold) such that $W^*$ is a neighborhood of $S^*$. Given a biholomorphism
$f\ :\ W\setminus S\to W^*\setminus S^*$, we may construct the well-defined complex manifold $M^*$ by cutting $S$ and pasting
$S^*$ by use of $f$. We say that $M^*$ is obtained from $M$ by a {\it holomorphic surgery} along $(S,W,S^*,W^*,f)$.
\enddefinition

Notice that if $f'$ is smoothly isotopic to $f$, the result of performing a holomorphic surgery along $(S,f')$ is diffeomorphic but
in general not biholomorphic to $M^*$. 

\proclaim{Holomorphic wall-crossing Theorem 13.3}
Let $N_{\Lambda}$ be a LV-M manifold. Let $N_{\Lambda'}$ be a LV-M manifold obtained from $N_\Lambda$ by crossing a wall. Then $N_{\Lambda'}$
is obtained from $N_\Lambda$ by a holomorphic surgery. 
\endproclaim

\demo{Proof}
Let $X_F$ be the smooth submanifold of $N_\Lambda$ along which the elementary surgery occurs. 
Using Section 1 and the standard identification of $\R^{2m}$ and $\C^m$, we
have that $X_F$ is the quotient space of the foliation $\Cal F$ restricted to 
$$
V\cap\{z_i=0\quad\vert\quad i\in I\}
$$
for the subset $I\subset \{1,\hdots n\}$ numbering $X_F$ (see \lastnum[-12]). 
Therefore it is a holomorphic submanifold of $N_\Lambda$ corresponding to the
admissible subconfiguration $(\Lambda_i)_{i\in I^c}$. By abuse of notations, we still call $X_F$ this complex manifold. 
On the other hand, we have $V'=V$ and the
submanifold $X'_{F'}$ is the quotient space of $\Cal F'$ restricted to the same $V\cap\{z_i=0\quad\vert\quad i\in I\}$. Define
$W=V\setminus \{z_i=0\quad\vert\quad i\in I\}$. As $\Lambda$ and $\Lambda'$ differ only by a translation factor, the open
complex manifolds $W/\Cal F=N_{\Lambda}\setminus X_F$ and $W/\Cal F'=N_{\Lambda'}\setminus X'_{F'}$ are biholomorphic. 
More precisely, the identity map of $W$ descends to a biholomorphism $f$ between
these two complex manifolds. As a consequence, $N_{\Lambda'}$ is obtained from $N_{\Lambda}$ by a holomorphic surgery along
$(X_F, N_{\Lambda}, X'_{F'}, N_{\Lambda'},f)$.
$\square$
\enddemo

\remark{Remark 13.4}
The holomorphic surgery described in the proof of Theorem 13.3 is a very particular case of Definition 13.2, since the neighborhood $W$ of the submanifold $X_F$
is in fact the whole manifold $N_{\Lambda}$. It is thus a global holomorphic transformation, whereas Definition 13.2 has a local flavour. It is perhaps better to say
that $N_\Lambda$ and $N_{\Lambda'}$ are holomorphic compactifications of the same open complex manifold $N_{\Lambda}\setminus X_F=N_{\Lambda'}\setminus X'_{F'}$.
\endremark

\head
{\bf 14. Topology of LV-M manifolds}
\endhead

As an application of Torsion Theorem 11.11, we have

\proclaim{Theorem 14.1}
The (co)homology groups of a $2$-connected LV-M manifold 
may have arbitrary amount of torsion. More precisely, let $G$ be any abelian finitely 
presented group. Then, there exists a $2$-connected LV-M manifold $N_{\Lambda}$ such that $H^i(N_{\Lambda},\Bbb Z)$ contains $G$
as a free summand (that is $H^i(N_{\Lambda},\Bbb Z)=G\oplus \hdots$) for some $2<i<2n-2m-4$.
\endproclaim

\demo{Proof}
Apply Torsion Theorem 11.11 to obtain a $2$-connected link $X$ with this property. If $X$ is even-dimensional, then we may conclude by 
Theorem 12.2. Otherwise, we perform a surgery of type $(1,n)$ on $X\times\Bbb S^1$. By Proposition 11.2, the resulting $2$-connected
link $X'$ still has the property that $G$ is a free summand of one of its cohomology groups. But now $X'$ is even-dimensional and
we conclude by Theorem 12.2.
$\square$
\enddemo

\remark{Remark 14.2}
As a consequence of a result of \cite{Ta}, every finitely presented group may appear as the fundamental group of a compact complex
non-k\"ahlerian $3$-fold. The previous Theorem is a sort of (much) weaker version of this result for higher dimensional homology groups.
Notice that it is not known if a similar statement is true for K\"ahler manifolds.
\endremark
\medskip
 
Before drawing an interesting consequence of this Theorem, we want to go back to the question asked in Remark 11.10. The ``holomorphic''
version of this question is

\proclaim{Question}
Let $N$ and $N'$ be two LV-M manifolds. Assume that they have isomorphic cohomology rings. Are they homotopically equivalent?
\endproclaim

In the case of two {\it K\"ahler} manifolds, the answer to this question is yes: two K\"ahler manifolds with isomorphic cohomology
rings are indeed homotopically equivalent (see \cite{D-G-M-S}). For non-K\"ahler manifolds, the answer is not in general. Counterexamples
exist yet in dimension two. Consider the open manifold
$$
W=\{(w_1,w_2,w_3)\in\C^3\setminus\{(0,0,0)\}\qquad\vert\qquad w_1^2+w_2^3+w_3^5=0\}
$$
The quotient space of $W$ by the group generated by a well-chosen weighted homothety is a compact complex surface which is diffeomorphic
to $\Sigma\times S^1$, where $\Sigma$ is the Poincar\'e sphere (see \cite{B-VdV} and \cite{Mi}). 
Thinking about the Hopf surfaces, this means that both $S^3\times S^1$ and $\Sigma\times S^1$ admit complex structures. Now they
have isomorphic cohomology rings but different homotopy type (since the Poincar\'e sphere is not simply-connected).

It seems plausible that the techniques of \cite{D-G-M-S} can be applied to the non-K\"ahler class of LV-M manifolds and would bring a
positive answer to the question.
\medskip
Going back to Theorem 14.1, we obtain easily the following surprising Corollary:

\proclaim{Corollary 14.3}
The (co)homology groups of a $2$-connected compact complex {\bf affi\-ne} manifold may have arbitrary amount of torsion (in the sense
of Theorem 14.1).
\endproclaim

\demo{Proof}
By use of Theorem 14.1 and Lemma 12.1, it is enough to prove that, given a LV-M manifold $N_{\Lambda}$ of dimensions $(m,n)$, there
exists a LV-M manifold $N_{\Lambda'}$ of dimensions $(m',n')$ such that

\noindent (i) The manifold $N_{\Lambda'}$ is diffeomorphic to a product of $N_{\Lambda}$ by circles.

\noindent (ii) The number of indispensable points of $N_{\Lambda'}$ is $m'+1$.

Let $\Lambda_l$ be the matrix with $n+2l$ rows
$$
\pmatrix
\Lambda_1 &\hdots &\Lambda_n &0 &&\hdots &&0 \\
-1-i &\hdots &-1-i &1 &i &\hdots &0 &0 \\
\vdots &&\vdots &&&\ddots &&\\
-1-i &\hdots &-1-i &&&&1&i
\endpmatrix
$$
It is straightforward to check that $\Lambda_l$ is admissible, that it has $2l$ indispensable points, and 
that $N_{\Lambda_l}$ is diffeomorphic to $N_\Lambda\times (\Bbb S^1)^{2l}$ (see Example 0.6). The equality $m'+1=2l$ is achieved for $l=m+1$.
$\square$
\enddemo

This means that it is not possible to classify affine complex manifolds or complex manifolds having a holomorphic affine
connection up to diffeomorphism. Notice that an affine compact {\it K\"ahler} manifold is covered by a compact complex torus (see \cite{K-W}).

The previous proof suggests to ask the following question.

\proclaim{Question}
Let $M$ be a compact complex manifold. Under which assumptions on $M$ does the smooth manifold $M\times (\Bbb S^1)^{2N}$ admit a complex affine structure for $N$ sufficiently large?
Is it enough to assume that the total Stiefel-Whitney class and the total Pontrjagin class of $M$ are equal to one?  
\endproclaim

We emphasize that the searched complex affine structure on $M\times (\Bbb S^1)^{2N}$ does not need to respect $M$, that is we {do not}
require that $M$ may be embedded as a holomorphic submanifold of $M\times (\Bbb S^1)^{2N}$ endowed with its affine complex
structure. 

Every compact Riemann surface satisfies the conditions of the second part of the question. Since only the elliptic curves admit affine complex structures, the question is interesting and 
non-trivial even in dimension one. Every compact complex surface which is spin and has signature zero satisfies the conditions of
the second part of the question. Other examples are given by complex manifolds with stably parallelizable smooth tangent bundle (i.e. such that the Whitney sum of 
the smooth tangent bundle with a trivial bundle of sufficiently large rank is trivial). Indeed, this is exactly the case for a link $X_A$, since it is smoothly embedded
in $\Bbb C^n$ with trivial normal bundle, so that
$$
TX_A\oplus E^{p+1}=T\Bbb R^{2n}
$$
where $TM$ denotes the tangent bundle of a smooth manifold $M$ and where $E^k$ denotes the trivial bundle over $X_A$ with fibre $\Bbb R^k$. 

Notice that the condition on the characteristic classes is necessary. For, if $M\times (\Bbb S^1)^{2N}$ admits a complex affine structure,
then the total Chern class of this structure is one (see \cite{K-W}), which implies the same property for the total Stiefel-Whitney and Pontrjagin classes of 
$M\times (\Bbb S^1)^{2N}$. But these classes coincide with the total Stiefel-Whitney and Pontrjagin classes of 
$M$. In particular, for any $n>1$ and for any $N\geq 0$, the smooth manifold $\Bbb P^n\times (\Bbb S^1)^{2N}$ does not admit any complex affine structure by computation of
its Pontrjagin total class (see \cite{M-S}, Example 15.6).

\vfill
\eject

\enddocument